\documentclass[11pt,reqno]{amsart}
\usepackage{amsmath}
\usepackage{cases}
\usepackage{mathrsfs}
\usepackage{amssymb}
\usepackage{amscd}
\usepackage{latexsym,amsthm,amsxtra}
\usepackage[all,cmtip]{xy}
\usepackage{multicol}
\usepackage{hyperref}
\usepackage{mathtools}
\usepackage[margin=1in]{geometry}

\hypersetup{hidelinks}
\usepackage{cite}
\usepackage{xcolor}
\definecolor{RSrevisionyellow}{RGB}{255,247,157}

\usepackage{zref-savepos,zref-abspage,eso-pic}
\makeatletter
\zref@addprop{savepos}{abspage}
\newcommand*\RSmarkstart[1]{\ifvmode\leavevmode\fi\zsavepos{rs-start-#1}}

\newcommand*\RSmarkpage[1]{\zref@extractdefault{#1}{abspage}{0}}
\makeatother
\ExplSyntaxOn
\clist_new:N \g_rs_revision_list_clist
\newcommand*\RSregister[3]{\clist_gput_right:Nn \g_rs_revision_list_clist {{#1}{#2}{#3}}}
\dim_new:N \l_rs_left_dim
\dim_new:N \l_rs_right_dim
\dim_new:N \l_rs_top_dim
\dim_new:N \l_rs_bottom_dim
\dim_new:N \l_rs_startx_dim
\dim_new:N \l_rs_endx_dim
\dim_new:N \l_rs_starty_dim
\dim_new:N \l_rs_endy_dim
\cs_new_protected:Npn \rs_rectangle:nnnn #1#2#3#4
  {
    \dim_compare:nNnT {#3} > {#1}
      {\dim_compare:nNnT {#4} > {#2}
        {\put(\dim_eval:n{#1},\dim_eval:n{#2})
           {\color{RSrevisionyellow}\rule{\dim_eval:n{(#3)-(#1)}}{\dim_eval:n{(#4)-(#2)}}}}}
  }
\cs_new_protected:Npn \rs_span:nnn #1#2#3
  {
    \int_set:Nn \l_tmpa_int {\RSmarkpage{rs-start-#1}}
    \int_set:Nn \l_tmpb_int {\RSmarkpage{rs-end-#1}}
    \int_compare:nT {0 < \l_tmpa_int <= \ReadonlyShipoutCounter <= \l_tmpb_int}
      {
        \dim_set:Nn \l_rs_startx_dim {\zposx{rs-start-#1}sp}
        \dim_set:Nn \l_rs_endx_dim {\zposx{rs-end-#1}sp}
        \dim_set:Nn \l_rs_starty_dim {\zposy{rs-start-#1}sp}
        \dim_set:Nn \l_rs_endy_dim {\zposy{rs-end-#1}sp}
        \int_compare:nT {\l_tmpa_int < \ReadonlyShipoutCounter}
          {\dim_set_eq:NN \l_rs_startx_dim \l_rs_left_dim
           \dim_set:Nn \l_rs_starty_dim {\l_rs_top_dim-#2}}
        \int_compare:nT {\ReadonlyShipoutCounter < \l_tmpb_int}
          {\dim_set_eq:NN \l_rs_endx_dim \l_rs_right_dim
           \dim_set:Nn \l_rs_endy_dim {\l_rs_bottom_dim+#3}}
        \dim_compare:nNnTF {\dim_abs:n{\l_rs_starty_dim-\l_rs_endy_dim}} < {2pt}
          {\rs_rectangle:nnnn {\l_rs_startx_dim}{\l_rs_endy_dim-#3}
                             {\l_rs_endx_dim}{\l_rs_starty_dim+#2}}
          {\rs_rectangle:nnnn {\l_rs_startx_dim}{\l_rs_starty_dim-3pt}
                             {\l_rs_right_dim}{\l_rs_starty_dim+#2}
           \rs_rectangle:nnnn {\l_rs_left_dim}{\l_rs_endy_dim+9pt}
                             {\l_rs_right_dim}{\l_rs_starty_dim-3pt}
           \rs_rectangle:nnnn {\l_rs_left_dim}{\l_rs_endy_dim-#3}
                             {\l_rs_endx_dim}{\l_rs_endy_dim+9pt}}
      }
  }
\AddToShipoutPictureBG{\AtPageLowerLeft{
  \dim_set:Nn \l_rs_left_dim {1in+\hoffset+\oddsidemargin}
  \dim_set:Nn \l_rs_right_dim {\l_rs_left_dim+\textwidth}
  \dim_set:Nn \l_rs_top_dim {\paperheight-1in-\voffset-\topmargin-\headheight-\headsep}
  \dim_set:Nn \l_rs_bottom_dim {\l_rs_top_dim-\textheight}
  \clist_map_inline:Nn \g_rs_revision_list_clist {\rs_span:nnn #1}
}}
\ExplSyntaxOff
\RSregister{weight1}{9pt}{3pt}
\RSregister{weight4}{5.8pt}{2pt}
\RSregister{weight2}{9pt}{3pt}
\RSregister{weight9}{5.8pt}{2pt}
\RSregister{weight11}{9pt}{3pt}
\RSregister{weight3}{9pt}{3pt}
\RSregister{weight5}{5.8pt}{2pt}
\RSregister{weight8}{5.8pt}{2pt}
\RSregister{weight6}{5.8pt}{2pt}
\RSregister{audit-v13-xiao-subfamily}{9pt}{3pt}
\RSregister{audit-v13-xiao-complement}{9pt}{3pt}
\RSregister{weight7}{5.8pt}{2pt}
\RSregister{weight10}{5.8pt}{2pt}
\RSregister{weight12}{9pt}{3pt}
\RSregister{weight13}{9pt}{3pt}
\RSregister{v15-cor-sum-prose}{9pt}{3pt}
\RSregister{weight14}{9pt}{3pt}
\RSregister{v15-cor-twist-sum}{9pt}{3pt}
\RSregister{weight15}{9pt}{3pt}
\RSregister{weight16}{9pt}{3pt}
\RSregister{weight24}{9pt}{3pt}
\RSregister{boost-plainendpoint}{9pt}{3pt}
\hypersetup{
    colorlinks,
    linkcolor={red!50!black},
    citecolor={blue!100!black},
    urlcolor={blue!100!black},
    bookmarksdepth=3,
    bookmarksopen=true,
    bookmarksopenlevel=3,
    pdftitle={Second moments and non-vanishing of Rankin\textendash Selberg L-functions },
    pdfauthor={Hui Wang and Xin Wang}
}
\newcommand{\bea}{\begin{eqnarray}}
\newcommand{\eea}{\end{eqnarray}}
\newcommand{\bna}{\begin{eqnarray*}}
\newcommand{\ena}{\end{eqnarray*}}

\numberwithin{equation}{section}

\theoremstyle{plain}
\newtheorem{lemma}{Lemma}[section]
\newtheorem{theorem}[lemma]{Theorem}
\newtheorem{corollary}[lemma]{Corollary}
\newtheorem{proposition}[lemma]{Proposition}

\theoremstyle{definition}

\newtheorem{remark}{Remark}

\renewcommand{\Re}{\operatorname{Re}}
\renewcommand{\Im}{\operatorname{Im}}

\newcommand{\GL}{\operatorname{GL}}

\newcommand{\sym}{\operatorname{sym}}

\newcommand{\Res}{\mathop{\operatorname{Res}}}

\newcommand{\dd}{\mathrm{d}}

\newcommand{\eps}{\varepsilon}

\title[Second moments and non-vanishing of Rankin--Selberg $L$-functions]{Second moments and non-vanishing of Rankin--Selberg $L$-functions in the weight aspect}

\author{Hui Wang}
\address{The University of Hong Kong Shenzhen Institute of Research and Innovation,
Shenzhen, Guangdong 518057, China}
\email{wh0315@mail.sdu.edu.cn}
\author{Xin Wang}
\address{Yau Mathematical Sciences Center, Tsinghua University, Beijing 100048, China}
\email{wang\_xin@tsinghua.edu.cn}



\date{\today}

\begin{document}

\begin{abstract}
Fix a normalized holomorphic Hecke cusp form $g$ of level one and let
$f$ range over such forms. We prove an asymptotic formula with a
power-saving error term for the twisted second moment of
$L(1/2,f\otimes g)$, averaged smoothly over even weights near $K$ on
the scale $L$. For each fixed $1/3<\delta\leq1$, the formula holds
uniformly for $K^\delta\leq L\leq K$ and for twist indices up to a
positive power of $K$ depending on $\delta$.
Using a mollifier constructed from the reciprocal of the
Rankin--Selberg Euler product, we show that the harmonic proportion
under smooth averaging and the natural proportion over $K<k\leq K+L$
of forms $f$ with $L(1/2,f\otimes g)\ne0$ have the same explicit
asymptotic lower bound $1-(1+\min\{(3\delta-1)/12,25/174\})^{-3}$.
\end{abstract}
\maketitle

\section{Introduction}

\subsection{Background}
For standard $L$-functions of holomorphic cusp forms on $\GL(2)$,
central non-vanishing in level or weight families is a basic problem.
For forms associated with elliptic curves, the Birch and Swinnerton-Dyer
conjecture relates it to rank zero. Duke~\cite{Duke1995} proved that
$\gg q/(\log q)^2$ newforms of weight $2$ and prime level $q$ have
non-zero central values, while Kowalski and Michel~\cite{KowalskiMichel1999}
obtained a proportion of at least $1/6-\eps$ in the same family.
Iwaniec and Sarnak~\cite{IwaniecSarnak2000} proved an asymptotic lower
bound of $1/2$ among newforms with sign $+1$ in the functional equation,
at fixed even weight as
the square-free level $q$ grows with $\varphi(q)\sim q$. Pi, Wang and
Zhang~\cite[Theorem~1.2]{PiWangZhang} obtained a lower bound of
$1/4-\eps$ for the harmonic proportion of non-zero central values at
fixed even weight at least $4$ and level $N^3$, as the square-free
integer $N$ tends to infinity.
For quadratic twists of a fixed cuspidal representation on $\GL(2)$,
Friedberg and Hoffstein~\cite{FriedbergHoffstein1995} established
non-vanishing for infinitely many twists with admissible local conditions.

At level one, Iwaniec and Sarnak~\cite{IwaniecSarnak2000} obtained
the same lower bound of $1/2$ for the natural proportions
averaged over weights divisible by $4$ up to $K$.
For a single sufficiently large weight divisible by four,
Luo~\cite{Luo2015} proved a positive proportion, and Balkanova and
Frolenkov~\cite{BalkanovaFrolenkov2021} gave the explicit lower bound
$1/5-\eps$. For symmetric-square $L$-functions of level-one holomorphic
Hecke eigenforms, Khan~\cite[Theorem~1.1]{Khan} obtained an asymptotic
harmonic proportion of at least $19/27$ for non-zero central values by
averaging smoothly over weights on the scale $K$, and Balkanova and
Frolenkov~\cite[Theorem~1.1]{BalkanovaFrolenkov2026} recently improved
the lower bound in the same family to $1685/2197$.

For Rankin--Selberg convolutions, Kowalski, Michel and
VanderKam~\cite[Theorem~1.7]{KMV} used mollification to
obtain positive-proportion non-vanishing in the prime-level aspect
at fixed even weights $2<k<12$, both for central values and, in the self-dual case with root number
$-1$, for first derivatives at the central point. In the weight aspect,
Sarnak~\cite{S} obtained subconvexity bounds, while Lau, Liu and
Ye~\cite[Theorem~1 and Corollary~2]{LLY2006IMRP} bounded the second moment for
$K^{1/3+\eps}\leq L\leq K^{1-\eps}$ and obtained
a Weyl-type bound.
Ghosh~\cite{Ghosh2025} obtained an asymptotic formula for mixed second
moments of Rankin--Selberg $L$-functions over long weight intervals and
deduced simultaneous non-vanishing.

We study central non-vanishing in short weight intervals. Let
$H_k$ be the set of normalized holomorphic cuspidal Hecke
eigenforms of level one and weight $k$. We fix $g\in H_\ell$ and
let $f$ range over $H_k$ for even $k>\ell$. The Rankin--Selberg
$L$-function $L(s,f\otimes g)$ has root number \(+1\) and analytic conductor of order \(k^4\) at the central point.
We average smoothly over weights $k$ near $K$ on the scale $L$
and seek a positive harmonic proportion of non-zero central values.
The Hecke expansion of the mollifier introduces twists that can grow
with $K$, so we need moment asymptotics uniform in both $L$ and the twist.

\subsection{Main results}
Fix a non-negative, non-zero function $h\in C_c^\infty((0,\infty))$.
For $K,L>0$, put
$$h_{K,L}(k):=h\left(\frac{k-1-K}{L}\right).$$ 
We write $\sum^h$ for the harmonic sum in the Petersson trace formula.
For $a>0$, put $X=K^{2a}$. Define $\mu_{f,g}(n)$ by
\begin{align*}
\frac1{L(s,f\otimes g)}= & ~
\sum_{n\geq1}\frac{\mu_{f,g}(n)}{n^s},
\qquad \Re(s)>1.
\end{align*}
For a fixed real polynomial $P$ with $P(0)=P'(0)=0$, define
the mollifier by
\begin{align*}
M_P(f):= & ~\sum_{n\leq X}\frac{\mu_{f,g}(n)}{\sqrt n}
P\left(\frac{\log(X/n)}{\log X}\right).
\end{align*}
The upper bound for $a$ depends on $\delta$ in $L\geq K^\delta$.
For $1/3<\delta\leq1$, put
\begin{align*}
\mathfrak a(\delta):= & ~\min\bigl\{(3\delta-1)/12,25/174\bigr\}\\
= & ~\begin{cases}
(3\delta-1)/12,&1/3<\delta\leq79/87,\\
25/174,&79/87<\delta\leq1.
\end{cases}
\end{align*}

The choice of $a$ determines the length of the mollifier and
the lower bound for the harmonic proportion of non-zero central
values.

\begin{theorem}\label{thm:intro-nonvanishing}
Fix $1/3<\delta\leq1$. Let $a$ be fixed and suppose that
\begin{align}\label{eq:intro-parameters}
K^\delta\leq & ~L\leq K,
\qquad 0<a<\mathfrak a(\delta).
\end{align}
Then, uniformly for $L$ in this range as $K\to\infty$,
\begin{align}\label{eq:intro-nonvanishing}
\sum_{k\equiv0\bmod2}h_{K,L}(k)
\sideset{}{^h}\sum_{\substack{f\in H_k\\L(1/2,f\otimes g)\ne0}}1
\geq & ~\left(1-\frac1{(1+a)^3}
+o(1)\right)
\sum_{k\equiv0\bmod2}h_{K,L}(k)
\sideset{}{^h}\sum_{f\in H_k}1.
\end{align}
For fixed $\delta$, every harmonic proportion strictly below
$1-(1+\mathfrak a(\delta))^{-3}$ is obtained by a fixed
choice of $a$ in \eqref{eq:intro-parameters}.
\end{theorem}

\medskip
\noindent\textbf{Remark.}
For fixed $g$, the theorem gives a positive harmonic proportion
of non-zero central values in weight intervals on the scale
$L=K^{1/3+\eps}$ for every sufficiently small fixed $\eps>0$.
In this range $\mathfrak a(1/3+\eps)=\eps/4$, so the lower bound
$1-(1+\eps/4)^{-3}$ is asymptotic to $3\eps/4$ as $\eps$ tends to zero.
In particular, for $\delta=2/5$ we have $\mathfrak a(2/5)=1/60$. Thus for
$K^{2/5}\leq L\leq K$ the theorem gives every harmonic proportion
strictly below $1-(60/61)^3=0.048378\ldots$.
For $\delta=9/10$, we have $\mathfrak a(9/10)=17/120$, and
the theorem gives every harmonic proportion strictly below
$1-(120/137)^3=0.327980\ldots$.
For $79/87\leq\delta\leq1$ we have $\mathfrak a(\delta)=25/174$, giving every
proportion strictly below $1-(174/199)^3=0.331519\ldots$.
\smallskip

The same lower bound for the proportion holds when the forms
are counted with equal weight.

\begin{corollary}\label{cor:ordinary-nonvanishing}
Fix $1/3<\delta\leq1$ and $0<a<\mathfrak a(\delta)$.
As $K\to\infty$, uniformly for $K^\delta\leq L\leq K$, we have
\begin{align*}
 &\sum_{\substack{K<k\leq K+L\\k\equiv0\bmod2}}
 \#\bigl\{f\in H_k:L(1/2,f\otimes g)\ne0\bigr\}\geq
 \bigl(1-(1+a)^{-3}-o(1)\bigr)
 \sum_{\substack{K<k\leq K+L\\k\equiv0\bmod2}}\#H_k.
\end{align*}
For fixed $\delta$, every proportion strictly below
$1-(1+\mathfrak a(\delta))^{-3}$ is obtained by a fixed
choice of $a$ in this range.
\end{corollary}

Theorem~\ref{thm:intro-nonvanishing} follows by applying Cauchy's
inequality to the first two mollified moments and optimizing the
polynomial $P$. For the main terms of these moments, put
\begin{align*}
\mu_1(P,a):= & ~P(1)+a^{-1}P'(1),
\end{align*}
and
\begin{align*}
\mu_2(P,a):= & ~\mu_1(P,a)^2+
\frac1{3a^3}\int_0^1P''(u)^2\dd u.
\end{align*}

\begin{proposition}\label{prop:mollified-moments}
Let $1/3<\delta\leq1$ and assume \eqref{eq:intro-parameters}.
For every fixed real polynomial $P$ with $P(0)=P'(0)=0$
and every fixed positive integer $t$, we have
\begin{align}\label{eq:mollified-first-moment}
 &\sum_{k\equiv0\bmod2}h_{K,L}(k)
 \sideset{}{^h}\sum_{f\in H_k}
 \lambda_f(t) L(1/2,f\otimes g)M_P(f)=\frac L2\mathbf 1_{t=1}\mu_1(P,a)
 \int_{\mathbb R}h(u)\dd u+o_t(L),
\end{align}
and
\begin{align}\label{eq:mollified-second-moment}
 &\sum_{k\equiv0\bmod2}h_{K,L}(k)
 \sideset{}{^h}\sum_{f\in H_k}
 \lambda_f(t) L(1/2,f\otimes g)^2M_P(f)^2=\frac L2\mathbf 1_{t=1}\mu_2(P,a)
 \int_{\mathbb R}h(u)\dd u+o_t(L).
\end{align}
The errors are uniform for $K^\delta\leq L\leq K$ when
$(\delta,a)$ ranges over a compact subset of
$\{(\delta,a):1/3<\delta\leq1,\ 0<a<\mathfrak a(\delta)\}$,
with $g,h,P,t$ fixed.
\end{proposition}

\subsection{Overview of the proof}

In \S\,\ref{sec:fmcor} we apply the approximate functional equation
and the Petersson trace formula to the twisted first moment. Averaging
the Bessel function over the even weights makes the off-diagonal
negligible in the stated twist range. On the diagonal, a Mellin
contour shift followed by Poisson summation over the weights gives
\eqref{eq:fmcor-digamma}.

For the second moment in \S\,\ref{sec:ndnz-general-reduction}, we
square the approximate functional equation, obtaining two sums of
effective length about $K^2$. After applying the Petersson trace
formula, we open the Kloosterman sums and use Voronoi summation to
obtain Ramanujan sums. At zero shift the two Bessel arguments
coincide, and the same Hecke coefficients occur as on the diagonal.
By the gamma factor identity \eqref{eq:balanced-gamma-identity}, the
zero shift cancels the integral on the shifted diagonal contour. The
residues from the diagonal term and the zero shift give the complete main term in \eqref{eq:general-PV-residue}.

For the non-zero shifts, averaging over the weights restricts the
moduli, and the Bessel phases restrict the shifts. We estimate the
resulting shifted convolution sums by splitting the moduli according
to the size of the stationary phase parameter. For small parameters we
apply Bernard's estimate~\cite{Bernard2015} at each modulus. For large
parameters a quadratic change of variables and a rescaling give a
common phase for the shifts and moduli in each block. We can then sum
over the shifts inside the spectral expansion of Lau, Liu and
Ye~\cite{LLY2006IMRP} and apply the spectral large sieve~\cite{BM2015} and Michel's
triple product estimate~\cite{Michel2004}. The resulting bound \eqref{v6:eq:joint-error}
includes the continuous and exceptional spectra and gives the
twisted second moment asymptotic of Theorem~\ref{thm:twisted-second-moment},
with a power saving for twists up to a small power of $K$ throughout
$K^\delta\leq L\leq K$ for fixed $1/3<\delta\leq1$.

In \S\,\ref{sec:mollifier} we use the Hecke expansion of the mollifier,
including the prime-square terms in the reciprocal Euler product.
With the fixed factor $\lambda_f(t)$, the first and second moments
require twists up to $tX$ and $tX^2$, respectively. Summing the errors
gives the admissible mollifier length (see \eqref{eq:general-mollified-error-sum}).
For $t=1$, Perron inversion and the residue in \eqref{eq:moll-kernel-residue}
give the main term of the mollified second moment. For the first moment,
one numerator and one denominator cancel in the same Euler product.
Cauchy's inequality and the optimization of the mollifier polynomial
as in Khan~\cite[\S\,4]{Khan} give the harmonic non-vanishing bound
\eqref{eq:intro-nonvanishing}. For fixed $t>1$ both main terms vanish.
Approximating $L(1,\sym^2 f)$ in harmonic mean square by non-negative
Hecke sums supported on squares then gives the same lower bound for
the ordinary count.
\medskip


\noindent\textbf{Notation.}
Unless otherwise stated, $\eps>0$ is arbitrarily small and
$A>0$ is arbitrarily large. Both may vary from line to line.
As usual, we write $e(z)=\exp(2\pi\mathrm{i}z)$, and $\gamma$ denotes
Euler's constant. Our Fourier transform convention is
\[
\widehat h(t)=\int_{\mathbb R}h(u)e(-ut)\dd u,
\qquad
h(u)=\int_{\mathbb R}\widehat h(t)e(ut)\dd t.
\]
We write $f=O(g)$ or $f\ll g$ if $|f|\leq Cg$ for some constant
$C>0$, and $f\asymp g$ if $c_1g\leq|f|\leq c_2g$ for positive
constants $c_1,c_2$. The notation $k\sim K$ means $K<k\leq2K$.

\section{Preliminaries}
In this section, we present essential facts and tools needed later.
\subsection{Automorphic forms and Rankin--Selberg $L$-functions}
Any $f\in H_k$ has a Fourier expansion
\bna
f(z)=\sum_{n\geq1}\lambda_f(n)n^{(k-1)/2}e(nz),
\ena
where
$\lambda_f(n)$ is the $n$-th normalized Hecke eigenvalue of $f$ such that $\lambda_f(1)=1$.
In particular, the Fourier coefficient $\lambda_f(n)$ is real and
satisfies the Hecke relation
\bea\label{Hecke relation}
\lambda_f(m)\lambda_f(n)=\sum_{d|(m,n)}\lambda_f\left(\frac{mn}{d^2}\right),
\qquad
\lambda_f(m n)=\sum_{d|(m ,n)}\mu(d)\lambda_f\left(\frac{m}{d}\right)\lambda_f\left(\frac{n}{d}\right),
\eea
and the Ramanujan--Petersson bound (see Deligne \cite{Deligne})
\bna
|\lambda_f(n)|\leq \tau(n),
\ena
where $\tau(n)$ is the divisor function.
For $\Re(s)>1$, the $L$-function associated to $f$ is
defined by
\bna
L(s,f)=\sum_{n\geq 1}\frac{\lambda_f(n)}{n^s}.
\ena
Its completed $L$-function satisfies the functional equation
\bna
\Lambda(s,f)=(2\pi)^{-s}\Gamma\left(s+\frac{k-1}{2}\right)L(s,f)=\mathrm{i}^{-k}\Lambda(1-s, f).
\ena

Let $f\in H_k$ and $g$ be a fixed holomorphic cusp form of weight $\ell<k$ with Hecke eigenvalues $\lambda_g(n)$.
We consider the corresponding Rankin--Selberg convolution $L$-function, given for $\Re(s)>1$ by
\bea\label{Dirichlet series}
L(s,f\otimes g)=\zeta(2s)\sum_{n\geq 1}\frac{\lambda_f(n)\lambda_g(n)}{n^s},
\eea
which extends to an entire function of $s$. The completed $L$-function is
\bna
\Lambda(s,f\otimes g)=\left(\frac{1}{4\pi^2}\right)^{s}\Gamma\left(s+\frac{k-\ell}{2}\right)
\Gamma\left(s+\frac{k+\ell}{2}-1\right)L(s,f\otimes g),
\ena
which extends analytically to $\mathbb{C}$ and satisfies the functional equation
\bna
\Lambda(s,f\otimes g)=\Lambda(1-s, f\otimes g).
\ena
Let
\bna
\gamma(s, k)=(4\pi^2)^{-s}\Gamma\left(s+\frac{k-\ell}{2}\right)
\Gamma\left(s+\frac{k+\ell}{2}-1\right)
\ena
and $G(s)=e^{s^2}$. Put
\begin{equation}\label{eq:gamma-ratio}
\Gamma_k(s):=\frac{\gamma(1/2+s,k)}{\gamma(1/2,k)},
\end{equation}
and
\begin{equation*}
\alpha_k:=\Gamma_k'(0)=
\psi\left(\frac{k+\ell-1}{2}\right)
+\psi\left(\frac{k-\ell+1}{2}\right)-\log(4\pi^2),
\end{equation*}
where $\psi=\Gamma'/\Gamma$. 

By Stirling's asymptotic formula (see \cite[\S\,8.4, Eq.\,(4.03)]{O}),
for $|\arg s|\leq \pi-\eps$, $|s|\gg 1$ and any $\eps>0$,
\bna
\log \Gamma(s)=\left(s-\frac{1}{2}\right)\log s-s+\frac{1}{2}\log(2\pi)
+\sum_{j=1}^{K_1}\frac{B_{2j}}{2j(2j-1)s^{2j-1}}+O\left(\frac{1}{|s|^{2K_1+1}}\right),
\ena
where $B_j$ is the $j$-th Bernoulli number.
For fixed $\ell$, Stirling's formula gives
$\alpha_k=2\log k+O(1)$ as $k\to\infty$.
Moreover, for $x>0$ and $H$ holomorphic near $0$, we have
\begin{equation}\label{eq:common-mellin-residue}
2\Res_{s=0}\frac{G(s)\Gamma_k(s)\zeta(1+2s)}{s\,x^s}H(s)
 =\bigl(\alpha_k+2\gamma-\log x\bigr)H(0)+H'(0).
\end{equation}
Indeed, $G(0)=\Gamma_k(0)=1$, $G'(0)=0$, and
\bna
\zeta(1+2s)=\frac{1}{2s}+\gamma+O(s).
\ena

We now state the approximate functional equation for $L(1/2,f\otimes g)$.
\begin{lemma}\label{afe}
With the notation and assumptions above, 
we have
\bna
L\left(\frac{1}{2}, f\otimes g\right)=
2\sum_{n\geq 1}\frac{\lambda_f(n)\lambda_g(n)}{n^{1/2}}
V_k\left(n\right),
\ena
where, for any $y>0$,
\bea\label{V_k definition}
V_k(y)=\frac{1}{2\pi \mathrm{i}}\int_{(3)}y^{-s}\zeta(1+2s)
\Gamma_k(s)G(s)
\frac{\dd s}{s}.
\eea
\end{lemma}
\begin{proof}
The proof is based on \cite[Theorem 5.3]{IK}.
\end{proof}

We then have the following bounds for $V_k(y)$.
\begin{lemma}\label{cutoff of Vk}
For $y>0$,
the derivatives of $V_k(y)$ with respect to $y$ satisfy
\bna
y^a\frac{\partial^a}{\partial y^a}V_k(y)
\ll\left(1+\frac{y}{k^2}\right)^{-A},
\ena
where $a\geq1$ is an integer and $A>0$. The implied constant
may depend on $g$, $a$ and $A$.
For $a=0$ the logarithmic growth at zero satisfies
\begin{align*}
|V_k(y)|\ll  & ~\left(1+\log^+\frac{k^2}{y}\right)\left(1+\frac{y}{k^2}\right)^{-A},
\end{align*}
where $\log^+t=\max(\log t,0)$.

\end{lemma}
\begin{proof}
The proof is based on \cite[Proposition 5.4]{IK} and
\cite[Lemma 2.2]{Khan}. Applying $y^a\frac{\partial^a}{\partial y^a}$
under the integral in \eqref{V_k definition} inserts the factor
$(-s)(-s-1)\cdots(-s-a+1)$, interpreted as $1$ for $a=0$.
By Euler's integral and Stirling's formula,
\[
|\Gamma_k(\sigma+\mathrm{i}t)|\leq\Gamma_k(\sigma)\ll_{g,A}k^{2\sigma},
\qquad\text{for}\,\,\,-1/2\leq\sigma\leq\max\{3,A\}.
\]
For $y\geq k^2$, moving the line to $\Re(s)=A$ gives the bound
$O_{g,a,A}((k^2/y)^A)$ by the rapid decay of $G(s)$.
For $0<y\leq k^2$, move the line to $\Re(s)=-1/2$. By
\eqref{eq:common-mellin-residue},
\[
V_k(y)=\tfrac12(\alpha_k+2\gamma-\log y)
+O_g\left(\frac{y^{1/2}}{k}\right).
\]
For $a\geq1$, the same residue formula gives
\[
y^a\frac{\partial^a}{\partial y^a}V_k(y)
=\frac{(-1)^a(a-1)!}{2}
+O_{g,a}\left(\frac{y^{1/2}}{k}\right).
\]
Since $\alpha_k=2\log k+O_g(1)$ and $y^{1/2}\leq k$, these give
the required bounds.
\end{proof}

\subsection{The Euler product}
The $L$-function $L(s,\pi)$ attached to an automorphic representation
$\pi$ of $\mathrm{GL}(m)$ is given by a Dirichlet series and Euler product
\bna
L(s, \pi)=\sum_{n\geq 1}\frac{\lambda_\pi(n)}{n^s}=\prod_{p}L_p(s, \pi),
\ena
where $L_p(s, \pi)$ is the local $L$-factor given by
\bna
L_p(s, \pi)=\sum_{j=0}^\infty\frac{\lambda_\pi (p^j)}{p^{js}}.
\ena
For $\zeta(s)$ and the standard, symmetric-square and
Rankin--Selberg $L$-functions of $f$ and $g$ in the unitary
normalization used here, the Dirichlet series, Euler products and
local series converge absolutely for $\Re(s)>1$.
\begin{lemma}\label{Euler product}
For any positive integer $r$ and any modular form $g\in H_\ell$,
we have, for $\Re(s)>1$,
\bna
\sum_{n\geq 1}\frac{\lambda_g(n^2r)}{n^s}
=\frac{L(s,{\sym}^2g)}{\zeta(2s)}
\prod_{p^\nu\parallel r}
\frac{\lambda_g(p^\nu)-\lambda_g(p^{\nu-2})p^{-s}}
{1+p^{-s}},
\ena
where $p^\nu\parallel r$ means $v_p(r)=\nu\geq1$,
$\lambda_g(p^{-1})=0$ by convention, and the product is empty for $r=1$.
\end{lemma}

\begin{proof}
For $\nu,j\geq1$, apply the inverse Hecke relation in
\eqref{Hecke relation} to $(p^\nu,p^{2j})$ and $(p^{\nu-1},p^{2j-1})$.
Adding the two identities gives
\[
\lambda_g(p^\nu)\lambda_g(p^{2j})
-\lambda_g(p^{\nu-2})\lambda_g(p^{2j-2})
=\lambda_g(p^{2j+\nu})+\lambda_g(p^{2j+\nu-2}).
\]
Multiplying by $p^{-js}$, summing over $j\geq1$ and including
the term $j=0$, we get
\[
(1+p^{-s})\sum_{j\geq0}\lambda_g(p^{2j+\nu})p^{-js}
=\bigl(\lambda_g(p^\nu)-\lambda_g(p^{\nu-2})p^{-s}\bigr)
\sum_{j\geq0}\lambda_g(p^{2j})p^{-js}.
\]
The symmetric square local factor gives
\[
\sum_{j\geq0}\lambda_g(p^{2j})p^{-js}
=(1-p^{-2s})L_p(s,\sym^2g).
\]
Multiplicativity of \(\lambda_g\) and absolute convergence for $\Re(s)>1$
now give the stated Euler product.
\end{proof}

\begin{remark}\label{rem:zero-free-region}
Let $g\in H_\ell$ be fixed and write $s=\sigma+\mathrm{i}t$.
There is a constant $c_g>0$ such that $\zeta(s)$ and
$L(s,\sym^2g)$ have no zeros in
\begin{align*}
\sigma\geq1-\frac{c_g}{\log(3+|t|)}.
\end{align*}
This follows from the zero-free region for $L(s,g\otimes g)$
in \cite[Lemma~9.1]{KMV} and the factorization
$L(s,g\otimes g)=\zeta(s)L(s,\sym^2g)$.
At $s=1$, the simple poles of $L(s,g\otimes g)$ and $\zeta(s)$
give $L(1,\sym^2g)\ne0$.

After decreasing $c_g$, the functions
\begin{align*}
\zeta(s)-\frac{1}{s-1},\qquad \frac{1}{\zeta(s)},\qquad
L(s,\sym^2g),\qquad \frac{1}{L(s,\sym^2g)}
\end{align*}
and their derivatives of any fixed order are bounded in this region
by powers of $\log(3+|t|)$, with values at $s=1$ given by
holomorphic continuation. The bounds follow from the Euler products on
$\sigma=1+1/\log(3+|t|)$ and convexity,
Borel--Carath\'eodory in a slightly wider zero-free region
for the reciprocals, and Cauchy's formula for the derivatives.
\end{remark}

\subsection{Petersson trace formula}
For any complex number $\alpha_f$ depending on $f$, let
\bna
\sideset{}{^h}\sum_{f\in H_k}\alpha_f=\frac{2\pi^2}{k-1}\sum_{f\in H_k}\frac{\alpha_f}{L(1,\text{sym}^2f)}
\ena
denote the `harmonic' average of $\alpha_f$.
Note that $\sideset{}{^h}\sum_{f\in H_k}1=1+O(2^{-k})$ by the following Petersson trace formula.
\begin{lemma}\label{Petersson}
For any $m,n\geq 1$, we have
\bna
\sideset{}{^h}\sum_{f\in H_k}\lambda_f(m)\lambda_f(n)=\delta_{m,n}+2\pi \mathrm{i}^k\sum_{c=1}^\infty
\frac{S(m,n;c)}{c}J_{k-1}\left(\frac{4\pi\sqrt{mn}}{c}\right),
\ena
where the value of $\delta_{m,n}$ is $1$ if $m=n$ and $0$ otherwise,
$J_{k-1}$ is the Bessel function of the first kind and
\bna
S(m,n;c)=\sideset{}{^\star}\sum_{d\bmod c}e\left(\frac{md+n\overline{d}}{c}\right)
\ena
is the classical Kloosterman sum.
Here $\sideset{}{^\star}\sum$ restricts the summation to the primitive residue classes and
$\overline{d}$ denotes the multiplicative inverse of $d$ modulo $c$.
\end{lemma}
\subsection{$J$-Bessel function}
For the Bessel function $J_{k-1}$, we use the standard estimate
(see \cite[\S\S\,2.11, 7.1]{Watson})
\bna
J_{k-1}(x)\ll_k\left\{
\begin{aligned}
&x^{k-1}, \ &&\textrm{if\ \ }x\leq 1,\\
&x^{-1/2}, \ && \textrm{if\ \ }x>1.
\end{aligned}
\right.
\ena
We have the Bessel differential equation
\begin{equation}\label{Bessel differential equation}
y^2J_\alpha''(y)+yJ_\alpha'(y)
+(y^2-\alpha^2)J_\alpha(y)=0,
\end{equation}
and the derivative recurrence 
\bea\label{derivative recurrence}
2J_\alpha'(y)=J_{\alpha-1}(y)-J_{\alpha+1}(y)
\eea
for every $\alpha$ (see \cite[8.401, 8.471(2)]{GR}).
Moreover, we have
\bna
|J_{k-1}(x)|\leq 1, \quad\text{for}\,\, \text{all}\,\, x\geq 0,
\ena
and, as $k\to\infty$, uniformly for positive $x$ in the indicated range,
\bna
J_{k-1}(x)\sim (1/\Gamma(k))(x/2)^{k-1},\quad \text{for } \,\, x\leq k^{1/2-\eps}.
\ena

With $h_{K,L}$ as in the introduction, we have the following exact
Poisson representation for the average over the even weights.
\begin{lemma}\label{lem:exact-bessel-average}
Let $K,L>0$, let $x>0$, and let $h\in C_c^\infty((0,\infty))$.
The average of the $J$-Bessel function has the exact representation
\begin{equation}\label{after Poisson}
\begin{split}
\sum_{k\equiv 0\bmod 2}2\mathrm{i}^k h_{K, L}(k)J_{k-1}(x)
=\frac{e^{-\mathrm{i}x}W_{K, L}(x)-e^{\mathrm{i}x}W_{K, L}(-x)}{\mathrm{i}},
\end{split}
\end{equation}
where
\bea\label{W_KL}
W_{K, L}(x)=\int_{-\infty}^\infty \widehat{h}(t)e\left(-\frac{K t}{L}+\frac{x}{2\pi}\left(1-\cos\frac{2\pi t}{L}\right)\right)\dd t.
\eea
\end{lemma}
\begin{proof}
Apply Poisson summation on the even integers to the Fourier
representation of $J_{k-1}$ (see \cite[Lemma~5.8]{I} and
\cite[Eqs.\,(53)--(55)]{S}). This gives the displayed identity
without a Taylor approximation to the trigonometric phase.
\end{proof}


\subsection{Voronoi summation formula}
We use the Voronoi summation formula in the form recorded in Blomer and Harcos \cite[Proposition 2]{BH}.
\begin{lemma}\label{Voronoi summation}
Let $a$ and $c$ be coprime positive integers, and let $F: (0, \infty)\rightarrow \mathbb{C}$ be a smooth function of compact support. Then
\bna
c\sum_{n= 1}^\infty\lambda_g(n)e\left(\frac{na}{c}\right)F(n)
=2\pi \mathrm{i}^\ell \sum_{n= 1}^\infty\lambda_g(n)e\left(-\frac{\overline{a} n}{c}\right)\int_0^\infty F(x)J_{\ell-1}\left(\frac{4\pi \sqrt{n x}}{c}\right)\dd x.
\ena
\end{lemma}

\subsection[Stationary phase]{Stationary phase}
We recall a nonstationary phase estimate from
\cite[Lemma~8.1]{BKY} in the language of inert functions.
Put
\begin{align*}
I:=\int_{\mathbb R}w(y)e^{\mathrm{i}\varrho(y)}\dd y.
\end{align*}

Let $\mathcal{F}$ be an index set and let
$Y:\mathcal{F}\to[1,\infty)$. Write $Y_T=Y(T)$.
A family $\{w_T\}_{T\in \mathcal{F}}$ of smooth
functions supported on a product of dyadic
intervals in $\mathbb{R}_{>0}^d$
is called $Y$-inert if for
each $j=(j_1,\ldots,j_d) \in \mathbb{Z}_{\geq 0}^d$
we have
$$
C(j_1,\ldots,j_d)
= \sup_{T \in \mathcal{F} }
\sup_{(y_1, \ldots, y_d) \in \mathbb{R}_{>0}^d}
Y_T^{-j_1- \cdots -j_d}\left| y_1^{j_1} \cdots y_d^{j_d}
w_T^{(j_1,\ldots,j_d)}(y_1,\ldots,y_d) \right| < \infty.
$$

\begin{lemma}
\label{lemma:exponentialintegral}
Suppose that $w=w_T(y)$ is a family of $Y$-inert functions
supported on $[Z,2Z]$. Let $\varrho$ be real valued and smooth.
Suppose that $\varrho^{(j)}(y)\ll_j H/Z^j$ for $j\geq1$ and
$|\varrho'(y)|\gg H/Z$ on the support of $w$, where
$H/Y^2\geq R\geq1$. Then $I\ll_A ZR^{-A}$ for every $A>0$.
\end{lemma}

\section{The first twisted moment}
\label{sec:fmcor}
With the notation of the introduction, put
\begin{align*}
\mathcal M_1(r):= & ~
\sum_{k\equiv0\bmod2}h_{K,L}(k)
\sideset{}{^h}\sum_{f\in H_k}
L(1/2,f\otimes g)\lambda_f(r).
\end{align*}

\begin{theorem}\label{thm:fmcor}
Fix $\eps>0$. Suppose that $K^\eps\leq L\leq K$ and that
$r$ is an integer with
\begin{align*}
1\leq r\leq L^2K^{-\eps}.
\end{align*}
Then for every $A>0$,
\begin{align}\label{eq:fmcor-digamma}
\mathcal M_1(r)
=\frac{\lambda_g(r)}{\sqrt r}L
\int_{\mathbb R}h(u)
\left(\gamma-\frac12\log r+
\frac12\alpha_{K+1+Lu}\right)\dd u
+O(K^{-A}).
\end{align}
\end{theorem}

\begin{proof}
On the support of $h_{K,L}$ we have $k\asymp K$.
The approximate functional equation (see Lemma \ref{afe}) and the Petersson trace formula (see Lemma~\ref{Petersson}) give
\begin{align}\label{eq:fmcor-trace}
\mathcal M_1(r)
= & ~\frac{2\lambda_g(r)}{\sqrt r}
\sum_{k\equiv0\bmod2}h_{K,L}(k)V_k(r)+\mathcal R_1(r),
\end{align}
where
\begin{align*}
\mathcal R_1(r)
= & ~2\pi\sum_{n,c\geq1}
\frac{\lambda_g(n)S(n,r;c)}{\sqrt n\,c}
\sum_{k\equiv0\bmod2}2\mathrm{i}^kh_{K,L}(k)V_k(n)
J_{k-1}\left(\frac{4\pi\sqrt{nr}}c\right).
\end{align*}
The bounds $|S(n,r;c)|\leq c$ and $J_{k-1}(y)\ll\min(1,y^2)$, the latter holding uniformly for integers \(k\geq3\), give
\[
\sum_{c\geq1}\frac{|S(n,r;c)|}{c}
\left|J_{k-1}\left(\frac{4\pi\sqrt{nr}}c\right)\right|
\ll \sqrt{nr}.
\]
Lemma~\ref{cutoff of Vk} and Deligne's bound therefore restricts $n\leq K^{2+\eps/4}$
with error $O(K^{-A})$.
On this range $nr\leq K^4$. For $c>K^3$ the Bessel series gives
$$|J_{k-1}(4\pi\sqrt{nr}/c)|\ll\frac{(2\pi K^2/c)^{k-1}}{\Gamma(k)}.$$
Since $k\asymp K$, summation over $c>K^3$ saves an arbitrary power of $K$. Put
\begin{align*}
x:=\frac{4\pi\sqrt{nr}}c
\ll & ~K^{1+\eps/8}\sqrt r
\leq K^{1-3\eps/8}L.
\end{align*}

For fixed $n$, put $H_n(v)=h(v)V_{K+1+Lv}(n)$.
On the support of $h$, write $k=K+1+Lv\asymp K$. We move the line in
\eqref{V_k definition} to $\Re(s)=\sigma=1/\log K$. On this line, Euler's integral for the
Gamma function and Stirling's formula give
\[
|\Gamma_k(s)|\leq\Gamma_k(\sigma)\ll K^{2\sigma}\ll1.
\]
For $j\geq1$ and $\Re(z)\geq1$, the polygamma series gives
\[
\psi^{(j)}(z)=(-1)^{j+1}j!\sum_{m=0}^{\infty}\frac{1}{(m+z)^{j+1}},
\qquad\text{and}\qquad
|\psi^{(j)}(z)|\ll_j(\Re(z))^{-j}.
\]
For $\Re(z)\asymp K$, integration along the segment from $z$ to $z+s$ gives
\[
|\psi^{(j-1)}(z+s)-\psi^{(j-1)}(z)|\ll_j |s|K^{-j}.
\]
Using this bound with $z=(k+\ell-1)/2$ and $z=(k-\ell+1)/2$ in
\eqref{eq:gamma-ratio}, we obtain
$\bigl|\frac{\dd^j}{\dd v^j}\log\Gamma_k(s)\bigr|
\ll_j (L/K)^j|s|\ll_j |s|$ since $L\leq K$.
The product rule and $|\Gamma_k(s)|\ll1$ then give
\[
\left|\frac{\dd^j}{\dd v^j}\Gamma_{K+1+Lv}(s)\right|
\ll_j (1+|\Im(s)|)^j,
\qquad\text{for}\,\,\,j\geq0.
\]

Write $s=\sigma+\mathrm{i}t$. For $|t|\leq1$, the expansion
$\zeta(1+2s)=1/(2s)+O(1)$ gives
\[
\left|\frac{\zeta(1+2s)}{s}\right|
\ll \frac1{\sigma^2+t^2},
\qquad\text{and}\qquad
\int_{-1}^1\frac{\dd t}{\sigma^2+t^2}\ll \sigma^{-1}=\log K.
\]
For $|t|>1$, we use $|\zeta(1+2s)|\leq\zeta(1+2\sigma)\ll\log K$
and the rapid decay of $G(s)$ to obtain the same bound. These
estimates justify differentiation under the integral and give
$\frac{\dd^j}{\dd v^j}V_{K+1+Lv}(n)\ll_j\log K$, uniformly in $n\geq1$.
It follows that $H_n^{(j)}(v)\ll_j\log K$.
Since $H_n$ has fixed compact support, differentiation of its
Fourier integral and repeated integration by parts give
\[
\left|\frac{\dd^j}{\dd t^j}\widehat H_n(t)\right|
\ll_{j,N}(\log K)\,(1+|t|)^{-N},
\qquad \text{for}\,\,\,j,N\geq0.
\]

Apply \eqref{after Poisson} with $h$ replaced by $H_n$. It suffices
to estimate $W_{K,L}(\pm x)$ in \eqref{W_KL} with this weight.
Suppose first that $K/L\geq K^{\eps/8}$. Insert a smooth cutoff
at scale $K^{\eps/8}$, equal to $1$ for $|t|\leq K^{\eps/8}$
and supported on $|t|\leq 2K^{\eps/8}$. By the rapid decay of
$\widehat H_n(t)$, this introduces an error of $O_B(K^{-B})$
for any $B>0$. Write the phase of $W_{K,L}(\pm x)$ as
\[
\phi_\pm(t)=-\frac{Kt}{L}
\pm\frac{x}{2\pi}\left(1-\cos\frac{2\pi t}{L}\right).
\]
Since $|x|\ll K^{1-3\eps/8}L$, on the support of the cutoff we have
\[
\left|\frac{x}{L}\sin\frac{2\pi t}{L}\right|
\ll \frac{|x|K^{\eps/8}}{L^2}
\ll \frac KL K^{-\eps/4}.
\]
It follows that
\[
\left|\frac{\dd\phi_\pm}{\dd t}\right|
=\left|-\frac KL\pm\frac{x}{L}\sin\frac{2\pi t}{L}\right|
\asymp \frac KL.
\]
For $j\geq2$, we have
$\bigl|\frac{\dd^j\phi_\pm}{\dd t^j}\bigr|\ll_j |x|/L^j\ll_j K/L$.
These phase bounds and the estimates for $\widehat H_n$ give
$O_N((L/K)^N\log K)$ for the truncated integral after $N$
integrations by parts. Since $L/K\leq K^{-\eps/8}$, taking $N$ large
gives $W_{K,L}(\pm x)=O_B(K^{-B})$ for any $B>0$.

Suppose now that $K/L<K^{\eps/8}$, then
$
|x|/L^2\ll K^{-\eps/4}.
$
Taylor expansion through degree $N-1$ gives
\begin{align*}
\exp\left(\pm\mathrm{i}x\left(1-\cos\frac{2\pi t}{L}\right)\right)
=\sum_{m=0}^{N-1}\frac{(\pm\mathrm{i}x)^m}{m!}
\left(1-\cos\frac{2\pi t}{L}\right)^m+R_N(t).
\end{align*}
Since $\cos(2\pi t/L)=\tfrac12(e(t/L)+e(-t/L))$, the power
$(1-\cos(2\pi t/L))^m$ is a finite linear combination of $e(jt/L)$
with $|j|\leq m$. By Fourier inversion and the positive support of $h$, the corresponding integrals are
\[
\int_{\mathbb R}\widehat H_n(t)e\left(\frac{j-K}{L}t\right)\dd t
=H_n\left(\frac{j-K}{L}\right)=0
\]
for $K>N$.
We have
$|R_N(t)|\ll_N |x|^N|1-\cos(2\pi t/L)|^N$ for all real $t$.
Using $1-\cos y\leq y^2/2$ and the decay of $\widehat H_n$, we obtain
\[
|W_{K,L}(\pm x)|
\ll_N \left(\frac{|x|}{L^2}\right)^N
\int_{\mathbb R}|\widehat H_n(t)|\,|t|^{2N}\dd t
\ll_N K^{-N\eps/4}\log K.
\]
Taking $N$ large gives $W_{K,L}(\pm x)=O_B(K^{-B})$ for any $B>0$.

Thus \eqref{after Poisson} gives $O_B(K^{-B})$ for the sum over $k$
in $\mathcal R_1(r)$ in both cases. On the truncated ranges, Deligne's
bound and $|S(n,r;c)|\leq c$ give
\[
\mathcal R_1(r)\ll_{A,B}
K^{-B}\sum_{c\leq K^3}\sum_{n\leq K^{2+\eps/4}}
\frac{\tau(n)}{\sqrt n}+K^{-A}
\ll_{A,B} K^{-B+5+\eps/4}+K^{-A},
\]
which is $O(K^{-A})$
by choosing $B$ sufficiently large.

For the diagonal term, Poisson summation on the even integers and
the rapid decay of $\widehat H_r$ give, for $N\geq2$,
\[
2\sum_{k\equiv0\bmod2}h_{K,L}(k)V_k(r)
=L\int_{\mathbb R}H_r(u)\dd u+O_N(L^{1-N}\log K).
\]
For $k=K+1+Lu$ on the support of $h$, move the line in
\eqref{V_k definition} to $\Re(s)=-C$ for fixed $C>0$. For large $K$
only the double pole at $s=0$ is crossed. By
\eqref{eq:common-mellin-residue} with $x=r$ and $H=1$, and using
$|\Gamma_k(-C+\mathrm{i}t)|\leq\Gamma_k(-C)\ll_C k^{-2C}$, we obtain
\[
2V_k(r)=\alpha_k+2\gamma-\log r+O_C\bigl((r/k^2)^C\bigr).
\]
Insert this into the integral above. Since $r/k^2\ll K^{-\eps}$,
the integrated remainder is $O_C(LK^{-C\eps})$. As $L\geq K^\eps$,
taking $C$ and $N$ sufficiently large in \eqref{eq:fmcor-trace}
gives \eqref{eq:fmcor-digamma}.
\end{proof}

\section{The twisted second moment}
\label{sec:ndnz-general-reduction}
For a positive integer $r$, put
\begin{equation*}
\mathcal M_2(r):=\sum_{k\equiv0\bmod2}h_{K,L}(k)
\sideset{}{^h}\sum_{f\in H_k}
L\left(\frac12,f\otimes g\right)^2\lambda_f(r).
\end{equation*}

\begin{theorem}\label{thm:twisted-second-moment}
Fix $1/3<\delta\leq1$, and put
\begin{align*}
\rho_*(\delta):= & ~
\min\left\{\delta-\frac13,\frac{50}{87}\right\}.
\end{align*}
Fix $0\leq\rho<\rho_*(\delta)$ and define
\begin{align*}
\sigma(\delta,\rho):= &\,
\min\left\{\frac{3\delta-1-3\rho}{2},
                                    \frac{50-87\rho}{64}\right\}>0.
\end{align*}
Then, for every fixed $\eps>0$, uniformly for
$K^\delta\leq L\leq K$ and all positive
integers $r\leq K^\rho$, we have
\begin{align}\label{eq:twisted-second-moment}
\mathcal M_2(r)=\sum_{k\equiv0\bmod2}h_{K,L}(k)\mathscr P_k(r)
+O\left(\frac{L}{\sqrt r}K^{-\sigma(\delta,\rho)+\eps}\right),
\end{align}
where $\mathscr P_k(r)$ is the complete main term in
\eqref{eq:general-PV-residue}.
\end{theorem}

\subsection{The main term}
Assume $1\leq r\leq K$ and $K^\eps\leq L\leq K$.
Squaring the approximate functional equation for $L(1/2,f\otimes g)$,
applying the Hecke relation and the Petersson trace formula, and
replacing $n$ by $nd$ in the terms with $d\mid(n,r)$, we obtain

\begin{equation}\label{the second moment}
\begin{split}
\mathcal M_2(r)={}&\frac4{\sqrt r}
\sum_{k\equiv0\bmod2}h_{K,L}(k)\sum_{d\mid r}
\sum_{n\geq1}\frac{\lambda_g(nr/d)\lambda_g(nd)}n
V_k(nr/d)V_k(nd)\\
&+8\pi\sum_{k\equiv0\bmod2}\mathrm{i}^kh_{K,L}(k)\sum_{d\mid r}
\sum_{m,n\geq1}\frac{\lambda_g(m)\lambda_g(nd)}{\sqrt{mnd}}
V_k(m)V_k(nd)\\
&\hspace{40mm}\times\sum_{c\geq1}\frac{S(m,nr/d;c)}c
J_{k-1}\left(\frac{4\pi\sqrt{mnr/d}}c\right).
\end{split}
\end{equation}
We denote the diagonal and off-diagonal contributions in
\eqref{the second moment} by $\mathcal D(r)$ and $\mathcal{ND}(r)$,
respectively.
Before Voronoi summation, the two weights give the separate
effective ranges $m,nd\ll K^{2+\eps}$.

\subsubsection{The diagonal term}
Put $q=r/d$. 
For $d\mid r$, define, initially for $\Re(s)>0$,
\begin{align}\label{eq:general-D-dq}
\mathscr D_{d,q}(s):= & ~\sum_{n\geq1}
\frac{\lambda_g(nd)\lambda_g(nq)}{n^{1+s}}.
\end{align}
For $p\mid r$, put $a=v_p(d)$ and $b=v_p(q)$, with
$\lambda_g(p^{-1})=0$. By \eqref{Hecke relation}, for $|x|<1$,
\begin{align*}
(1-x)\sum_{j\geq0}
\lambda_g(p^{j+a})\lambda_g(p^{j+b})x^j
=\sum_{j\geq0}\lambda_g(p^{2j+a+b})x^j
+\lambda_g(p^{a-1})\lambda_g(p^{b-1}).
\end{align*}
Using the local factors in Lemma~\ref{Euler product} and taking Euler products,
we obtain
\begin{align}\label{eq:general-D-euler}
\mathscr D_{d,q}(s)= & ~
\frac{\zeta(1+s)L(1+s,\sym^2g)}{\zeta(2+2s)}
\prod_{p\mid r}F_{p;v_p(d),v_p(q)}(p^{-1-s}),
\end{align}
where, for $a+b\geq1$,
\begin{align*}
F_{p;a,b}(x):= 
\frac{\lambda_g(p^{a+b})-\lambda_g(p^{a+b-2})x}{1+x}
+\lambda_g(p^{a-1})\lambda_g(p^{b-1})
\frac{1-(\lambda_g(p)^2-2)x+x^2}{1+x}.
\end{align*}
Let $\mathcal D_k(r)$ denote the diagonal term before averaging
over $k$. Substituting \eqref{V_k definition} into the diagonal
term of \eqref{the second moment} gives
\begin{align}\label{eq:balanced-diagonal-Mellin}
\mathcal D_{k}(r)= \frac4{\sqrt r}\sum_{d\mid r}
\frac1{(2\pi \mathrm{i})^2}\int_{(1)}\int_{(2)}& ~
G(u)\Gamma_k(u)G(t)\Gamma_k(t)q^{-u}d^{-t}\notag\\
& ~\times\zeta(1+2u)\zeta(1+2t)
\mathscr D_{d,q}(u+t)\frac{\dd t}{t}\frac{\dd u}{u}. 
\end{align}
The sums and integrals converge absolutely. We combine this
integral with the zero shift of the off-diagonal term before
moving the remaining contour.

\subsubsection[The off-diagonal term]{The off-diagonal term}
We return to the off-diagonal term in \eqref{the second moment}.
Choose a non-negative smooth function $U$ supported on $[1/2,2]$
such that
\begin{equation}\label{eq:ndnz-dyadic-partition}
\sum_{j\in\mathbb Z}U(x/2^j)=1,\qquad\text{for}\,\,\, x>0.
\end{equation}
Insert this partition of unity into the original $m$-sum,
so that the block indexed by $M=2^j$ carries the factor $U(m/M)$.
After opening the Kloosterman sum, fix $k,n,c,d$ and $M$.
Apply Lemma~\ref{Voronoi summation}
to this block with
\[
F(x)=x^{-1/2}V_k(x)
J_{k-1}\left(\frac{4\pi\sqrt{x qn}}c\right)U(x/M)
\in C_c^\infty((0,\infty)).
\]
The resulting Bessel transform is
\begin{equation}\label{TKM}
\mathcal T_{k,M}(m,qn;c):=
\int_0^\infty F(x)
J_{\ell-1}\left(\frac{4\pi\sqrt{xm}}c\right)\dd x.
\end{equation}
Here $M$ denotes the scale of the variable before Voronoi summation.
For fixed $k>4$, $q,n,c$ and every $A>0$,
Lemma~\ref{cutoff of Vk}, the power series for the Bessel function in
\cite[8.402]{GR}, and two applications of
\eqref{Bessel differential equation}, followed by integration by parts, give
\begin{equation*}
|\mathcal T_{k,M}(m,qn;c)|
\ll_{k,\ell,qn,c,A}m^{-2}
\begin{cases}
M^{(k-4)/2}(1+|\log M|),&0<M\leq1,\\
M^{-A},&M>1.
\end{cases}
\end{equation*}
Hence the dual sums over $m$ and $M=2^j$ converge absolutely.
The original sum is also absolutely summable, so summing the
Voronoi identities over $M$ gives
\begin{align}\label{eq:balanced-Voronoi}
\mathcal{ND}(r)= 16\pi^2\mathrm{i}^\ell
\sum_{k\equiv0\bmod2}\mathrm{i}^kh_{K,L}(k)\sum_{d\mid r}
\sum_{n,m,c\geq1}
\frac{\lambda_g(nd)\lambda_g(m)}{\sqrt{nd}}V_k(nd)
\frac{R_c(qn-m)}{c^2}\mathcal T_k(m,qn;c),
\end{align}
where \(R_c(qn-m)\) is the Ramanujan sum
\[
R_c(qn-m)=\sideset{}{^\star}\sum_{\alpha\bmod c}
e\left(\frac{(qn-m)\overline\alpha}{c}\right),
\qquad R_c(0)=\varphi(c),
\]
and the full kernel, obtained after summing the dyadic blocks, is
\begin{align*}
\mathcal T_k(m,qn;c)
:={}&\sum_{j\in\mathbb Z}\mathcal T_{k,2^j}(m,qn;c)\notag\\
={}&\int_0^\infty x^{-1/2}V_k(x)
J_{k-1}\left(\frac{4\pi\sqrt{x qn}}c\right)
J_{\ell-1}\left(\frac{4\pi\sqrt{xm}}c\right)\dd x.
\end{align*}
The outer $n,c$-sums are first taken with smooth cutoffs whose
normalized derivatives are uniformly bounded. For $m=qn$,
absolute convergence in \eqref{eq:balanced-zero-Mellin}
removes these cutoffs. For $m\neq qn$, they are kept until
the uniform estimates in \S\,\ref{sec:nonzero-reduction} apply.
Separating $m=qn$ from $m\neq qn$ in
\eqref{eq:balanced-Voronoi}, we write
\begin{align}\label{eq:ndnz-general-balanced}
\mathcal M_2(r)= & ~\mathcal D(r)+\mathcal Z(r)+\mathcal E(r),
\end{align}
where $\mathcal Z(r)$ is the zero shift \((m=nq)\) and
\begin{align*}
\mathcal E(r)= & ~16\pi^2\mathrm{i}^\ell
\sum_{k\equiv0\bmod2}\mathrm{i}^kh_{K,L}(k)\sum_{d\mid r}\sum_{\substack{n,m,c\geq1\\m\neq qn}}
\frac{\lambda_g(nd)\lambda_g(m)}{\sqrt{nd}}V_k(nd)
\frac{R_c(qn-m)}{c^2}\mathcal T_k(m,qn;c).
\end{align*}
Here and below, the sum defining $\mathcal E(r)$ is interpreted
using the common outer $n,c$ cutoffs above. At each finite stage
it equals the difference $\mathcal{ND}(r)-\mathcal Z(r)$ with
the same cutoffs applied to both terms. Absolute convergence of the
original off-diagonal sum in \eqref{the second moment} and of the
zero shift on the initial lines in \eqref{eq:balanced-zero-Mellin}
makes this limit independent of the cutoffs. The estimates in
\S\,\ref{sec:nonzero-reduction} are uniform in these cutoffs and
therefore pass to the limit.
\par

\subsubsection{The zero shift $\mathcal Z(r)$}
For $m=qn$, make the change of variables
\bna
y=\frac{4\pi\sqrt{xqn}}c,
\qquad x=\frac{c^2y^2}{16\pi^2qn},
\qquad x^{-1/2}\dd x=\frac{c}{2\pi\sqrt{qn}}\dd y.
\ena
Writing $\mathcal Z_{k}(r)$ for the term before averaging over $k$,
we obtain
\begin{equation*}
\mathcal Z_{k}(r)=\frac{8\pi \mathrm{i}^{k+\ell}}{\sqrt r}
\sum_{d\mid r}\sum_{n,c\geq1}
\frac{\lambda_g(nd)\lambda_g(nq)}n V_k(nd)\frac{\varphi(c)}c
\int_0^\infty
V_k\left(\frac{c^2y^2}{16\pi^2qn}\right)
J_{k-1}(y)J_{\ell-1}(y)\dd y.
\end{equation*}
Using \eqref{V_k definition} for both weights and
\[
\sum_{c\geq1}\frac{\varphi(c)}{c^{1+2s}}
=\frac{\zeta(2s)}{\zeta(1+2s)},
\]
we get
\begin{multline}\label{eq:balanced-zero-Mellin}
\mathcal Z_{k}(r)=\frac{8\pi \mathrm{i}^{k+\ell}}{\sqrt r}
\sum_{d\mid r}\frac1{(2\pi \mathrm{i})^2}\int_{(1)}\int_{(2)}
G(s)\Gamma_k(s)G(t)\Gamma_k(t)
(16\pi^2q)^s d^{-t}\\
\times\zeta(2s)\zeta(1+2t)
\mathscr D_{d,q}(t-s)\mathcal J_{k,\ell}(s)
\frac{\dd t}{t}\frac{\dd s}{s},
\end{multline}
where
\[
\mathcal J_{k,\ell}(s):=\int_0^\infty
y^{-2s}J_{k-1}(y)J_{\ell-1}(y)\dd y.
\]
The series over $c,n$ and the Bessel integral converge absolutely
on $\Re(s)=1$, $\Re(t)=2$, for sufficiently large $k$.

\subsubsection{The combined archimedean factor}
Let
\[
A_+=\frac{k+\ell-1}{2},\qquad A_-=\frac{k-\ell+1}{2},\qquad
\rho_k=\frac{(-1)^{(k-\ell)/2}}{2\pi}.
\]
The Weber--Schafheitlin formula (see \cite{Watson} or \cite[Eq.\,(6.574.2)]{GR}) gives, initially
for $0<2\Re(s)<k+\ell-1$,
\begin{equation*}
\mathcal J_{k,\ell}(s)
=2^{-2s}\frac{\Gamma(2s)\Gamma(A_+-s)}
{\Gamma(A_-+s)\Gamma(A_++s)\Gamma(1-A_-+s)}.
\end{equation*}
Since $A_-$ is a half integer, the reflection formula gives
\[
\frac1{\Gamma(1-A_-+s)}
=\frac{(-1)^{(k-\ell)/2}}\pi\cos(\pi s)\Gamma(A_--s).
\]
Using the functional equation
$\zeta(1-2s)=2(2\pi)^{-2s}\Gamma(2s)\cos(\pi s)\zeta(2s)$, the reflection formula above,
and \eqref{eq:gamma-ratio}, we obtain
\begin{align}\label{eq:balanced-gamma-identity}
\Gamma_k(s)(16\pi^2)^s\zeta(2s)\mathcal J_{k,\ell}(s)
=&~\rho_k\Gamma_k(-s)\zeta(1-2s).
\end{align}

\subsubsection{Combining the diagonal and the zero shift}
Apply \eqref{eq:balanced-gamma-identity} in
\eqref{eq:balanced-zero-Mellin} and put $u=-s$. Since $G$ is even
and $\mathrm{i}^{k+\ell}\rho_k=(2\pi)^{-1}$, the prefactor agrees
with that of \eqref{eq:balanced-diagonal-Mellin}. Reversing the
$u$-contour supplies a minus sign, so $\mathcal Z_k(r)$ is the
negative of \eqref{eq:balanced-diagonal-Mellin} with the $u$-integral
on $\Re(u)=-1$. Thus we move the $u$-contour in
\eqref{eq:balanced-diagonal-Mellin} from $\Re(u)=1$ to $-1$,
keeping $\Re(t)=2$.
Only $u=0$ is crossed, since $u=-t$ and all gamma poles lie
to the left of this strip. Applying \eqref{eq:common-mellin-residue} with $x=q$ and $H(u)=\mathscr D_{d,q}(u+t)$ gives
\begin{align}\label{eq:ndnz-general-combination}
\mathcal D_{k}(r)+\mathcal Z_{k}(r)
= \frac4{\sqrt r}\sum_{d\mid r}& ~\frac1{2\pi \mathrm{i}}
\int_{(2)}G(t)\Gamma_k(t)\zeta(1+2t)d^{-t}\notag\\
& ~\times\left\{
\left(\gamma+\frac{\alpha_k-\log q}{2}\right)\mathscr D_{d,q}(t)
+\frac12\mathscr D_{d,q}'(t)\right\}\frac{\dd t}{t}.
\end{align}
The term $-\frac12\log q$ comes from $q^{-u}$.
\begin{proposition}
Suppose that $K^\eps\leq L\leq K$. For every integer
$1\leq r\leq K$, one has
\begin{align}\label{eq:general-integer-second}
\mathcal M_2(r)= & ~\sum_{k\equiv0\bmod2}h_{K,L}(k)
\mathscr P_k(r)+\mathcal E(r)+O(LK^{-1+\eps}),
\end{align}
where \begin{align}\label{eq:general-PV-residue}
\mathscr P_k(r):= \frac4{\sqrt r}\sum_{d\mid r}& ~
\Res_{s=0}\frac{G(s)\Gamma_k(s)\zeta(1+2s)}{sd^{s}}
\notag\\
& ~\times\left\{
\left(\gamma+\frac{\alpha_k-\log(r/d)}2\right)
\mathscr D_{d,r/d}(s)+\frac12\mathscr D_{d,r/d}'(s)
\right\}.
\end{align}
The residue is taken of the entire product and keeps the cubic main term, all lower logarithmic terms,
and the local factors in \eqref{eq:general-D-euler}.
\end{proposition}

\begin{proof}
Move the contour in \eqref{eq:ndnz-general-combination}
to $\Re(t)=-1/2+\eps/8$.
By \eqref{eq:general-D-euler}, only $t=0$ is crossed.
The finite Euler product in \eqref{eq:general-D-euler}
and its fixed derivatives are $O_\eps(r^{\eps/8})$ on this line.
Stirling's formula and the rapid decay of $G$ show that the
shifted integral is
\[
\ll_\eps k^{-1+\eps/4}r^{\eps/8}
\frac{\log(2kr)}{\sqrt r}\sum_{d\mid r}d^{1/2}
\ll_\eps k^{-1+\eps},
\]
uniformly for $r\leq K$ and $k\asymp K$.
The residue is \eqref{eq:general-PV-residue}.
Averaging over $k$ and using \eqref{eq:ndnz-general-balanced}
proves the proposition.
\end{proof}
\subsection{Reduction to shifted convolution sums}
\label{sec:nonzero-reduction}
Fix $0<\eps<1/2$ and suppose $K^{2\eps}\leq L\leq K$ and
$1\leq r\leq K$.
For each dyadic scale $M=2^j$, define the averaged kernel
from \eqref{TKM} by
\begin{align}\label{eq:ndnz-averaged-kernel}
\mathcal{K}_M(m,qn;c):= & ~
\sum_{k\equiv0\bmod2}\mathrm{i}^kh_{K,L}(k)V_k(nd)
\mathcal T_{k,M}(m,qn;c).
\end{align}
The additional dependence on $d,n$ through $V_k(nd)$ is left implicit.
Setting $\Delta=qn-m$ gives
\begin{equation}\label{eq:ndnz-shift-representation}
\mathcal E_{M}(r)
=16\pi^2\mathrm{i}^\ell\sum_{d\mid r}\sum_{c\geq1}\frac1{c^2}
\sum_{\Delta\in\mathbb Z\setminus\{0\}}R_c(\Delta)\sum_{\substack{n\geq1\\qn-\Delta\geq1}}
\frac{\lambda_g(nd)\lambda_g(qn-\Delta)}{\sqrt{nd}}
\mathcal{K}_M(qn-\Delta,qn;c),
\end{equation}
where $\mathcal E(r)=\sum_{j\in\mathbb Z}\mathcal E_{2^j}(r)$.

\subsubsection{Uniform truncation}

We first truncate $M$ while keeping the full $n$ sum.

\begin{lemma}\label{lem:ndnz-truncation}
For every $A>0$, we have 
\[
\mathcal E(r)=
\sum_{1\leq M \ll K^{2+\eps}}\mathcal E_{M}(r)+O(K^{-A}).
\]
\end{lemma}

\begin{proof}
All estimates below are uniform in the outer cutoffs
following \eqref{eq:balanced-Voronoi}.
For dyadic $M<1$, we have $M\leq1/2$ and $U(m/M)=0$
for every integer $m\geq1$.
The Voronoi summation formula for the original $m$-sum in
\eqref{the second moment} with this factor gives
\[
\sum_{\substack{m\geq1\\m\ne qn}}
\lambda_g(m)R_c(qn-m)\mathcal T_{k,M}(m,qn;c)
=-\lambda_g(qn)\varphi(c)\mathcal T_{k,M}(qn,qn;c).
\]
Substituting this into \eqref{eq:ndnz-shift-representation} gives
\begin{align}\label{small M}
\mathcal E_{M}(r)
=-16\pi^2\mathrm{i}^\ell
\sum_{k\equiv0\bmod2}\mathrm{i}^kh_{K,L}(k)\sum_{d\mid r}\sum_{n\geq1}
\frac{\lambda_g(nd)\lambda_g(qn)V_k(nd)}{\sqrt{nd}}
\sum_{c\geq1}\frac{\varphi(c)}{c^2}
\mathcal T_{k,M}(qn,qn;c).
\end{align}
Making the substitution $y=4\pi\sqrt{xqn}/c$ in \eqref{TKM}, we obtain
\[
\mathcal T_{k,M}(qn,qn;c)
=\int_0^\infty\phi_M(y)J_{k-1}(y)J_{\ell-1}(y)
\frac{\dd y}{y},
\]
where
\begin{align*}
\phi_M(y)&=\frac{cy}{2\pi\sqrt{qn}}
V_k\left(\frac{c^2y^2}{16\pi^2qn}\right)
U\left(\frac{c^2y^2}{16\pi^2qnM}\right).
\end{align*}
Apply \eqref{Bessel differential equation} separately to
$J_{k-1}(y)$ and $J_{\ell-1}(y)$, multiply each equation by
the other Bessel factor, and subtract. Integration by parts gives
\begin{align*}
&\bigl((k-1)^2-(\ell-1)^2\bigr)
\mathcal T_{k,M}(qn,qn;c)\\
&\qquad=-\int_0^\infty y^2\frac{\dd\phi_M(y)}{\dd y}
\bigl(J_{k-1}'(y)J_{\ell-1}(y)
-J_{k-1}(y)J_{\ell-1}'(y)\bigr)\frac{\dd y}{y}.
\end{align*}
For fixed $N\geq1$, $k\asymp K$ and fixed $\ell$ imply
\[
\left|(k-1+\rho)^2-(\ell-1+\sigma)^2\right|
\asymp_N K^2,
\qquad\text{for}\quad\rho,\sigma\in\mathbb Z\,\,\,\text{and}\,\,\, |\rho|+|\sigma|\leq N.
\]
Applying \eqref{derivative recurrence} and the product rule
at each step, then using $|J_n(y)|\leq1$, we obtain after
$N$ integrations by parts
\[
|\mathcal T_{k,M}(qn,qn;c)|
\ll_N K^{-2N}\sum_{j=1}^N
\int_0^\infty y^{N+j}
\left|\frac{\dd^j\phi_M(y)}{\dd y^j}\right|
\frac{\dd y}{y}.
\]
Since $c^2y^2/(16\pi^2qn)\asymp M$ on their support,
Lemma~\ref{cutoff of Vk} yields
\[
y^j\left|\frac{\dd^j\phi_M(y)}{\dd y^j}\right|
\ll_N\frac{cy}{\sqrt{qn}}(\log K+|\log M|),
\qquad\text{for}\quad0\leq j\leq N.
\]
Substituting this bound gives
\begin{align*}
|\mathcal T_{k,M}(qn,qn;c)|
&\ll_N K^{-2N}\frac{c(\log K+|\log M|)}{\sqrt{qn}}
\int_{4\pi\sqrt{qnM/2}/c}^{4\pi\sqrt{2qnM}/c}
y^N\,\dd y\\
&\ll_N K^{-2N}\frac{(qn)^{N/2}}{c^N}
M^{(N+1)/2}(\log K+|\log M|).
\end{align*}
Thus the $c$-sum in \eqref{small M} satisfies
\[
\sum_{c\geq1}\frac{\varphi(c)}{c^2}
\frac{(qn)^{N/2}}{c^N}
=(qn)^{N/2}
\sum_{c\geq1}\frac{\varphi(c)}{c^{N+2}}
\ll_N(qn)^{N/2}.
\]
Since $q=r/d$, partial summation and
Lemma~\ref{cutoff of Vk} give
\begin{align}\label{eq:ndnz-low-coefficient-sum}
&\sum_{d\mid r}\sum_{n\geq1}
\frac{|\lambda_g(nd)\lambda_g(qn)V_k(nd)|}{\sqrt{nd}}\,
(qn)^{N/2}\notag\\
&\quad\leq r^{N/2}
\sum_{d\mid r}\frac{\tau(d)\tau(r/d)}{d^{(N+1)/2}}
\sum_{n\geq1}n^{(N-1)/2}\tau(n)^2|V_k(nd)|\notag\\
&\quad\ll_N r^{N/2}K^{N+1}(\log K)^3
\sum_{d\mid r}\frac{\tau(d)\tau(r/d)}{d^{N+1}}\notag\\
&\quad\ll_{N,\eps} K^{3N/2+1+\eps}.
\end{align}
Since $\sum_{k\equiv0\bmod2}|h_{K,L}(k)|\ll L$,
combining \eqref{small M}--\eqref{eq:ndnz-low-coefficient-sum}, and summing over
$M=2^j<1$, we obtain
\[
\sum_{j<0}|\mathcal E_{2^j}(r)|
\ll_{N,\eps} L K^{1-N/2+\eps}.
\]
Since $L\leq K$, this is $O_{A,\eps}(K^{-A})$ for
$N>2(A+3+\eps)$.

For $M>K^{2+\eps}$ and any fixed $B>1$,
Lemma~\ref{cutoff of Vk} gives
\[
\max_{0\leq i\leq4}\sup_{x\asymp M}
x^i\left|\frac{\dd^i}{\dd x^i}
\left(V_k(x)U\left(\frac{x}{M}\right)\right)\right|
\ll_B\left(\frac{K^2}{M}\right)^B.
\]
Applying the Bessel equation \eqref{Bessel differential equation} twice to
$J_{\ell-1}(4\pi\sqrt{xm}/c)$ in \eqref{TKM} and
integrating by parts gives
\[
|\mathcal T_{k,M}(m,qn;c)|
\ll_\ell\frac{c^4}{m^2M}
\max_{0\leq i\leq4}\sup_{x\asymp M}
x^i\left|\frac{\dd^iF(x)}{\dd x^i}\right|.
\]
The product rule, \eqref{derivative recurrence}, and the
Bessel power series \cite[8.402]{GR} imply, for a fixed
sufficiently large $C$,
\begin{equation}\label{eq:ndnz-large-M-kernel}
|\mathcal T_{k,M}(m,qn;c)|
\ll_B\frac{c^4M^{-3/2}}{m^2}
\left(\frac{K^2}{M}\right)^B
\begin{cases}
\left(1+\sqrt{M qn}/c\right)^4,
&c\leq C\sqrt{M qn},\\
\left(\sqrt{M qn}/c\right)^6,
&c>C\sqrt{M qn}.
\end{cases}
\end{equation}
Summing over $m$ and $c$ gives
\begin{align}\label{eq:ndnz-large-M-dual-sum}
&\sum_{c\geq1}\frac1{c^2}
\sum_{\substack{m\geq1\\m\ne qn}}
|\lambda_g(m)R_c(qn-m)\mathcal T_{k,M}(m,qn;c)|\notag\\
&\quad\ll_B M^{-3/2}\left(\frac{K^2}{M}\right)^B
\left(
(M qn)^2\sum_{c\leq C\sqrt{M qn}}\frac1c
+(M qn)^3\sum_{c>C\sqrt{M qn}}\frac1{c^3}
\right)\notag\\
&\quad\ll_B M^{1/2}(qn)^2
\bigl(1+\log(2+M qn)\bigr)
\left(\frac{K^2}{M}\right)^B.
\end{align}
Reindexing $nd$ as $n$ and using $q=r/d$, we obtain
\begin{align}\label{eq:ndnz-large-M-coefficient-sum}
&\sum_{d\mid r}\sum_{n\geq1}
\frac{|\lambda_g(nd)V_k(nd)|}{\sqrt{nd}}\,
(qn)^2\bigl(1+\log(2+M qn)\bigr)\notag\\
&\quad=r^2\sum_{d\mid r}\frac1{d^4}
\sum_{\substack{n\geq1\\d\mid n}}
n^{3/2}|\lambda_g(n)V_k(n)|
\left(1+\log\left(2+\frac{Mrn}{d^2}\right)\right)\notag\\
&\quad\ll K^7\log K\log(2MK).
\end{align}
Combining \eqref{eq:ndnz-shift-representation},
\eqref{eq:ndnz-large-M-dual-sum}, and
\eqref{eq:ndnz-large-M-coefficient-sum} gives
\[
|\mathcal E_{M}(r)|
\ll_B K^8M^{1/2}\log K\log(2MK)
\left(\frac{K^2}{M}\right)^B,
\]
and hence
\begin{align*}
\sum_{2^j>K^{2+\eps}}|\mathcal E_{2^j}(r)|
&\ll_B K^8\log K\sum_{2^j>K^{2+\eps}}
2^{j/2}\log(2^{j+1}K)\left(\frac{K^2}{2^j}\right)^B\\
&\ll_{B,\eps} K^{9+\eps/2-\eps B}(\log K)^2
\ll_{B,\eps} K^{9+\eps-\eps B}.
\end{align*}
Taking $B>(A+10+\eps)/\eps$ completes the proof.
\end{proof}

Apply \eqref{eq:ndnz-dyadic-partition} to the $n$-sum in
\eqref{eq:ndnz-shift-representation}, with weight $U(n/N)$.
By Lemma~\ref{cutoff of Vk} and the estimates used in
\eqref{eq:ndnz-large-M-dual-sum} and
\eqref{eq:ndnz-large-M-coefficient-sum}, the parts
$dN>K^{2+\eps}$ contribute $O_{A,\eps}(K^{-A})$.
Together with Lemma~\ref{lem:ndnz-truncation}, this gives
\begin{align}
1\leq M\ll K^{2+\eps},\qquad
1\leq N\ll K^{2+\eps}/d,
\label{eq:ndnz-original-lengths}
\end{align}
and
\begin{align}
nd\ll K^{2+\eps},\qquad
qn=\frac{r}{d^2}(nd)\ll\frac{rK^{2+\eps}}{d^2}\ll K^{3+\eps}.
\label{eq:ndnz-truncated-original-ranges}
\end{align}

We next truncate $c$ and the dual variable $m$.

\begin{lemma}\label{lem:ndnz-polynomial-truncation}
Under \eqref{eq:ndnz-original-lengths} and
\eqref{eq:ndnz-truncated-original-ranges}, one may restrict
\eqref{eq:ndnz-shift-representation} to
\begin{align*}
c\leq K^4,\qquad m=qn-\Delta\leq K^{B(A)},
\end{align*}
where $B(A)=2(A+30)$, with total error $O_A(K^{-A})$.

\end{lemma}
\begin{proof}
By Lemma~\ref{cutoff of Vk}, the argument proving
\eqref{eq:ndnz-large-M-kernel} applies in the ranges
\eqref{eq:ndnz-original-lengths}, with $(K^2/M)^B$
replaced by $\log K$.
Lemma~\ref{lem:ndnz-truncation} and
\eqref{eq:ndnz-truncated-original-ranges} give
$Mqn\ll K^{5+2\eps}$.
For $c>K^4$, taking any fixed $B>4$ in the small argument
estimate used for \eqref{eq:ndnz-large-M-kernel} yields
\begin{align*}
&\sum_{c>K^4}\frac1{c^2}
\sum_{\substack{m\geq1\\m\ne qn}}
|\lambda_g(m)R_c(qn-m)\mathcal T_{k,M}(m,qn;c)|\notag\\
&\quad\ll_B M^{-3/2}\log K\,(M qn)^{B/2}K^{16-4B}.
\end{align*}
By \eqref{eq:ndnz-shift-representation},
\eqref{eq:ndnz-large-M-coefficient-sum}, and the above inequality,
the total contribution of
$c>K^4$ is $O_B(K^{25-(3/2-\eps)B})=O_A(K^{-A})$
when $B>A+26$.

For $c\leq K^4$, \eqref{eq:ndnz-large-M-kernel} and
\eqref{eq:ndnz-truncated-original-ranges} give
\begin{equation*}
|\mathcal T_{k,M}(m,qn;c)|
\ll K^{17}m^{-2}.
\end{equation*}
Hence
\bea\label{eq:ndnz-dual-sum-tail}
\sum_{c\leq K^4}\frac1{c^2}
\sum_{\substack{m>K^{B(A)}\\m\ne qn}}
|\lambda_g(m)R_c(qn-m)\mathcal T_{k,M}(m,qn;c)|
\ll_A K^{17-B(A)}(\log K)^2.
\eea
Combining \eqref{eq:ndnz-shift-representation},
\eqref{eq:ndnz-large-M-coefficient-sum}, and
\eqref{eq:ndnz-dual-sum-tail}, the total contribution of
$m>K^{B(A)}$ is $O_A(K^{-A})$.
\end{proof}

We impose \eqref{eq:ndnz-original-lengths},
\eqref{eq:ndnz-truncated-original-ranges}, and
Lemma~\ref{lem:ndnz-polynomial-truncation} below.

\subsubsection{Effective ranges}

The factor $U(x/M)$ in \eqref{TKM} restricts $x$ to $[M/2,2M]$,
so we write
\begin{equation}\label{eq:ndnz-bessel-scale}
Y:=\frac{4\pi\sqrt{M qn}}c,\qquad
y=\frac{4\pi\sqrt{x qn}}c\asymp Y.
\end{equation}

\begin{lemma}
\label{lem:ndnz-localization}
Fix $B>0$ and $\eps>0$. Suppose that $K^{2\eps}\leq L\leq K$ and
\[
1\leq M,m,qn,c,nd\leq K^B.
\]
For every $A>0$, the kernel $\mathcal{K}_M(m,qn;c)$ in \eqref{eq:ndnz-averaged-kernel} is
$O_{A,B,\eps}(K^{-A})$ unless
\begin{align}\label{eq:ndnz-y-range}
Y\geq K^{1-\eps}L.
\end{align}
In this range it is also $O_{A,B,\eps}(K^{-A})$ unless
\begin{align}\label{eq:ndnz-shift-window}
m\asymp qn,\qquad |m-qn|\ll K^\eps
\left(\frac{qn}{Y}+\frac{qn K^2}{Y^2}\right).
\end{align}
\end{lemma}

\begin{proof}
For $x\in[M/2,2M]$, put
\begin{equation*}
H_{x,nd}(v):=h(v)V_{K+1+Lv}(x)V_{K+1+Lv}(nd),
\end{equation*}
extended by zero outside the support of $h$.
Substituting \eqref{TKM} into \eqref{eq:ndnz-averaged-kernel}
gives the sum
\begin{equation}\label{eq:ndnz-weighted-bessel-sum}
\sum_{k\equiv0\bmod2}\mathrm{i}^k
H_{x,nd}\left(\frac{k-1-K}{L}\right)J_{k-1}(y)
\end{equation}
inside the $x$ integral.
The spectral weight estimates in the proof of
Theorem~\ref{thm:fmcor} apply to $H_{x,nd}(v)$
with $(\log K)^2$ in place of $\log K$.
Lemma~\ref{lem:exact-bessel-average} and the averaging argument
therefore make \eqref{eq:ndnz-weighted-bessel-sum}
$O(K^{-A})$ for $y\ll K^{1-\eps}L$.
After the $x$-integration, this proves \eqref{eq:ndnz-y-range}.

Assume now $Y\geq K^{1-\eps}L\geq K^{1+\eps}$.
Put $\kappa=k-1\asymp K$ and $\beta=\sqrt{m/(qn)}$.
The integral representation \cite[pp.~19--21]{Watson}
has stationary points $t=\pm\arccos(\kappa/y)$, where the second
derivative of its phase has magnitude $\sqrt{y^2-\kappa^2}\asymp Y$.
Since $\kappa/y\ll K^{-\eps}$, the stationary phase gives uniformly
\begin{align}\label{eq:ndnz-large-order-bessel}
J_\kappa(y)
=\frac1{2\pi}\int_{-\pi}^{\pi}e^{\mathrm{i}(y\sin t-\kappa t)}\,\dd t
=\frac{e^{\mathrm{i}\theta_\kappa(y)}+e^{-\mathrm{i}\theta_\kappa(y)}}
{\sqrt{2\pi}(y^2-\kappa^2)^{1/4}}+O(y^{-3/2}),
\end{align}
where
\begin{equation*}
\theta_\kappa(y):=
\sqrt{y^2-\kappa^2}-\kappa\arccos(\kappa/y)-\frac\pi4.
\end{equation*}
For $\beta y\gg1$, the expansion for fixed $\ell$ in
\cite[pp.~237--242]{O} gives
\begin{equation}\label{eq:ndnz-fixed-order-bessel}
J_{\ell-1}(\beta y)
=\frac{e^{\mathrm{i}\beta y-\mathrm{i}\pi(2\ell-1)/4}
+e^{-\mathrm{i}\beta y+\mathrm{i}\pi(2\ell-1)/4}}{\sqrt{2\pi\beta y}}
+O_\ell((\beta y)^{-3/2}).
\end{equation}
Absorbing the constant phases into the amplitudes, the product
in \eqref{TKM} has phases
\begin{equation*}
\pm\bigl(\theta_\kappa(y)+\beta y\bigr),\qquad
\pm\bigl(\theta_\kappa(y)-\beta y\bigr).
\end{equation*}
Only the leading terms are displayed.
We use expansions to arbitrary order, with fixed normalized
amplitude derivatives of sizes $O(Y^{-1/2})$ and
$O_\ell((\beta Y)^{-1/2})$, respectively.
By Lemma~\ref{cutoff of Vk}, every fixed normalized derivative of
$V_k(c^2y^2/(16\pi^2qn))U(y^2/Y^2)$ is $O_j(\log K)$.
Thus the full weights satisfy the derivative conditions
needed for Lemma~\ref{lemma:exponentialintegral}.

First suppose that $m<qn/4$ or $m>4qn$, so that
$\beta<1/2$ or $\beta>2$.
If $\beta Y\leq K^{\eps/10}$, keep $J_{\ell-1}(\beta y)$
in the amplitude. By \eqref{derivative recurrence}, its
$j$-th normalized derivative is $O_j(K^{j\eps/10})$,
whereas $|y\theta_\kappa'(y)|\asymp Y$.
Lemma~\ref{lemma:exponentialintegral} makes these terms negligible.
Otherwise apply \eqref{eq:ndnz-fixed-order-bessel}. For either sign,
\begin{align}\label{eq:ndnz-separated-phase-derivatives}
\left|y\frac{\dd}{\dd y}
\bigl(\theta_\kappa(y)\pm\beta y\bigr)\right|
=\left|\sqrt{y^2-\kappa^2}\pm\beta y\right|\gg Y\max(1,\beta),
\qquad\text{for}\,\,\,\beta<1/2\,\,\, \text{or}\,\,\, \beta>2.
\end{align}
For the minus sign, the difference is at least $y/4$ when
$\beta<1/2$, since $\sqrt{y^2-\kappa^2}\geq3y/4$ for sufficiently
large $K$. When $\beta>2$, its absolute value is at least
$(\beta-1)y\geq\beta y/2$.
For the plus sign, the lower bound holds for every $\beta>0$.
The higher derivatives satisfy
\begin{equation}\label{eq:ndnz-bessel-phase-higher-derivatives}
y^j\left|\frac{\dd^j}{\dd y^j}
\bigl(\theta_\kappa(y)\pm\beta y\bigr)\right|
\ll_j Y\max(1,\beta),\qquad\text{for}\,\,\,j\geq2.
\end{equation}
Lemma~\ref{lemma:exponentialintegral} again gives negligible terms.

For $qn/4\leq m\leq4qn$, use only the plus sign in
\eqref{eq:ndnz-separated-phase-derivatives}, together with
\eqref{eq:ndnz-bessel-phase-higher-derivatives}, to discard
the phases $\pm(\theta_\kappa(y)+\beta y)$.
For the difference $\Phi_\kappa(y):=\theta_\kappa(y)-\beta y$, we have
\begin{equation}\label{eq:ndnz-phase-derivative}
y\Phi_\kappa'(y)
=\sqrt{y^2-\kappa^2}-\beta y
=\frac{(qn-m)y^2/(qn)-\kappa^2}
{\sqrt{y^2-\kappa^2}+\beta y},
\end{equation}
and differentiation gives
\begin{equation}\label{eq:ndnz-phase-higher-derivatives}
y^j\left|\frac{\dd^j\Phi_\kappa(y)}{\dd y^j}\right|
\ll_j |1-\beta|Y+\frac{K^2}{Y}\qquad(j\geq1).
\end{equation}
If the bound on $|m-qn|$ in \eqref{eq:ndnz-shift-window} fails, then
\[
|1-\beta|Y=\frac{|m-qn|Y}{qn(1+\beta)}
\asymp\frac{|m-qn|Y}{qn}
\gg K^\eps(1+K^2/Y).
\]
\eqref{eq:ndnz-phase-derivative} then gives
$|y\Phi_\kappa'(y)|\asymp|1-\beta|Y$, so
\eqref{eq:ndnz-phase-higher-derivatives} and
Lemma~\ref{lemma:exponentialintegral} give the required decay.

Finally, $Y\geq K^{1+\eps}$, and $\beta Y>K^{\eps/10}$
whenever \eqref{eq:ndnz-fixed-order-bessel} is used.
Choosing the expansion and integration orders in terms of
$A,B,\eps$ absorbs the polynomial losses from integration and summation.
\end{proof}

\begin{corollary}\label{cor:ndnz-effective-ranges}
Assume $K^{2\eps}\leq L\leq K$.
For every $A>0$, the sum for $\mathcal E(r)$ may be
restricted, with total error $O_A(K^{-A})$, to
\eqref{eq:ndnz-original-lengths}, \eqref{eq:ndnz-y-range},
\eqref{eq:ndnz-shift-window}, and
\begin{align}\label{eq:ndnz-effective-ranges}
c\ll\frac{K^{1+2\eps}\sqrt r}{dL},
\qquad
0<|\Delta|\ll\frac{qn K^{3\eps}}{L^2}.
\end{align}
\end{corollary}

\begin{proof}
The ranges \eqref{eq:ndnz-original-lengths} and
\eqref{eq:ndnz-truncated-original-ranges}, together with
Lemma~\ref{lem:ndnz-polynomial-truncation}, put all variables
in fixed polynomial ranges. Lemma~\ref{lem:ndnz-localization}
therefore applies after summation.
Equations \eqref{eq:ndnz-truncated-original-ranges}, \eqref{eq:ndnz-bessel-scale}, and
\eqref{eq:ndnz-y-range} give the bound for $c$ in
\eqref{eq:ndnz-effective-ranges}.
The bound for $\Delta$ follows from
\eqref{eq:ndnz-shift-window} and
\[
K^\eps\left(\frac{qn}{Y}+\frac{qn K^2}{Y^2}\right)
\ll\frac{qn K^{3\eps}}{L^2}.
\]
On each block $n\asymp N$, use $qN$ in place of $qn$ in the
bounds for $c$ and $\Delta$, with constants large enough to
include all effective points. We then sum over the full support
of $U(n/N)$. By Lemma~\ref{lem:ndnz-localization}, the added terms
contribute $O_A(K^{-A})$ in total.
\end{proof}

\subsubsection{Application of the Hecke relation}
For $\Delta\neq0$, the Ramanujan sum has the exact divisor expansion
\begin{align}\label{eq:ndnz-ramanujan}
R_c(\Delta)= & ~\sum_{v\mid(c,\Delta)}v\mu(c/v).
\end{align}
In the $n$-sum of \eqref{eq:ndnz-shift-representation}, the Hecke relation \eqref{Hecke relation} gives
\[
 \lambda_g(nd)=\sum_{b\mid(n,d)}\mu(b)\lambda_g(n/b)\lambda_g(d/b).
\]
Writing $n=bt$ on each dyadic block, we obtain, for every $A>0$,
\begin{align}\label{eq:ndnz-hecke-error-representation}
\mathcal E(r)
={}&16\pi^2\mathrm{i}^\ell
\sum_{d\mid r}\sum_{b\mid d}\mu(b)\lambda_g(d/b)
\sum_{M,N}\sum_{c\geq1}\frac1{c^2}\notag\\
&\times\sum_{\Delta\neq0}R_c(\Delta)
\sum_{\substack{t,m\geq1\\qbt-m=\Delta}}
\frac{\lambda_g(t)\lambda_g(m)U(bt/N)}{\sqrt{bdt}}\,
\mathcal{K}_M(m,qbt;c)
+O_{A,\eps}(K^{-A}).
\end{align}
Here $q=r/d$, and $M,N$ run over the dyadic ranges in
 \eqref{eq:ndnz-original-lengths}.
\subsection{The spectral estimate}
Under the restrictions in Corollary~\ref{cor:ndnz-effective-ranges},
we estimate $\mathcal E(r)$ in \eqref{eq:ndnz-hecke-error-representation}
by combining the sum over $\Delta=qn-m$ with the spectral expansion
of the shifted convolution sums, following
\cite[\S\S\,5--9]{LLY2006IMRP}.
We take $\theta=7/64$.

\begin{proposition}\label{v6:prop:joint-error}
Fix $1/3<\delta\leq1$. Uniformly for $K^\delta\leq L\leq K$ and
all integers $1\leq r\leq K$, we have
\begin{align}\label{v6:eq:joint-error}
\mathcal E(r)\ll & ~\frac L{\sqrt r}K^\eps
\left(
r^{3/2}\left(\frac K{L^3}\right)^{1/2}
+r^{87/64}K^{-25/32}
\right).
\end{align}
\end{proposition}

We split the dyadic blocks according to the size of $T$ in
\eqref{eq:joint-TH-definitions} below. Throughout this subsection, we use $\varepsilon_2$ as the exponent in \eqref{eq:ndnz-original-lengths} and in
Corollary~\ref{cor:ndnz-effective-ranges}.

\subsubsection{Spectral expansion and coefficient bounds}
For a positive integer $\mathfrak{q}$, put
\[
U_{\mathfrak{q}}(z)=\mathfrak{q}^{(\ell-1)/2}y^\ell
g(\mathfrak{q} z)\overline{g(z)}.
\]
This function is $\Gamma_0(\mathfrak{q})$ invariant and square integrable.
We use the inner product
\[
\langle F_1,F_2\rangle
=\int_{\Gamma_0(\mathfrak{q})\backslash\mathbb H}
F_1(z)\overline{F_2(z)}\,\frac{\dd x\,\dd y}{y^2}.
\]
Let $\{u_\nu\}$ be an orthonormal basis of Maass cusp forms of
level $\mathfrak{q}$ and trivial character, including oldforms, with
Laplace eigenvalues $1/4+t_\nu^2$.
For each singular cusp $\mathfrak a$, let
$E_{\mathfrak a}(z,s)$ be the Eisenstein series normalized by the
term $y^s$ at $\mathfrak a$.
We use the Fourier expansions
\begin{align*}
u_\nu(x+\mathrm{i}y)= & ~\sqrt y\sum_{n\neq0}\rho_\nu(n)
K_{\mathrm{i}t_\nu}(2\pi|n|y)e(nx),
\end{align*}
and 
\begin{align*}
E_{\mathfrak a}(x+\mathrm{i}y,1/2+\mathrm{i}t)= & ~\text{constant terms}
+\sqrt y\sum_{n\neq0}\rho_{\mathfrak a}(n,t)
K_{\mathrm{i}t}(2\pi|n|y)e(nx).
\end{align*}
Let $\mathfrak S_{\mathfrak{q}}$ consist of the cusp forms with
$t_\nu\in\mathbb R$ and the Eisenstein spectrum, with measure
\[
\int_{\mathfrak S_{\mathfrak{q}}}f(\omega)\,\dd\mu(\omega)
=\sum_{t_\nu\in\mathbb R}f(u_\nu)
+\frac1{4\pi}\sum_{\mathfrak a}\int_{\mathbb R}
f(E_{\mathfrak a}(\cdot,1/2+\mathrm{i}t))\,\dd t.
\]
For each $\omega\in\mathfrak S_{\mathfrak{q}}$, write $t_\omega$
and $\rho_\omega(n)$ for its spectral parameter and Fourier
coefficients.
The exceptional cusp forms have $t_\nu=\mathrm{i}\beta_\nu$ with
$0<\beta_\nu\leq\theta$ and are included separately below.

\begin{lemma}\label{lem:joint-spectral-inputs}
Uniformly for integers $\mathfrak{q}\geq1$ and real $T\geq1$,
\begin{align}\label{v6:eq:uniform-Good}
\int_{\substack{\omega\in\mathfrak S_{\mathfrak{q}}\\|t_\omega|\leq T}}
|\langle U_{\mathfrak{q}},\omega\rangle|^2
e^{\pi|t_\omega|}\,\dd\mu(\omega)
\ll & ~\mathfrak{q} T^{2\ell}.
\end{align}
For $H\geq1$ and complex coefficients $a_n$ supported on
$H<|n|\leq2H$,
\begin{align}\label{v6:eq:uniform-large-sieve}
&\int_{\substack{\omega\in\mathfrak S_{\mathfrak{q}}\\|t_\omega|\leq T}}
\frac{1}{\cosh(\pi t_\omega)}
\biggl|\sum_{H<|n|\leq2H}a_n\rho_\omega(n)\biggr|^2\,\dd\mu(\omega)
+\sum_{t_\nu=\mathrm{i}\beta_\nu}
\frac{1}{\cos(\pi\beta_\nu)}
\biggl|\sum_{H<|n|\leq2H}a_n\rho_\nu(n)\biggr|^2\notag\\
&\qquad\ll (\mathfrak{q}HT)^\eps
\biggl(T^2+\frac{H}{\mathfrak{q}}\biggr)
\sum_{H<|n|\leq2H}|a_n|^2.
\end{align}
The exceptional projections satisfy
\begin{align}\label{v6:eq:exceptional-projections}
\sum_{t_\nu=\mathrm{i}\beta_\nu}
|\langle U_{\mathfrak{q}},u_\nu\rangle|^2
\leq\|U_{\mathfrak{q}}\|_2^2\ll & ~\mathfrak{q}^\eps.
\end{align}
For $\Delta\in\mathbb Z\setminus\{0\}$ and $\Re(s)>1$, define
\begin{equation}\label{eq:shifted-dirichlet-definition}
D_g(s;\mathfrak{q},1,\Delta)
=\sum_{\substack{u,v\geq1\\\mathfrak{q} u-v=\Delta}}
\lambda_g(u)\lambda_g(v)
\left(\frac{\sqrt{\mathfrak{q} uv}}{\mathfrak{q} u+v}\right)^{\ell-1}
(\mathfrak{q} u+v)^{-s}.
\end{equation}
Its spectral expansion is
\begin{align}\label{v6:eq:shifted-spectral-expansion}
D_g(s;\mathfrak{q},1,\Delta)
={}&|\Delta|^{1/2-s}\biggl(
\int_{\mathfrak S_{\mathfrak{q}}}
\mathcal B_\ell(s,t_\omega)\rho_\omega(\Delta)
\langle U_{\mathfrak{q}},\omega\rangle\,
\dd\mu(\omega)+
\sum_{t_\nu=\mathrm{i}\beta_\nu}
\mathcal B_\ell(s,t_\nu)\rho_\nu(\Delta)
\langle U_{\mathfrak{q}},u_\nu\rangle\biggr),
\end{align}
where
\begin{align}\label{v6:eq:continuous-Mellin-kernel}
\mathcal B_\ell(s,t):= & ~
\frac{2^{s+\ell-3}\pi^{\ell-1/2}}{\Gamma(s+\ell-1)}
\Gamma\left(\frac{s-1/2+\mathrm{i}t}{2}\right)
\Gamma\left(\frac{s-1/2-\mathrm{i}t}{2}\right).
\end{align}
\end{lemma}

\begin{proof}
In \cite[Eq.\,(4.24)]{Michel2004}, the product is
\[
(\mathfrak{q} y)^{\ell/2}g(\mathfrak{q} z)\,
y^{\ell/2}\overline{g(z)}
=\mathfrak{q}^{1/2}U_{\mathfrak{q}}(z).
\]
Its squared spectral bound is $O_g(\mathfrak{q}^2T^{2\ell})$.
Dividing by $\mathfrak{q}$ gives \eqref{v6:eq:uniform-Good}.

The Whittaker function used in
\cite[Eq.\,(2.1)]{BM2015} satisfies
\[
W_{0,\mathrm{i}t}(4\pi|n|y)
=2\sqrt{|n|y}\,K_{\mathrm{i}t}(2\pi|n|y).
\]
Since the large sieve in \cite[Lemma~11]{BM2015} weights each
Fourier coefficient by $\sqrt{|n|}$, its coefficient sums are
one half of the corresponding sums in our normalization.
Applying that lemma separately to positive and negative indices
and combining the bounds by Cauchy's inequality gives
\eqref{v6:eq:uniform-large-sieve}, including the exceptional spectrum.
No condition $(n,\mathfrak{q})=1$ is required.
The endpoint version $H\leq|n|\leq2H$ of the same lemma
covers $\Delta=\pm1$ by taking $H=1$.

Since $y^{\ell/2}|g(z)|$ is bounded on $\mathbb H$,
\[
|U_{\mathfrak{q}}(z)|
=\mathfrak{q}^{-1/2}
\big|(\mathfrak{q} y)^{\ell/2}g(\mathfrak{q} z)\big|\,
\big|y^{\ell/2}g(z)\big|
\ll_g\mathfrak{q}^{-1/2}.
\]
Consequently,
\[
\|U_{\mathfrak{q}}\|_2^2
\ll_g\mathfrak{q}^{-1}
\operatorname{vol}(\Gamma_0(\mathfrak{q})\backslash\mathbb H)
=\frac\pi3\prod_{p\mid\mathfrak{q}}(1+p^{-1})
\ll_\eps\mathfrak{q}^\eps.
\]
Bessel's inequality proves \eqref{v6:eq:exceptional-projections}.

For the spectral expansion, we use the argument of
\cite[\S\,3]{LLY2006IMRP} in the Fourier normalization above.
The Fourier coefficient of $U_{\mathfrak{q}}$ at index $\Delta$ is
\begin{align*}
\int_0^1U_{\mathfrak{q}}(x+\mathrm{i}y)e(-\Delta x)\,\dd x
={}&y^\ell
\sum_{\substack{u,v\geq1\\\mathfrak{q} u-v=\Delta}}
\lambda_g(u)\lambda_g(v)(\mathfrak{q} uv)^{(\ell-1)/2}e^{-2\pi(\mathfrak{q} u+v)y}.
\end{align*}
Hence, for $\Re(s)>1$,
\begin{align*}
D_g(s;\mathfrak{q},1,\Delta)
={}&\frac{(2\pi)^{s+\ell-1}}{\Gamma(s+\ell-1)}
\int_0^\infty
\left(\int_0^1U_{\mathfrak{q}}(x+\mathrm{i}y)e(-\Delta x)\,\dd x\right)
y^{s-2}\,\dd y.
\end{align*}
Insert the spectral decomposition of $U_{\mathfrak{q}}$.
First insert a smooth compactly supported cutoff in $y$ and apply
Parseval to the Poincar\'e series with seed
$y^{\overline s}e(\Delta x)$ times this cutoff. The bound
\[
|K_{\mathrm{i}t}(x)|\leq e^{-\pi|t|/4}K_0(x/\sqrt2),
\qquad\text{for}\,\,\, x>0,\ t\in\mathbb R,
\]
together with Cauchy's inequality, \eqref{v6:eq:uniform-Good}, and
the single coefficient large sieve above, permits removal of the
cutoff for $\Re(s)>1$. The finite exceptional sum is integrable there.
The constant eigenfunction contributes zero since $\Delta\neq0$.
For every remaining component, the Mellin integral is
\begin{align*}
\int_0^\infty y^{s-3/2}K_{\mathrm{i}t}(2\pi|\Delta|y)\,\dd y
=2^{s-5/2}(2\pi|\Delta|)^{1/2-s}
\Gamma\left(\frac{s-1/2+\mathrm{i}t}{2}\right)
\Gamma\left(\frac{s-1/2-\mathrm{i}t}{2}\right).
\end{align*}
This gives \eqref{v6:eq:shifted-spectral-expansion} with
\eqref{v6:eq:continuous-Mellin-kernel} for both the cuspidal and
continuous spectra.
\end{proof}
\subsubsection{The dyadic kernel}
Recall that $q=r/d$ and that $n=bt$. Fix $d\mid r$, $b\mid d$, and a
dyadic block $x\asymp M$, $n\asymp N$ in
\eqref{eq:ndnz-hecke-error-representation}.
Take $\mathfrak{q}=qb$ in Lemma~\ref{lem:joint-spectral-inputs}.
By Corollary~\ref{cor:ndnz-effective-ranges}, the relation is
$\mathfrak{q}t-m=\Delta$, with $t\asymp N/b$ and
$\mathfrak{q}t,m\asymp qN$. By \eqref{eq:ndnz-original-lengths},
\begin{equation}\label{v6:eq:balanced-new-ranges}
qN\ll\frac{rK^{2+\varepsilon_2}}{d^2}.
\end{equation}

In \eqref{TKM}, use the normalized variables
\[
z:=\frac nN=\frac{bt}{N},\qquad
w:=\sqrt{\frac{xn}{MN}},\qquad
c_*:=\sqrt{\frac{qN}{M}}\,c.
\]
Then $z,w\asymp1$. For fixed $n$, we have $x=Mw^2/z$, and
\begin{equation*}
\frac{\dd x}{c^2\sqrt x}
=\frac{2qN}{c_*^2\sqrt{Mz}}\,\dd w,\qquad
c=c_*\sqrt{M/(qN)}\in\mathbb N.
\end{equation*}

For a dyadic block $c_*\asymp C_*$, define
\begin{equation}\label{eq:joint-TH-definitions}
T:=\frac{K^2C_*}{qN},\qquad
H:=C_*T=\frac{K^2C_*^2}{qN}.
\end{equation}
Since $qn\asymp qN$, \eqref{eq:ndnz-bessel-scale} and
\eqref{eq:ndnz-y-range} give
\begin{align}\label{v6:eq:joint-dyadic-ranges}
C_*\ll\frac{qN}{K^{1-\varepsilon_2}L},\qquad
T\ll\frac{K^{1+\varepsilon_2}}L.
\end{align}
Consequently,
\[
\frac{H}{qN}=\frac{T^2}{K^2}\ll
\frac{K^{2\varepsilon_2}}{L^2}=o(1).
\]
We first treat $T\geq K^{\varepsilon_0}$
for a fixed $\varepsilon_0>0$.
\begin{lemma}\label{lem:joint-dyadic-kernel}
Suppose that $K^\delta\leq L\leq K^{1-\varepsilon_1}$, with
$\delta>1/3$, and that $T\geq K^{\varepsilon_0}$.
Up to a total error $O(K^{-A})$ for every fixed $A>0$, each
dyadic block in \eqref{eq:ndnz-hecke-error-representation} is
a sum over $\eta=\pm1$, restricted to $\Delta>0$ and
$\Delta\asymp H$.
After the change of variables above and the Bessel expansions,
define the normalized kernels by
\begin{align}\label{v6:eq:normalized-W}
W_{\Delta,c_*}^{\eta}(z)
:={}&\int_0^\infty A_{\Delta,c_*}^{\eta}(z,w)
e\left(\eta\left(
\frac{2qNw}{c_*}\left(1-\sqrt{1-\frac{\Delta}{qNz}}\right)
+\frac{K^2c_*}{16\pi^2qNw}\right)\right)\dd w.
\end{align}
Here $A_{\Delta,c_*}^{\eta}(z,w)$ is a smooth complex-valued
amplitude, independent of $s$. It is obtained by averaging over $k$
the product of the cutoffs, the $V_k$ weights, the retained
Bessel amplitudes and the remaining exponential factor.
Its support is contained in
\[
\frac12\leq z\leq2,\qquad
\sqrt{z/2}\leq w\leq\sqrt{2z}.
\]
For each prescribed integer $J\geq0$, the auxiliary exponents
can be chosen sufficiently small so that
\begin{equation}\label{eq:joint-amplitude-derivatives}
\sup_{z,w>0}
\left|\frac{\partial^{j_1+j_2}A_{\Delta,c_*}^{\eta}(z,w)}
{\partial z^{j_1}\partial w^{j_2}}\right|
\ll_{J,\eps}K^\eps,
\qquad\text{for}\quad j_1,j_2\geq0,\ j_1+j_2\leq J.
\end{equation}
For $n=Nz=bt$, these kernels represent the smooth weight in
\eqref{eq:ndnz-hecke-error-representation},
\begin{align}\label{eq:joint-original-weight-to-W}
\frac{16\pi^2\mathrm{i}^\ell U(n/N)}{c^2\sqrt{nd}}\,
\mathcal{K}_M(qn-\Delta,qn;c)=
\frac L{\sqrt r}\frac{q}{c_*\sqrt{MqN}}
\sum_{\eta=\pm1}W_{\Delta,c_*}^{\eta}(n/N)
+O_A(K^{-A}).
\end{align}
Fix $\eta$ in the remaining formulas and suppress its dependence
in the transforms. Define
\begin{align}\label{v6:eq:explicit-Gtilde}
\widetilde G_{\Delta,c_*}(s)
:={}&\frac K{\sqrt{qN}}(2qN)^s
\int_{1/2}^{2} W_{\Delta,c_*}^{\eta}(z)
\left(\frac{\sqrt{z(z-\Delta/(qN))}}{2z-\Delta/(qN)}\right)^{1-\ell}
\left(z-\frac{\Delta}{2qN}\right)^{s-1}\dd z
\end{align}
and
\begin{align}\label{v6:eq:explicit-P}
\mathcal P(C_*,T)
:={}&\sum_{c_*\asymp C_*}\frac1{c_*}
\sum_{\Delta\asymp H}R_c(\Delta)
\frac1{2\pi \mathrm{i}}\int_{(2)}
D_g(s;\mathfrak{q},1,\Delta)
\widetilde G_{\Delta,c_*}(s)\dd s,
\end{align}
where $D_g$ is the shifted convolution Dirichlet series
in \eqref{eq:shifted-dirichlet-definition}, with
$\mathfrak{q}u-v=\Delta$.
For each $\eta$, the Mellin inversion gives
\begin{align}\label{v6:eq:block-normalization}
\frac{Lq\mu(b)\lambda_g(d/b)}{\sqrt{rMqN}}
\sum_{c_*\asymp C_*}\frac1{c_*}
\sum_{\Delta\asymp H}R_c(\Delta)
\sum_{\substack{t,m\geq1\\\mathfrak{q}t-m=\Delta}}
\lambda_g(t)\lambda_g(m)W_{\Delta,c_*}^{\eta}(bt/N)
=\frac{Lq\mu(b)\lambda_g(d/b)}{K\sqrt{rM}}\,
\mathcal P(C_*,T).
\end{align}
Summing \eqref{v6:eq:block-normalization} over $\eta$ and using
\eqref{eq:joint-original-weight-to-W} gives the contribution of
this dyadic block to \eqref{eq:ndnz-hecke-error-representation},
up to $O_A(K^{-A})$.
Both $\widetilde G_{\Delta,c_*}$ and the normalized transform
\begin{equation}\label{v6:eq:uniform-transform-class}
\mathcal G_{\Delta,c_*}(s):=\Delta^{1-s}\widetilde G_{\Delta,c_*}(s)
\end{equation}
are entire. The normalized transform satisfies the exact identity
\begin{align}\label{eq:joint-exact-transform}
\mathcal G_{\Delta,c_*}(s)
={}&\frac H{\sqrt{c_*}}\left(\frac{2qN}{H}\right)^s
\int_{H/(2\Delta)}^{2H/\Delta}
F'_{\Delta,c_*}(\zeta)e_\eta(\zeta,s)\dd\zeta,
\end{align}
where $F_{\Delta,c_*}$ is independent of $s$, smooth and supported in
$[H/(2\Delta),2H/\Delta]$. For each prescribed order its derivatives
are $O(K^\eps)$, in particular
\begin{equation}\label{eq:joint-F-derivative-bound}
\|F'_{\Delta,c_*}\|_\infty+\|F'_{\Delta,c_*}\|_1\ll K^\eps.
\end{equation}
Choose a fixed $\zeta_0>0$ to the left of all these supports and put
\begin{align*}
e_\eta(\zeta,s)
:=\int_{\zeta_0}^{\zeta}
\left(u-\frac H{2qN}\right)^{s-1}e\bigl(\eta\Psi_H(u)\bigr)\dd u,
\end{align*}
where
\begin{align*}
\Psi_H(u):=\frac{T}{\pi\sqrt{2u}}
\left(1+\sqrt{1-\frac{H}{qNu}}\right)^{-1/2}.
\end{align*}
For sufficiently large $K$, we have $H/(qN)<\zeta_0$.
The function $e_\eta$ is independent of $\Delta$ and $c_*$, and
\begin{equation}\label{v6:eq:primitive-bound}
e_\eta(\zeta,s)\ll_C T^{-1/2},
\qquad \text{for}\,\,\,|\Re(s)|\leq C,\ |\Im(s)|\asymp T
\end{equation}
for each fixed $C>0$, uniformly in $\zeta$ on these supports.
There are fixed $0<c_-<c_+$ such that, for
$|\Im(s)|\notin[c_-T,c_+T]$ and any fixed $A,B,C>0$, we have
\[
|\widetilde G_{\Delta,c_*}(s)|+|\mathcal G_{\Delta,c_*}(s)|
\ll_{A,B,C} K^{-A}(1+|\Im(s)|)^{-B},
\qquad \text{for}\,\,\,|\Re(s)|\leq C.
\]
The tails remain negligible after summing the polynomial parameter ranges.
\end{lemma}

\begin{proof}
We use the Bessel expansions
\eqref{eq:ndnz-large-order-bessel} and
\eqref{eq:ndnz-fixed-order-bessel}, followed by the stationary
phase argument of \cite[\S\,5]{LLY2006IMRP}.

Use $\kappa=k-1$ and $\beta=\sqrt{m/(qn)}$ as in
Lemma~\ref{lem:ndnz-localization}. With $n=Nz$, $m=qn-\Delta$ and
$x=Mw^2/z$, the argument $y$ in \eqref{eq:ndnz-bessel-scale} becomes
\[
y=\frac{4\pi qNw}{c_*},
\qquad
\beta=\sqrt{1-\frac{\Delta}{qNz}}.
\]
The two Bessel arguments are $y$ and $\beta y$, both of size
$qN/C_*\asymp K^2/T\gg K^{1-\varepsilon_2}L$.
By \eqref{TKM} and \eqref{eq:ndnz-averaged-kernel}, the change of
variables gives
\begin{align}\label{eq:joint-weight-after-substitution}
&\frac{16\pi^2\mathrm{i}^\ell U(n/N)}{c^2\sqrt{nd}}\,
\mathcal K_M(qn-\Delta,qn;c)\notag\\
&\quad=
\frac{32\pi^2\mathrm{i}^\ell qN\,U(z)}{c_*^2z\sqrt{MdN}}
\sum_{k\equiv0\bmod2}\mathrm{i}^k h_{K,L}(k)V_k(dNz)\notag\\
&\qquad\times\int_0^\infty
U(w^2/z)V_k(Mw^2/z)
J_\kappa(y)J_{\ell-1}(\beta y)\,\dd w.
\end{align}
As in Lemma~\ref{lem:ndnz-localization}, apply
\eqref{eq:ndnz-separated-phase-derivatives} for the plus sign and
\eqref{eq:ndnz-bessel-phase-higher-derivatives} to discard the phases
$\pm(\theta_\kappa(y)+\beta y)$ by Lemma~\ref{lemma:exponentialintegral}.
The remaining phases are $\pm\Phi_\kappa(y)$. For fixed $z$,
$\beta$ is constant in $w$ and $\dd y/\dd w=y/w$, so
\eqref{eq:ndnz-phase-derivative} gives
\begin{align*}
\frac{\dd}{\dd w}\Phi_\kappa(y)
=\frac{\sqrt{y^2-\kappa^2}-\beta y}{w}=\frac{16\pi^2qNw^2\Delta/(c_*^2z)-\kappa^2}
 {w\bigl(\sqrt{y^2-\kappa^2}+\beta y\bigr)}.
\end{align*}
A stationary point therefore requires
\[
\Delta=\frac{\kappa^2c_*^2z}{16\pi^2qNw^2}
\asymp\frac{K^2C_*^2}{qN}=H>0.
\]
Outside a fixed enlargement of this range,
\[
\left|\frac{\dd}{\dd w}\Phi_\kappa(y)\right|
\gg T+\frac{|\Delta|}{c_*}.  
\]
Using \eqref{eq:ndnz-phase-higher-derivatives} and the amplitude
bounds following \eqref{eq:ndnz-fixed-order-bessel}, apply
Lemma~\ref{lemma:exponentialintegral} in $w$, with the Bessel
expansions taken to sufficiently high order. Uniformly for $k$
in the support of $h_{K,L}$,
\[
\left|\mathcal T_{k,M}(qn-\Delta,qn;c)\right|
\ll_B\frac{c\log K}{\sqrt{qn}}\,T^{-B}
\qquad(B\geq1).
\]
Since $T\geq K^{\varepsilon_0}$, summing over $k$ and the polynomial
ranges in Corollary~\ref{cor:ndnz-effective-ranges} gives a total
contribution $O_B(K^{C-\varepsilon_0 B})$, where $C$ is fixed independently
of $B$. Taking $B>(A+C)/\varepsilon_0$ makes this $O_{A,\eps}(K^{-A})$.

To obtain \eqref{v6:eq:normalized-W}, write
$\Phi_\kappa(y)=(1-\beta)y+\theta_\kappa(y)-y$ and expand
\[
\theta_\kappa(y)-y
=-\frac{\pi\kappa}{2}-\frac\pi4+\frac{\kappa^2}{2y}
+O\left(\frac{K^4}{y^3}\right).
\]
More precisely, for every required fixed $j$,
\[
\left|\frac{\dd^j}{\dd y^j}
\left(\theta_\kappa(y)-y+\frac{\pi \kappa}{2}+\frac\pi4
-\frac{\kappa^2}{2y}\right)\right|
\ll_j\frac{K^4}{y^{3+j}}.
\]
Since $\kappa-K=O(L)$, replacing $\kappa^2/(2y)$ by $K^2/(2y)$
leaves a phase whose fixed normalized derivatives are bounded by
\[
\frac{|\kappa^2-K^2|}{y}+\frac{K^4}{y^3}
\ll\frac{TL}{K}+\frac{T^3}{K^2}.
\]
By \eqref{v6:eq:joint-dyadic-ranges} and $L\geq K^\delta$,
\[
\frac{TL}{K}\ll K^{\varepsilon_2},\qquad
\frac{T^3}{K^2}\ll K^{1-3\delta+3\varepsilon_2}=o(1)
\]
for sufficiently small $\varepsilon_2$.
The chain rule gives, for each fixed $j\geq0$,
\begin{align*}
\left|\frac{\partial^j}{\partial w^j}
\exp\left(\mathrm{i}\eta\left(
\theta_\kappa(y)-y+\frac{\pi\kappa}{2}+\frac\pi4
-\frac{K^2}{2y}\right)\right)\right|\ll_j
\left(1+\frac{TL}{K}+\frac{T^3}{K^2}\right)^j.
\end{align*}
This factor is independent of $z$, since $y=4\pi qNw/c_*$.
For the prescribed derivative orders, choose the auxiliary
exponents so that these bounds are $O(K^{\eps/2})$.
For each fixed $k$, multiply this exponential and the retained
Bessel amplitude terms by
\[
V_k(nd)V_k(x)=V_k(dNz)V_k(Mw^2/z).
\]
Lemma~\ref{cutoff of Vk} bounds the required $z,w$ derivatives
of this product of weights by $O((\log K)^2)$.
For even $k$ and $\kappa=k-1$,
\[
\mathrm{i}^k\exp\left(\mathrm{i}\eta\left(
-\frac{\pi \kappa}{2}-\frac\pi4+\frac{\pi(2\ell-1)}4
\right)\right)
=\exp\left(\frac{\mathrm{i}\eta\pi\ell}{2}\right).
\]
This constant is also included in the amplitude.
The terms $(1-\beta)y$ and $K^2/(2y)$ give the phase in
\eqref{v6:eq:normalized-W}, since
\[
\frac{(1-\beta)y}{2\pi}=\frac{2qNw}{c_*}
 \left(1-\sqrt{1-\frac{\Delta}{qNz}}\right),
\qquad
\frac{K^2}{4\pi y}=\frac{K^2c_*}{16\pi^2qNw}.
\]
The leading Bessel amplitudes satisfy
\[
\frac{1}{2\pi(y^2-\kappa^2)^{1/4}(\beta y)^{1/2}}
=\frac{c_*}{8\pi^2qNw}
 \left(1-\frac{\kappa^2}{y^2}\right)^{-1/4}\beta^{-1/2}.
\]
After extracting $c_*/(qN)$, the remaining factors and the
retained higher terms have bounded derivatives of every fixed
order. Factor $L$ out of the sum over $k$, leaving the weights
$h_{K,L}(k)/L$. Using $r=qd$ and $n=Nz$, the exterior factors are
\[
\frac{2qN}{c_*^2\sqrt{Mz}}\,
\frac1{\sqrt{dNz}}\,
\frac{c_*}{qN}\,L
=\frac L{\sqrt r}\,
\frac{q}{c_*\sqrt{MqN}}\,\frac2z.
\]
The weighted sum of the remaining phase, Bessel and $V_k$
factors over even $k$, including all cutoffs and
$32\pi^2\mathrm{i}^\ell/z$, defines $A_{\Delta,c_*}^{\eta}(z,w)$.
The latter combines $2/z$ with the factor $16\pi^2\mathrm{i}^\ell$
in \eqref{eq:ndnz-hecke-error-representation}.
Since $L^{-1}\sum_{k\equiv0\bmod2}|h_{K,L}(k)|\ll1$,
the preceding estimates give \eqref{eq:joint-amplitude-derivatives}.
Substituting the Bessel expansions in
\eqref{eq:joint-weight-after-substitution} gives
\eqref{eq:joint-original-weight-to-W}.

On the support, $z\in[1/2,2]$ and
$\Delta/(qN)=O(H/(qN))=o(1)$.
Hence $z-\Delta/(qN)$ and $z-\Delta/(2qN)$ are bounded away
from zero. The smooth compact support then makes
\eqref{v6:eq:explicit-Gtilde} entire in $s$ and rapidly decreasing
on each vertical line.
In \eqref{v6:eq:explicit-Gtilde}, set $u=2qNz-\Delta$. Then
\begin{align*}
(2qN)^s\left(z-\frac{\Delta}{2qN}\right)^{s-1}\dd z
=(2qN)^s\left(\frac{u}{2qN}\right)^{s-1}
\frac{\dd u}{2qN}=u^{s-1}\dd u,
\end{align*}
and
\begin{align*}
\frac{\sqrt{z(z-\Delta/(qN))}}{2z-\Delta/(qN)}
=\frac{\sqrt{u^2-\Delta^2}}{2u}.
\end{align*}
Consequently,
\begin{align}\label{eq:joint-Gtilde-Mellin-form}
\widetilde G_{\Delta,c_*}(s)
={}&\frac K{\sqrt{qN}}
\int_{qN-\Delta}^{4qN-\Delta}
W_{\Delta,c_*}^{\eta}\left(\frac{u+\Delta}{2qN}\right)
\left(\frac{\sqrt{u^2-\Delta^2}}{2u}\right)^{1-\ell}
 u^{s-1}\,\dd u.
\end{align}
Since $\mathfrak{q}=qb$ and $\mathfrak{q}t-m=\Delta$,
\[
\frac{\mathfrak{q}t+m+\Delta}{2qN}=\frac{bt}{N},\qquad
(\mathfrak{q}t+m)^2-\Delta^2=4\mathfrak{q}tm.
\]
The Mellin inversion of \eqref{eq:joint-Gtilde-Mellin-form}
at $u=\mathfrak{q}t+m$ therefore gives
\begin{align}\label{eq:joint-Mellin-inversion}
\frac1{2\pi\mathrm{i}}\int_{(2)}
(\mathfrak{q}t+m)^{-s}\widetilde G_{\Delta,c_*}(s)\,\dd s=\frac K{\sqrt{qN}}
\left(\frac{\sqrt{\mathfrak{q}tm}}{\mathfrak{q}t+m}\right)^{1-\ell}
W_{\Delta,c_*}^{\eta}(bt/N).
\end{align}
Multiplying \eqref{eq:joint-Mellin-inversion} by the coefficients
in \eqref{eq:shifted-dirichlet-definition} cancels the factor
with exponent $1-\ell$.
The series is absolutely convergent
on $\Re(s)=2$, so summation gives
\begin{align*}
\frac1{2\pi \mathrm{i}}\int_{(2)}
D_g(s;\mathfrak{q},1,\Delta)
\widetilde G_{\Delta,c_*}(s)\dd s
={}&\frac K{\sqrt{qN}}
\sum_{\substack{t,m\geq1\\\mathfrak{q}t-m=\Delta}}
\lambda_g(t)\lambda_g(m)
W_{\Delta,c_*}^{\eta}(bt/N).
\end{align*}
Multiplying by
$Lq\mu(b)\lambda_g(d/b)R_c(\Delta)/(K\sqrt{rM}\,c_*)$
and summing over $c_*,\Delta$ gives
\eqref{v6:eq:block-normalization} by \eqref{v6:eq:explicit-P}.

The two transforms are entire by \eqref{v6:eq:explicit-Gtilde}, since
the powers have positive real bases on the support. Substituting
$z=\Delta\zeta/H$ gives
\begin{align}\label{eq:joint-rescaled-Mellin-integral}
\mathcal G_{\Delta,c_*}(s)
={}\frac{K\Delta}{\sqrt{qN}}\left(\frac{2qN}{H}\right)^s
\int_{H/(2\Delta)}^{2H/\Delta}&
W_{\Delta,c_*}^{\eta}(\Delta\zeta/H)\notag\\
&\times
\left(\frac{\sqrt{\zeta(\zeta-H/(qN))}}
{2\zeta-H/(qN)}\right)^{1-\ell}
\left(\zeta-\frac H{2qN}\right)^{s-1}\dd\zeta .
\end{align}
At $z=\Delta\zeta/H$, the phase multiplying $\eta$ in
\eqref{v6:eq:normalized-W} has its unique critical point at
\[
w_0(\zeta)
=\frac{c_*}{4\pi C_*}
\sqrt{\frac{\zeta}{2}\left(1+\sqrt{1-\frac H{qN\zeta}}\right)}.
\]
Using \eqref{eq:joint-TH-definitions}, we have
\[
\frac{2qN}{c_*}\left(1-\sqrt{1-\frac H{qN\zeta}}\right)
=\frac{K^2c_*}{16\pi^2qNw_0(\zeta)^2},
\qquad
\Psi_H(\zeta)=\frac{K^2c_*}{8\pi^2qNw_0(\zeta)}.
\]
Hence, for every $w>0$,
\begin{align*}
&\frac{2qNw}{c_*}\left(1-\sqrt{1-\frac H{qN\zeta}}\right)
+\frac{K^2c_*}{16\pi^2qNw}-\Psi_H(\zeta)\\
&\qquad=\frac{K^2c_*}{16\pi^2qN}\,
\frac{(w-w_0(\zeta))^2}{w_0(\zeta)^2w}.
\end{align*}
We make the change of variables
\[
v=\frac{w-w_0(\zeta)}{w_0(\zeta)\sqrt w},
\qquad
\frac{\partial v}{\partial w}
=\frac{w+w_0(\zeta)}{2w_0(\zeta)w^{3/2}}.
\]
For each $\zeta$ this maps $(0,\infty)$ increasingly onto $\mathbb R$.
On the support of the amplitude, $w,w_0(\zeta)\asymp1$ and
$\partial v/\partial w\asymp1$. The displayed formula for $w_0$
bounds each of its fixed derivatives. Differentiating the change
of variables shows that the inverse $w=w(\zeta,v)$ has bounded
mixed derivatives of each fixed order on the support. Consequently
\begin{align}\label{eq:joint-exact-kernel}
W_{\Delta,c_*}^{\eta}(\Delta\zeta/H)
={}&e\bigl(\eta\Psi_H(\zeta)\bigr)
\int_{\mathbb R}
A_{\Delta,c_*}^{\eta}(\Delta\zeta/H,w(\zeta,v))
\frac{\partial w}{\partial v}\,
e\left(\frac{\eta K^2c_*}{16\pi^2qN}v^2\right)\dd v .
\end{align}
The transformed amplitude is supported in a fixed compact interval
in $v$. By \eqref{eq:joint-amplitude-derivatives}, its required mixed
derivatives are $O(K^\eps)$. Since the quadratic phase is independent
of $\zeta$ and $K^2c_* /(qN)\asymp T$, differentiation in $\zeta$
followed by the second derivative test gives, for each fixed $j\geq0$,
\begin{align*}
&\left|\frac{\dd^j}{\dd\zeta^j}
\int_{\mathbb R}A_{\Delta,c_*}^{\eta}(\Delta\zeta/H,w(\zeta,v))
\frac{\partial w}{\partial v}\,
e\left(\frac{\eta K^2c_*}{16\pi^2qN}v^2\right)\dd v\right|
\ll_{j,\eps} K^\eps\left(\frac{K^2c_*}{qN}\right)^{-1/2}.
\end{align*}
Here the mixed derivatives through order $j+1$ bound the supremum
and total variation of the differentiated amplitude.

Inserting \eqref{eq:joint-exact-kernel} into
\eqref{eq:joint-rescaled-Mellin-integral}, define
\begin{align*}
F_{\Delta,c_*}(\zeta)
:={}&-\frac{\Delta}{H}\frac{K\sqrt{c_*}}{\sqrt{qN}}
\left(\frac{\sqrt{\zeta(\zeta-H/(qN))}}
{2\zeta-H/(qN)}\right)^{1-\ell}\\
&\times\int_{\mathbb R}
A_{\Delta,c_*}^{\eta}(\Delta\zeta/H,w(\zeta,v))
\frac{\partial w}{\partial v}\,
e\left(\frac{\eta K^2c_*}{16\pi^2qN}v^2\right)\dd v .
\end{align*}
The preceding estimate gives the asserted derivative bounds for
$F_{\Delta,c_*}$. Its support follows from $\Delta\zeta/H\in[1/2,2]$.
By the definition of $e_\eta$, the integral for
$\mathcal G_{\Delta,c_*}(s)$ is
\[
-\frac H{\sqrt{c_*}}\left(\frac{2qN}{H}\right)^s
\int_{H/(2\Delta)}^{2H/\Delta}
F_{\Delta,c_*}(\zeta)\frac{\dd}{\dd\zeta}e_\eta(\zeta,s)\dd\zeta .
\]
Integration by parts gives \eqref{eq:joint-exact-transform}, since
$F_{\Delta,c_*}$ vanishes at both endpoints.

For \eqref{v6:eq:primitive-bound}, set $x=\log(u-H/(2qN))$.
The amplitude becomes $e^{\Re(s)x}$ and the phase becomes
\[
\frac{\Im(s)}{2\pi}x+\eta\Psi_H\left(e^x+\frac H{2qN}\right).
\]
On the fixed compact interval in $x$, the formula for $\Psi_H$ gives
\[
\frac{\dd^2}{\dd x^2}\Psi_H\left(e^x+\frac H{2qN}\right)
=\frac{T}{8\pi}e^{-x/2}+O\left(\frac{TH}{qN}\right).
\]
Indeed, $\Psi_H(e^x+H/(2qN))/T$ is smooth in $x$ and $H/(qN)$
near $H/(qN)=0$, where it equals $e^{-x/2}/(2\pi)$.
Since $H/(qN)=o(1)$, the absolute second derivative of the phase
is $\asymp T$. The amplitude has bounded supremum and total
variation for $|\Re(s)|\leq C$. The second derivative test gives
\eqref{v6:eq:primitive-bound}.

For the tails, insert \eqref{v6:eq:normalized-W} into
\eqref{eq:joint-rescaled-Mellin-integral}. At fixed $w$, the
derivative of the phase in $\zeta$ is
\[
\frac{\Im(s)}{2\pi(\zeta-H/(2qN))}
-\eta\frac{Hw}{c_*\zeta^2\sqrt{1-H/(qN\zeta)}}.
\]
The second term has absolute value $\asymp H/c_*\asymp T$ on
the support. Hence the displayed derivative has absolute value
$\gg T+|\Im(s)|$ for $|\Im(s)|\notin[c_-T,c_+T]$ with suitable
fixed $c_-,c_+$. Each fixed higher derivative is
$O_j(T+|\Im(s)|)$. The remaining exponential factor in each term of the $k$-sum defining
$A_{\Delta,c_*}^{\eta}$ depends only on $w$, as shown above, so the
bounds for the $V_k$ weights and the Bessel amplitudes give
$O_j((\log K)^2)$ for each fixed $\zeta$ derivative of the amplitude
in \eqref{eq:joint-rescaled-Mellin-integral} on $|\Re(s)|\leq C$.
After extraction of the exterior factors, $J$ integrations by
parts in $\zeta$ give a bound
\[
\ll_{J,C}(\log K)^2(T+|\Im(s)|)^{-J}.
\]
Since $T\geq K^{\varepsilon_0}$ and the exterior factors for both
transforms are polynomially bounded in $K$ on fixed strips, the
asserted tail estimates follow by taking $J$ sufficiently large.
\end{proof}

\subsubsection{The estimate for a dyadic block}
\begin{lemma}\label{lem:joint-spectral-block}
Under the hypotheses of Lemma~\ref{lem:joint-dyadic-kernel},
fix $\eta\in\{\pm1\}$ as in \eqref{v6:eq:explicit-P}.
Up to $O_A(K^{-A})$ for every fixed $A>0$, the block $\mathcal P(C_*,T)$
is the sum of the residues of the individual spectral integrands at $s=1/2\pm\mathrm{i}t_\omega$ for
real $|t_\omega|\asymp T$ and the integrals on $\Re(s)=-1/2$.
The residues are summed over the cusp forms and integrated over the
Eisenstein spectrum. The integrals on $\Re(s)=-1/2$ also include
the exceptional spectrum. If $\mathcal P_{\rm pole}$ and
$\mathcal P_{\rm line}$ denote the absolute values of these two
contributions, given in \eqref{eq:joint-Ppole-explicit} and
\eqref{eq:joint-Pline-explicit} below, respectively, then we have
\begin{align*}
\mathcal P_{\rm pole}
\ll K^\eps\sqrt M\sqrt{\mathfrak{q}+qN/K^2}\,T^{3/2},
\end{align*}
and
\begin{align*}
\mathcal P_{\rm line}
\ll K^\eps\sqrt M\sqrt{\mathfrak{q}+qN/K^2}\,
\frac H{qN}T^{5/2}.
\end{align*}
\end{lemma}

\begin{proof}
By \eqref{v6:eq:shifted-spectral-expansion} and
\eqref{v6:eq:uniform-transform-class}, the Mellin integrand for each
spectral component of \eqref{v6:eq:explicit-P} is
\begin{align*}
\mathcal B_\ell(s,t_\omega)\langle U_{\mathfrak{q}},\omega\rangle
\sum_{\Delta\asymp H}\frac{\rho_\omega(\Delta)}{\sqrt\Delta}
\sum_{c_*\asymp C_*}\frac{R_c(\Delta)}{c_*}\mathcal G_{\Delta,c_*}(s).
\end{align*}
The range of $\Delta$ can be split into $O(1)$ dyadic intervals
when applying \eqref{v6:eq:uniform-large-sieve}, with the endpoint
$\Delta=1$ treated as in Lemma~\ref{lem:joint-spectral-inputs}.
If $H<1$ and the block is non-empty, then $H\asymp1$.

For $\Delta\ne0$ and $R\geq1$, Cauchy's inequality, \eqref{v6:eq:uniform-Good}
and \eqref{v6:eq:uniform-large-sieve} give
\begin{align}\label{eq:joint-spectral-coefficient-mass}
&\int_{\substack{\omega\in\mathfrak S_{\mathfrak{q}}\\|t_\omega|\leq R}}
|\rho_\omega(\Delta)\langle U_{\mathfrak{q}},\omega\rangle|\,\dd\mu(\omega)\notag\\
&\quad\leq
\left(\int_{\substack{\omega\in\mathfrak S_{\mathfrak{q}}\\|t_\omega|\leq R}}
\frac{|\rho_\omega(\Delta)|^2}{\cosh(\pi t_\omega)}
\,\dd\mu(\omega)\right)^{1/2}
\left(\int_{\substack{\omega\in\mathfrak S_{\mathfrak{q}}\\|t_\omega|\leq R}}
|\langle U_{\mathfrak{q}},\omega\rangle|^2\cosh(\pi t_\omega)
\,\dd\mu(\omega)\right)^{1/2}\notag\\
&\quad\ll(\mathfrak{q}|\Delta| R)^\eps\sqrt{\mathfrak{q}}\,R^\ell
\left(R^2+\frac{|\Delta|}{\mathfrak{q}}\right)^{1/2},
\end{align}
since $\cosh(\pi t_\omega)\leq e^{\pi|t_\omega|}$.
The large sieve with $R=1$ and \eqref{v6:eq:exceptional-projections}
likewise give
\begin{align}\label{eq:joint-exceptional-coefficients}
\sum_{t_\nu=\mathrm{i}\beta_\nu}
|\rho_\nu(\Delta)\langle U_{\mathfrak{q}},u_\nu\rangle|
&\ll(\mathfrak{q}|\Delta|)^\eps
\left(1+\frac{|\Delta|}{\mathfrak{q}}\right)^{1/2}.
\end{align}
On $s=2+\mathrm{i}t$, Stirling's formula in
\eqref{v6:eq:continuous-Mellin-kernel} gives
\[
|\mathcal B_\ell(s,\tau)|\ll_\ell
\begin{cases}
(1+|t|)^{-\ell},&|\tau|\leq4(1+|t|),\\
(1+|\tau|)^{1/2}e^{-3\pi|\tau|/8},&|\tau|>4(1+|t|),
\end{cases}
\]
where $\tau$ is real or $\tau=\mathrm{i}\beta$ with $0<\beta\leq\theta$.
Apply \eqref{eq:joint-spectral-coefficient-mass} with $R=4(1+|t|)$
and sum the exponential tail over dyadic intervals. Together with
\eqref{eq:joint-exceptional-coefficients}, this gives
\begin{align*}
&\int_{\mathfrak S_{\mathfrak{q}}}
\left|\mathcal B_\ell(s,t_\omega)
\rho_\omega(\Delta)\langle U_{\mathfrak{q}},\omega\rangle\right|
\,\dd\mu(\omega)
+\sum_{t_\nu=\mathrm{i}\beta_\nu}
\left|\mathcal B_\ell(s,t_\nu)
\rho_\nu(\Delta)\langle U_{\mathfrak{q}},u_\nu\rangle\right|\\
&\quad\ll_\eps
\bigl(\mathfrak{q}|\Delta|(1+|t|)\bigr)^\eps
\sqrt{\mathfrak{q}}
\left((1+|t|)^2+\frac{|\Delta|}{\mathfrak{q}}\right)^{1/2}.
\end{align*}
By the rapid decay of $\mathcal G_{\Delta,c_*}$ in
Lemma~\ref{lem:joint-dyadic-kernel} and the polynomial ranges in
Corollary~\ref{cor:ndnz-effective-ranges}, we may restrict the Mellin
integral in \eqref{v6:eq:explicit-P} to $c_-T\leq|t|\leq c_+T$,
with error $O_A(K^{-A})$.

On this band, \eqref{eq:joint-exact-transform},
\eqref{eq:joint-F-derivative-bound} and \eqref{v6:eq:primitive-bound}
give, for some fixed $C>0$,
\begin{align*}
&\sum_{c_*\asymp C_*}\frac1{c_*}
\sum_{\Delta\asymp H}\frac{|R_c(\Delta)|}{\sqrt\Delta}
\int_{c_-T\leq|t|\leq c_+T}
|\mathcal G_{\Delta,c_*}(2+\mathrm{i}t)|\,\dd t
\ll K^C.
\end{align*}
For $j\geq0$, apply \eqref{eq:joint-spectral-coefficient-mass} on
$2^{j+2}(1+c_+T)<|t_\omega|\leq2^{j+3}(1+c_+T)$, using the
exponential Gamma bound there. Absorbing the polynomial factors into
$K^C$ and the exponential, the total tail beyond $4(1+c_+T)$ is
\[
\ll K^C\sum_{j\geq0}\exp\bigl(-\pi 2^j(1+c_+T)\bigr)
\ll_{A,\eps}K^{-A},
\]
since $T\geq K^{\varepsilon_0}$. We may thus restrict the real
spectrum to $|t_\omega|\leq4(1+c_+T)$.

We now shift to $\Re(s)=-1/2$. For
$c_-T/2\leq |t_\omega|\leq2c_+T$, first extend the segments to
$c_-T/4\leq|\Im(s)|\leq4c_+T$. The only poles crossed are
$s=1/2\pm\mathrm{i}t_\omega$ by \eqref{v6:eq:continuous-Mellin-kernel}.
For the other real parameters these poles lie outside the original
height range, and the exceptional poles are real. In both cases
we shift the original segments.

Take the original cutoffs strictly in the decay range of
Lemma~\ref{lem:joint-dyadic-kernel}. The added segments and
horizontal sides then lie in this range, and the horizontal
sides stay $\gg T$ from the poles. Stirling's formula,
\eqref{eq:joint-spectral-coefficient-mass} and
\eqref{eq:joint-exceptional-coefficients} give $O_A(K^{-A})$ for
the total contribution after summing and integrating over the spectrum.

For $|\tau|\asymp T$, \eqref{v6:eq:continuous-Mellin-kernel} gives
\begin{align*}
\mathop{\rm Res}_{s=1/2\pm\mathrm{i}\tau}\mathcal B_\ell(s,\tau)
&=\frac{2^{\ell-3/2\pm\mathrm{i}\tau}\pi^{\ell-1/2}
\Gamma(\pm\mathrm{i}\tau)}{\Gamma(\ell-1/2\pm\mathrm{i}\tau)}
\ll T^{1/2-\ell}.
\end{align*}
Thus the absolute value of the residue contribution is
\begin{align}\label{eq:joint-Ppole-explicit}
\mathcal P_{\rm pole}
={}&\Bigg|
\int_{\substack{\omega\in\mathfrak S_{\mathfrak{q}}\\
c_-T/2\leq|t_\omega|\leq2c_+T}}
\langle U_{\mathfrak{q}},\omega\rangle
\sum_{\pm}
\frac{2^{\ell-3/2\pm\mathrm{i}t_\omega}\pi^{\ell-1/2}
\Gamma(\pm\mathrm{i}t_\omega)}
{\Gamma(\ell-1/2\pm\mathrm{i}t_\omega)}
\notag\\
&\qquad\times
\sum_{\Delta\asymp H}\frac{\rho_\omega(\Delta)}{\sqrt\Delta}
\sum_{c_*\asymp C_*}\frac{R_c(\Delta)}{c_*}
\mathcal G_{\Delta,c_*}(1/2\pm\mathrm{i}t_\omega)
\,\dd\mu(\omega)\Bigg|.
\end{align}
On $\Re(s)=-1/2$, the Gamma poles are at distance at least $1-\theta$
from the line. The tails in Lemma~\ref{lem:joint-dyadic-kernel}
therefore allow all shifted integrals to be extended to
$c_-T/4\leq |t|\leq4c_+T$ with total error $O_A(K^{-A})$.
Since $\dd s=\mathrm{i}\dd t$ on this line, the absolute value of the
shifted contribution is
\begin{align}\label{eq:joint-Pline-explicit}
\mathcal P_{\rm line}
={}&\frac1{2\pi}\Bigg|
\int_{c_-T/4\leq|t|\leq4c_+T}
\sum_{\Delta\asymp H}\frac1{\sqrt{\Delta}}
\sum_{c_*\asymp C_*}\frac{R_c(\Delta)}{c_*}
\mathcal G_{\Delta,c_*}(-1/2+\mathrm{i}t)
\notag\\
&\quad\times\Bigg(
\int_{\substack{\omega\in\mathfrak S_{\mathfrak{q}}\\
 |t_\omega|\leq4(1+c_+T)}}
\mathcal B_\ell(-1/2+\mathrm{i}t,t_\omega)
\rho_\omega(\Delta)\langle U_{\mathfrak{q}},\omega\rangle
\,\dd\mu(\omega)\notag\\
&\qquad\quad+
\sum_{t_\nu=\mathrm{i}\beta_\nu}
\mathcal B_\ell(-1/2+\mathrm{i}t,\mathrm{i}\beta_\nu)
\rho_\nu(\Delta)\langle U_{\mathfrak{q}},u_\nu\rangle
\Bigg)\,\dd t\Bigg|.
\end{align}
Consequently
$|\mathcal P(C_*,T)|\leq \mathcal P_{\rm pole}
+\mathcal P_{\rm line}+O_A(K^{-A})$.

The divisor formula \eqref{eq:ndnz-ramanujan} and
$c=c_*\sqrt{M/(qN)}\in\mathbb N$ give
\begin{align}\label{v6:eq:arithmetic-rescaling-cancellation}
\sum_{c_*\asymp C_*}|R_c(\Delta)|
&\leq \sum_{v\mid\Delta}v
\sum_{\substack{c_*\asymp C_*\\v\mid c}}1
\ll C_*\sqrt{\frac M{qN}}\,\tau(\Delta).
\end{align}
By \eqref{eq:joint-F-derivative-bound},
\eqref{v6:eq:arithmetic-rescaling-cancellation} and the divisor bound,
we have, uniformly in $\zeta$,
\begin{align}
\sum_{\Delta\asymp H}\frac1{\Delta}
\left|\sum_{c_*\asymp C_*}
\frac{R_c(\Delta)}{c_*^{3/2}}F'_{\Delta,c_*}(\zeta)\right|^2
&\ll K^{\eps/2}\frac{M}{qNC_*}
\sum_{\Delta\asymp H}\frac{\tau(\Delta)^2}{\Delta}\notag\\
&\ll K^\eps\frac M{qNC_*}.
\label{eq:joint-residue-coefficient-norm}
\end{align}
By \eqref{eq:joint-exact-transform}, the coefficient of
$\rho_\omega(\Delta)$ in the sums over $\Delta$ in
\eqref{eq:joint-Ppole-explicit} and \eqref{eq:joint-Pline-explicit}
has the representation
\begin{align}\label{eq:joint-explicit-Mellin-coefficient}
\frac1{\sqrt\Delta}\sum_{c_*\asymp C_*}
 \frac{R_c(\Delta)}{c_*}\mathcal G_{\Delta,c_*}(s)
={}&H\left(\frac{2qN}H\right)^s
\int_{\zeta_0}^{\infty} e_\eta(\zeta,s)
\frac1{\sqrt\Delta}\sum_{c_*\asymp C_*}
\frac{R_c(\Delta)}{c_*^{3/2}}F'_{\Delta,c_*}(\zeta)\,\dd\zeta .
\end{align}
Here the integrand is extended by zero outside the fixed compact
interval containing all supports of $F'_{\Delta,c_*}$.

At $s=1/2\pm\mathrm{i}t_\omega$, the factor
$e_\eta(\zeta,s)$ is independent of $\Delta$ and $c_*$
and is $O(T^{-1/2})$ by \eqref{v6:eq:primitive-bound}.
For each $\zeta$, the large sieve
\eqref{v6:eq:uniform-large-sieve} and
\eqref{eq:joint-residue-coefficient-norm} give
\begin{align*}
&\int_{\substack{\omega\in\mathfrak S_{\mathfrak{q}}\\
c_-T/2\leq|t_\omega|\leq2c_+T}}
\frac1{\cosh(\pi t_\omega)}
\left|\sum_{\Delta\asymp H}
\frac{\rho_\omega(\Delta)}{\sqrt\Delta}
\sum_{c_*\asymp C_*}
\frac{R_c(\Delta)}{c_*^{3/2}}F'_{\Delta,c_*}(\zeta)\right|^2
\,\dd\mu(\omega)\\
&\qquad\ll K^\eps\left(T^2+\frac H{\mathfrak{q}}\right)
\frac M{qNC_*}.
\end{align*}
In \eqref{eq:joint-Ppole-explicit}, the product of the Gamma
residue, the primitive, and $H(2qN/H)^s$ has absolute value
$\ll \sqrt{qNH}\,T^{-\ell}$. By Cauchy's inequality with weights
$\cosh(\pi t_\omega)^{\pm1}$ and \eqref{v6:eq:uniform-Good},
we obtain
\begin{align*}
\mathcal P_{\rm pole}
&\ll K^\eps\sqrt{qNH}\,T^{-\ell}(\mathfrak{q}T^{2\ell})^{1/2}
\left(\left(T^2+\frac H{\mathfrak{q}}\right)
\frac M{qNC_*}\right)^{1/2}
=K^\eps\sqrt M\sqrt{\mathfrak{q}+qN/K^2}\,T^{3/2}.
\end{align*}
Here and below we use
\begin{align*}
\frac H{C_*}=T,\qquad\frac H{T^2}=\frac{qN}{K^2}.
\end{align*}

For $s=-1/2+\mathrm{i}t$ with $|t|\asymp T$, Minkowski's inequality
in \eqref{eq:joint-explicit-Mellin-coefficient} gives
\begin{align}\label{eq:joint-mellin-coefficient-norm}
&\sum_{\Delta\asymp H}\frac1{\Delta}
\left|\sum_{c_*\asymp C_*}\frac{R_c(\Delta)}{c_*}
\mathcal G_{\Delta,c_*}(s)\right|^2\notag\\
&\quad\leq H^2\frac H{2qN}
\left(\int_{\zeta_0}^{\infty}|e_\eta(\zeta,s)|
\left(\sum_{\Delta\asymp H}\frac1{\Delta}
\left|\sum_{c_*\asymp C_*}
\frac{R_c(\Delta)}{c_*^{3/2}}F'_{\Delta,c_*}(\zeta)\right|^2
\right)^{1/2}\dd\zeta\right)^2\notag\\
&\quad\ll K^\eps H^2\frac H{qN}T^{-1}\frac M{qNC_*}
=K^\eps M\left(\frac H{qN}\right)^2,
\end{align}
by \eqref{v6:eq:primitive-bound},
\eqref{eq:joint-residue-coefficient-norm} and $H=C_*T$.
On the shifted segments, Stirling's formula gives
\[
|\mathcal B_\ell(-1/2+\mathrm{i}t,\tau)|
\ll_\ell T^{1-\ell}\bigl(1+\bigl||t|-|\tau|\bigr|\bigr)^{-1}
\]
for real $|\tau|\ll T$, and $O_\ell(T^{-\ell})$ for
$\tau=\mathrm{i}\beta$ with $0<\beta\leq\theta$. In both cases,
\begin{align}\label{eq:joint-gamma-integral-left}
\int_{c_-T/4\leq|t|\leq4c_+T}
|\mathcal B_\ell(-1/2+\mathrm{i}t,\tau)|^2\,\dd t
&\ll T^{2-2\ell}.
\end{align}
By \eqref{eq:joint-gamma-integral-left}, Cauchy's inequality in $t$
gives a factor $O(T^{1-\ell})$ for the Gamma kernel in each
spectral component of \eqref{eq:joint-Pline-explicit}. Applying
Cauchy's inequality over the spectrum, we use
\begin{align*}
&\int_{\substack{\omega\in\mathfrak S_{\mathfrak{q}}\\
|t_\omega|\leq4(1+c_+T)}}
|\langle U_{\mathfrak{q}},\omega\rangle|^2
\cosh(\pi t_\omega)\,\dd\mu(\omega)
+\sum_{t_\nu=\mathrm{i}\beta_\nu}
|\langle U_{\mathfrak{q}},u_\nu\rangle|^2\cos(\pi\beta_\nu)\\
&\qquad\ll K^\eps\mathfrak{q}T^{2\ell}
\end{align*}
for the squared projection norm, by \eqref{v6:eq:uniform-Good}
and \eqref{v6:eq:exceptional-projections}.

For the other squared norm, interchange the integral in $t$ with
the spectral sums and integrals. At each fixed $s=-1/2+\mathrm{i}t$
in the common height range, \eqref{v6:eq:uniform-large-sieve} and
\eqref{eq:joint-mellin-coefficient-norm} give
\begin{align*}
&\int_{\substack{\omega\in\mathfrak S_{\mathfrak{q}}\\
|t_\omega|\leq4(1+c_+T)}}
\frac{1}{\cosh(\pi t_\omega)}
\left|\sum_{\Delta\asymp H}\frac{\rho_\omega(\Delta)}{\sqrt\Delta}
\sum_{c_*\asymp C_*}\frac{R_c(\Delta)}{c_*}
\mathcal G_{\Delta,c_*}(s)\right|^2\,\dd\mu(\omega)\\
&\quad+
\sum_{t_\nu=\mathrm{i}\beta_\nu}\frac{1}{\cos(\pi\beta_\nu)}
\left|\sum_{\Delta\asymp H}\frac{\rho_\nu(\Delta)}{\sqrt\Delta}
\sum_{c_*\asymp C_*}\frac{R_c(\Delta)}{c_*}
\mathcal G_{\Delta,c_*}(s)\right|^2\\
&\qquad\ll K^\eps\left(T^2+\frac H{\mathfrak{q}}\right)
\sum_{\Delta\asymp H}\frac1{\Delta}
\left|\sum_{c_*\asymp C_*}\frac{R_c(\Delta)}{c_*}
\mathcal G_{\Delta,c_*}(s)\right|^2\\
&\qquad\ll K^\eps\left(T^2+\frac H{\mathfrak{q}}\right)
M\left(\frac H{qN}\right)^2.
\end{align*}
The $t$ range has length $O(T)$.

Combining the three factors gives
\begin{align*}
\mathcal P_{\rm line}
&\ll K^\eps T^{1-\ell}(\mathfrak{q}T^{2\ell})^{1/2}
\left(\left(T^2+\frac H{\mathfrak{q}}\right)
TM\left(\frac H{qN}\right)^2\right)^{1/2}\\
&=K^\eps\sqrt M\sqrt{\mathfrak{q}+qN/K^2}\,
\frac H{qN}T^{5/2}.\qedhere
\end{align*}
\end{proof}

\subsubsection{The small stationary parameter}
\begin{lemma}\label{lem:joint-small-parameter}
Fix $1/3<\delta\leq1$ and assume $K^\delta\leq L\leq K$ and
$1\leq r\leq K$. Let $d\mid r$, $b\mid d$, $q=r/d$, and $\mathfrak{q}=qb$.
With $M,N$ in the dyadic ranges \eqref{eq:ndnz-original-lengths},
take $\varepsilon_0,\varepsilon_2>0$ sufficiently
small in terms of $\eps$. Then
\begin{align}\label{eq:joint-small-parameter-bound}
&\sum_{M,N}\Bigg|16\pi^2\mathrm{i}^{\ell}
\sum_{\substack{c\geq1\\K^2c/\sqrt{MqN}\leq K^{2\varepsilon_0}}}
\frac1{c^2}\sum_{\Delta\neq0}R_c(\Delta)
\notag\\
&\qquad\qquad\times
\sum_{\substack{u,v\geq1\\\mathfrak{q}u-v=\Delta}}
\frac{\lambda_g(u)\lambda_g(v)U(bu/N)}{\sqrt{bdu}}
\mathcal{K}_M(v,\mathfrak{q}u;c)\Bigg|\notag\\
&\qquad\qquad\qquad\ll \frac L{\sqrt r}K^{-25/32+\eps}r^{87/64}.
\end{align}
\end{lemma}

\begin{proof}
Fix $M,N$. By Corollary~\ref{cor:ndnz-effective-ranges} and
\eqref{eq:ndnz-shift-window}, we may restrict to
\begin{align*}
c&\leq C:=\sqrt{MqN}\,K^{-2+2\varepsilon_0},\qquad
0<|\Delta|\leq H_c,
\end{align*}
where
\begin{align*}
H_c&=C_1K^{\varepsilon_2}c\sqrt{\frac{qN}{M}}
 \left(1+\frac{K^2c}{\sqrt{MqN}}\right)
 \ll qNK^{-2+\varepsilon_2+4\varepsilon_0},
\end{align*}
with total error $O_A(K^{-A})$, where $C_1$ is a sufficiently large
fixed constant. For each $c$ take $H=\lfloor H_c\rfloor$ and omit
the moduli with $H=0$.

For $N/2\leq n\leq2N$, put
\begin{equation*}
W_{\Delta,c}(n)=\frac{U(n/N)}{c^2\sqrt{nd}}\,
 \mathcal K_M(qn-\Delta,qn;c),
\end{equation*}
and extend it by zero. Use $\kappa$, $Y$, and $\beta$ as in
Lemma~\ref{lem:ndnz-localization}, with $m=qn-\Delta$. Rescaling
\eqref{TKM} gives
\[
\mathcal T_{k,M}(qn-\Delta,qn;c)
=2\sqrt M\int_0^\infty V_k(Mx^2)U(x^2)
J_\kappa(xY)J_{\ell-1}(\beta xY)\dd x.
\]
On the support of $U(x^2)$,
\[
xY\asymp\beta xY\asymp\frac{\sqrt{MqN}}{c}
\geq K^{2-2\varepsilon_0},\qquad
\frac{|\Delta|}{qN}\ll K^{-2+\varepsilon_2+4\varepsilon_0}.
\]
Apply the Bessel expansions \eqref{eq:ndnz-large-order-bessel}
and \eqref{eq:ndnz-fixed-order-bessel} to sufficiently high order.
For the sum phases $\pm(\theta_\kappa(xY)+\beta xY)$, the first
$x$ derivative has size $\asymp Y$ by the identity in
\eqref{eq:ndnz-separated-phase-derivatives} with the plus sign,
since $\beta\asymp1$. The higher $x$ derivatives are $O_b(Y)$
by \eqref{eq:ndnz-bessel-phase-higher-derivatives}. After $j$ derivatives
in $n$ and multiplication by $N^j$, the amplitudes and their fixed
$x$ derivatives are $O_{j,b}(\sqrt M(\log K)\,Y^{j-1})$. Thus $B$
integrations by parts in $x$ give $O_{j,B}(\sqrt M(\log K)\,Y^{j-1-B})$.
Taking $B$ and the orders of the Bessel expansions sufficiently large,
we obtain $O_{A,j}(K^{-A})$ for the integrals with the sum phases
and for the contributions of the differentiated remainders.

For the difference phase, write
\[
\Phi_\kappa(xY)=\theta_\kappa(xY)-xY+
\frac{\Delta xY}{qn(1+\beta)}.
\]
Since $Y$ is proportional to $\sqrt n$,
\eqref{eq:ndnz-phase-higher-derivatives} with $\beta=1$ bounds
the normalized $n$ derivatives of the first difference by
$O_j(K^2/Y)$. The last term and its normalized $n$ derivatives
are $O_j(\sqrt M\,|\Delta|/(c\sqrt{qN}))$, since $n\asymp N$
and $\beta=\sqrt{1-\Delta/(qn)}\asymp1$. Hence, for $j\geq1$,
\[
N^j\left|\frac{\partial^j}{\partial n^j}\Phi_\kappa(xY)\right|
\ll_j \frac{\sqrt M\,|\Delta|}{c\sqrt{qN}}+\frac{K^2}{Y}
\ll_j K^{\varepsilon_2}\left(1+\frac{K^2c}{\sqrt{MqN}}\right).
\]
In the rescaled integral, let $\mathcal A_{k,0}^{\pm}(x,n)$ 
denote the amplitudes of $e^{\pm\mathrm{i}\Phi_\kappa(xY)}$. 
Both amplitudes include the factor $2\sqrt M\,V_k(Mx^2)U(x^2)$. The Bessel amplitude
bounds and Lemma~\ref{cutoff of Vk} give, for $j\geq0$,
\[
N^j\left|\frac{\partial^j}{\partial n^j}\mathcal A_{k,0}^{\pm}(x,n)\right|
\ll_j \frac{\sqrt M\log K}{Y}
\ll_{\varepsilon_2} \frac{cK^{\varepsilon_2}}{\sqrt{qN}}.
\]
Apply Lemma~\ref{cutoff of Vk} to $V_k(nd)$.
The product rule and $\sum_k|h_{K,L}(k)|\ll L$ then give
\begin{equation}\label{eq:joint-small-weight-derivatives}
N^j\left|\frac{\dd^j}{\dd n^j}W_{\Delta,c}(n)\right|
\ll_{j,\varepsilon_2}\frac{LK^{2\varepsilon_2}}{N\sqrt r\,c}\,P_0^j, \quad\text{for}\quad j\geq0,
\end{equation}
where
\begin{equation*}
P_0=2K^{\varepsilon_2}\left(1+\frac{K^2c}{\sqrt{MqN}}\right).
\end{equation*}

Choose $\chi\in C_c^\infty((1/8,8))$ equal to $1$ on $[1/4,4]$.
On $\mathfrak{q}u-v=\Delta$, the weight in the sum, including
$c^{-2}$, equals $W_{\Delta,c}(bu)\chi(v/(qN))$. After a fixed
dyadic partition, the slopes are $qb,1$, the lengths are comparable
to $N/b,qN$, and the derivative bounds give smoothness parameters
comparable to $P_0,2$. Since $H\geq1$, we have
\[
qN\gg K^{2-\varepsilon_2-4\varepsilon_0},\qquad
\frac Nb=\frac{qN}{qb}\gg K^{1-\varepsilon_2-4\varepsilon_0}>1.
\]
Using $N\ll K^{2+\varepsilon_2}/d$ and $qb\leq r\leq K$, we have
\[
\frac{HP_0^2}{qN}\ll K^{-2+3\varepsilon_2+8\varepsilon_0},\qquad
\frac H{qb}\ll K^{2\varepsilon_2+4\varepsilon_0},\qquad
q^2N^2P_0\ll K^7.
\]
It follows that $HP_0^2(q^2N^2P_0)^{\varpi}\ll qN$ for a
sufficiently small fixed $\varpi>0$. This verifies the short shift
condition of \cite[Theorem~3]{Bernard2015}.

Apply Bernard's theorem with $h_1=1$ and coefficients $R_c(\Delta)$.
By \eqref{eq:ndnz-ramanujan},
\[
\sum_{1\leq\Delta\leq H}|R_c(\Delta)|^2
\leq \tau(c)\sum_{v\mid c}v^2\lfloor H/v\rfloor
\leq Hc\,\tau(c)^2.
\]
For negative shifts interchange the variables and use
$R_c(-\Delta)=R_c(\Delta)$. Since $\theta=7/64$,
\[
(P_0+2)^{3/2}
\left(\sqrt{P_0+2}
+\left(\frac{qN}{P_0}\right)^\theta
 \left(1+\sqrt{\frac H{qb}}\right)\right)
\ll (qN)^\theta K^{O(\varepsilon_0+\varepsilon_2)}.
\]
Using \eqref{eq:joint-small-weight-derivatives} to restore the
size of the weight and taking the loss in Bernard's theorem
sufficiently small, the contribution for fixed $M,N$, summed in
absolute value over $c$, is
\begin{align*}
&\ll K^{O(\varepsilon_0+\varepsilon_2)}
 \frac L{\sqrt r}\frac{(qN)^{1/2+\theta}}{N}
 \sum_{c\leq C}\tau(c)\sqrt{\frac{H_c}{c}}\\
&\ll \frac L{\sqrt r}K^{-2+\eps/2}
 \frac{M^{1/4}(qN)^{5/4+\theta}}{N}
\ll \frac L{\sqrt r}K^{-25/32+3\eps/4}
 r^{87/64}d^{-3/2-2\theta}.
\end{align*}
Here $H_c/c\ll K^{\varepsilon_2+2\varepsilon_0}\sqrt{qN/M}$
and $\sum_{c\leq C}\tau(c)\ll C\log(2C)$ for $C\geq1$.
The last bound uses \eqref{eq:ndnz-original-lengths} and $q=r/d$.
Summing over $M,N$ proves \eqref{eq:joint-small-parameter-bound}.
\end{proof}

\begin{proof}[Proof of Proposition~\ref{v6:prop:joint-error}]
Fix $\eps>0$ and choose $0<\varepsilon_0<\min\{1/20,\eps/100\}$
sufficiently small for Lemma~\ref{lem:joint-small-parameter} with
loss $\eps/4$. Put $\varepsilon_1=\varepsilon_0/3$ and take
$0<\varepsilon_2<\varepsilon_0/3$ sufficiently small in terms of
$\varepsilon_0$, $\delta-1/3$ and the required derivative orders.
For $T\geq K^{\varepsilon_0}$, \eqref{v6:eq:joint-dyadic-ranges} gives
\[
L\ll K^{1-\varepsilon_0+\varepsilon_2}\leq K^{1-\varepsilon_1}
\]
for sufficiently large $K$. Thus Lemmas~\ref{lem:joint-dyadic-kernel}
and \ref{lem:joint-spectral-block} apply. Use loss $K^{\varepsilon_0}$
in Lemma~\ref{lem:joint-spectral-block}.
By \eqref{eq:joint-TH-definitions},
\eqref{v6:eq:joint-dyadic-ranges}, and
\eqref{v6:eq:balanced-new-ranges},
\begin{align*}
T\ll\frac{K^{1+\varepsilon_0}}L,\qquad
C_*=\frac{qN}{K^2}T,\qquad
\frac H{qN}=\frac{T^2}{K^2},\qquad
\frac{qN}{K^2}\ll rK^{\varepsilon_0},\qquad
q,\mathfrak{q}\leq r.
\end{align*}
For $T\geq K^{\varepsilon_0}$, \eqref{v6:eq:joint-dyadic-ranges}
and $L\geq K^\delta$ give
$HT/(qN)=T^3/K^2\ll K^{1-3\delta+3\varepsilon_2}=o(1)$.
We can therefore use \eqref{v6:eq:block-normalization} and
Lemma~\ref{lem:joint-spectral-block} to bound the two contributions together, and get
\begin{align*}
\frac q{K\sqrt M}(\mathcal P_{\rm pole}+\mathcal P_{\rm line})
&\ll K^{\varepsilon_0}\frac qK
\left(\mathfrak{q}+\frac{qN}{K^2}\right)^{1/2}T^{3/2}\\
&\ll K^{3\varepsilon_0}r^{3/2}
\left(\frac K{L^3}\right)^{1/2}.
\end{align*}
For sufficiently large $K$, the remaining blocks satisfy
\[
T<K^{\varepsilon_0}\quad\Longrightarrow\quad
\frac{K^2c}{\sqrt{MqN}}
=\frac{c_*}{C_*}T
\ll K^{\varepsilon_0}\leq K^{2\varepsilon_0}.
\]
The proof of Lemma~\ref{lem:joint-small-parameter} takes absolute
values before summing over $c$, so it applies to these restricted
moduli. With loss $\eps/4$, it gives the second term of
\eqref{v6:eq:joint-error}.

Finally, Deligne's bound gives
\[
\sum_{d\mid r}\sum_{b\mid d}|\mu(b)\lambda_g(d/b)|
\leq\sum_{m\mid r}\tau(m)\tau(r/m)
\ll_{\varepsilon_0}r^{\varepsilon_0}
\leq K^{\varepsilon_0}.
\]
There are $O((\log K)^3)\ll K^{\varepsilon_0}$ dyadic
$M,N,C_*$ blocks and two choices $\eta=\pm1$.
Since
\[
\max\{3\varepsilon_0,\eps/4\}+2\varepsilon_0<\eps,
\]
summing these estimates proves \eqref{v6:eq:joint-error}.
\end{proof}

\subsection{Proof of Theorem~\ref{thm:twisted-second-moment}}
For $L\geq K^\delta$ and $r\leq K^\rho$, the terms in parentheses
in \eqref{v6:eq:joint-error} satisfy
\begin{align*}
r^{3/2}\left(\frac K{L^3}\right)^{1/2}
&\leq K^{-(3\delta-1-3\rho)/2},
\end{align*}
and
\begin{align*}
r^{87/64}K^{-25/32}
&\leq K^{-(50-87\rho)/64}.
\end{align*}
It follows that
\[
\mathcal E(r)\ll_\eps\frac L{\sqrt r}K^{-\sigma(\delta,\rho)+\eps}.
\]
Apply \eqref{eq:general-integer-second} with
$0<\varepsilon'<\min\{\eps,\delta\}$ in place of $\eps$.
After division by $L/\sqrt r$, its error term is
$O(K^{-1+\rho/2+\varepsilon'})$. Since
\[
\sigma(\delta,\rho)\leq\frac{50-87\rho}{64}<1-\frac\rho2,
\]
this proves \eqref{eq:twisted-second-moment}.

\section{Mollified moments and non-vanishing}
\label{sec:mollifier}

\subsection{The Hecke expansion of the mollifier}
Multiplying the four inverse Satake factors gives
\begin{align}\label{eq:moll-local-inverse}
L_p(s,f\otimes g)^{-1}=
1-\lambda_f(p)\lambda_g(p)p^{-s}
+\bigl(\lambda_f(p^2)+\lambda_g(p^2)\bigr)p^{-2s}
-\lambda_f(p)\lambda_g(p)p^{-3s}+p^{-4s},
\end{align}
where we used $\lambda_{F}(p)^2=\lambda_F(p^2)+1$ for $F\in\{f,g\}$.
To record the coefficient of each Hecke term in
\eqref{eq:moll-local-inverse}, define
\begin{align*}
c_g(1,1)=c_g(p^2,1)=c_g(1,p^2)=1,\quad
c_g(p,1)=c_g(p,p)= -\lambda_g(p),\quad
c_g(1,p)=\lambda_g(p^2)
\end{align*}
for every prime $p$. All other local coefficients
$c_g(p^i,p^j)$, with $i,j\geq0$, are zero. Extend these
coefficients to all positive integers $r,m$ by
\begin{align*}
c_g(r,m):= & ~
\prod_{p\mid rm}c_g\bigl(p^{v_p(r)},p^{v_p(m)}\bigr),
\end{align*}
where the empty product is $1$. These coefficients depend
only on $g,r,m$.
In each local term, the index of $\mu_{f,g}$ is $rm^2$.
By multiplicativity of the Hecke eigenvalues,
\begin{align*}
\mu_{f,g}(n)= & ~\sum_{rm^2=n}c_g(r,m)\lambda_f(r).
\end{align*}
Substituting this identity into the definition of $M_P(f)$ gives
\begin{align*}
M_P(f)= & ~\sum_{rm^2\leq X}
\frac{c_g(r,m)\lambda_f(r)}{m\sqrt r}
P\left(\frac{\log(X/(rm^2))}{\log X}\right).
\end{align*}
For $r\leq X$, put
\begin{align*}
x_r:= & ~\sum_{m\leq\sqrt{X/r}}\frac{c_g(r,m)}m
P\left(\frac{\log(X/(rm^2))}{\log X}\right).
\end{align*}
Grouping the finite sum by $r$, we obtain
\begin{align}\label{eq:moll-Hecke-expansion}
M_P(f)= & ~\sum_{r\leq X}\frac{x_r\lambda_f(r)}{\sqrt r},
\qquad x_r\ll X^\eps.
\end{align}
Deligne's bound gives $|c_g(r,m)|\leq\tau(rm^2)$.
Since $P$ is bounded on $[0,1]$, it follows that
\begin{align*}
|x_r|\ll_{P,\eps} & ~
\sum_{m\leq\sqrt{X/r}}\frac{\tau(rm^2)}m
\ll X^{\eps/2}\log(2X)\ll X^\eps.
\end{align*}

\subsection[Euler products for the mollified moments]{Euler products for the mollified moments}
\label{sec:moll-euler}
Write $P(u)=\sum_{j\geq2}p_ju^j$, where the sum is finite,
and put
\begin{align*}
\mathcal H_P(w):= & ~
\sum_{j\geq2}\frac{p_jj!}{(\log X)^jw^{j+1}}.
\end{align*}
For $j\geq2$ and $\sigma>0$, Perron's formula gives
\begin{align*}
\frac1{2\pi\mathrm{i}}\int_{(\sigma)}
\frac{(X/n)^w}{w^{j+1}}\dd w
= & ~\begin{cases}
\dfrac{(\log(X/n))^j}{j!},&n<X,\\[4pt]
0,&n\geq X.
\end{cases}
\end{align*}
Multiplying by $p_jj!/(\log X)^j$ and summing over $j$ gives
\begin{align}\label{eq:moll-Perron}
\frac1{2\pi\mathrm{i}}\int_{(\sigma)}
\left(\frac{X}{n}\right)^w\mathcal H_P(w)\dd w
= & ~\begin{cases}
P\left(\dfrac{\log(X/n)}{\log X}\right),&n\leq X,\\[4pt]
0,&n>X.
\end{cases}
\end{align}
At $n=X$ the weight is zero since $P(0)=0$.
The kernel $\mathcal H_P$ is integrable on every vertical line
$\Re(w)=\sigma>0$.
For $n=1$, shifting the line to the left gives the residue at
$w=0$. Applying the same argument to $w\mathcal H_P(w)$ and
expanding $X^w$ yields
\begin{align}\label{eq:moll-first-residue-P}
\Res_{w=0}X^w\mathcal H_P(w)
= & ~\sum_{j\geq2}p_j=P(1),
\end{align}
and
\begin{align}\label{eq:moll-first-residue-dP}
\Res_{w=0}X^ww\mathcal H_P(w)
= & ~\frac1{\log X}\sum_{j\geq2}jp_j
=\frac{P'(1)}{\log X}.
\end{align}

For a local Hecke expression, let $\mathfrak C_p$ denote the
coefficient of $\lambda_f(1)$ after applying the Hecke relation.
In particular,
\begin{align*}
\mathfrak C_p\bigl(\lambda_f(p^j)\bigr)=\delta_{j,0},
\qquad
\mathfrak C_p\bigl(\lambda_f(p^i)\lambda_f(p^j)\bigr)
=\delta_{i,j}.
\end{align*}
Fix an integer $t\geq1$. For
$\Re(s_1),\Re(s_2),\Re(w),\Re(z)>1/2$, define
\begin{align*}
F_{2,t}(s_1,s_2,w,z):= & ~\prod_p\mathfrak C_p\left(
\lambda_f(p^{v_p(t)})
\frac{L_p(1/2+s_1,f\otimes g)L_p(1/2+s_2,f\otimes g)}
{L_p(1/2+w,f\otimes g)L_p(1/2+z,f\otimes g)}\right),
\end{align*}
and write $F_{2,t,p}$ for its $p$-factor.

For an integer $r\geq1$ and $\Re(s_1),\Re(s_2)>1/2$, let
$C_r(s_1,s_2)$ denote the coefficient of $\lambda_f(1)$ in the
Hecke expansion of
$L(1/2+s_1,f\otimes g)L(1/2+s_2,f\otimes g)\lambda_f(r)$.
By \eqref{Hecke relation},
\begin{align*}
\lambda_f(m)\lambda_f(n)\lambda_f(r)
= & ~\sum_{d\mid(n,r)}\lambda_f(m)\lambda_f(nr/d^2).
\end{align*}
In the expansion by \eqref{Dirichlet series}, a summand contributes
to $C_r(s_1,s_2)$ exactly when $m=nr/d^2$. Writing $n=dv$ and
$m=vr/d$ for $d\mid r$, we obtain
\begin{align*}
C_r(s_1,s_2)= & ~\zeta(1+2s_1)\zeta(1+2s_2)
\sum_{d\mid r}\sum_{v\geq1}
\frac{\lambda_g(vr/d)\lambda_g(vd)}
{(vr/d)^{1/2+s_1}(vd)^{1/2+s_2}}\\
= & ~\frac{\zeta(1+2s_1)\zeta(1+2s_2)}{\sqrt r}
\sum_{d\mid r}(r/d)^{-s_1}d^{-s_2}
\mathscr D_{d,r/d}(s_1+s_2).
\end{align*}
Here we used \eqref{eq:general-D-dq}.
By \eqref{eq:general-D-euler}, the last expression continues
meromorphically to a neighbourhood of $s_1=s_2=0$.
Expanding the reciprocal factors by \eqref{eq:moll-local-inverse}
and applying the Hecke relation to their eigenvalues and
$\lambda_f(t)$ gives
\begin{align*}
F_{2,t}(s_1,s_2,w,z)
= & ~\sum_{u,v,m,n\geq1}
\frac{c_g(u,m)c_g(v,n)}
{u^{1/2+w}v^{1/2+z}m^{1+2w}n^{1+2z}}
\sum_{d\mid(u,v)}\sum_{e\mid(uv/d^2,t)}
C_{uvt/(d^2e^2)}(s_1,s_2).
\end{align*}

For $p\nmid t$, the inserted Hecke factor is $1$. By
\eqref{Dirichlet series}, for $i=1,2$,
\begin{align*}
L_p(1/2+s_i,f\otimes g)
= & ~1+\lambda_f(p)\lambda_g(p)p^{-1/2-s_i}
+\bigl(1+\lambda_f(p^2)\lambda_g(p^2)\bigr)p^{-1-2s_i}
+\cdots.
\end{align*}
By \eqref{eq:moll-local-inverse} and \(\mathfrak C_p(\lambda_f(p)^2)=1\), the four quadratic terms
and the six products of linear terms give
\begin{align*}
F_{2,t,p}
= & ~1+p^{-1-2s_1}+p^{-1-2s_2}
+\lambda_g(p^2)\bigl(p^{-1-2w}+p^{-1-2z}\bigr)\\
& ~+\lambda_g(p)^2\bigl(
p^{-1-s_1-s_2}+p^{-1-w-z}-p^{-1-s_1-w}-p^{-1-s_1-z}\\
& ~\hspace{37mm}
-p^{-1-s_2-w}-p^{-1-s_2-z}\bigr)
+O_g(p^{-2+4\eta})
\end{align*}
when $|\Re(s_1)|,|\Re(s_2)|,|\Re(w)|,|\Re(z)|\leq\eta<1/16$.

At these primes the Hecke relation preserves parity of the
degree in the four local powers, so odd degrees vanish under
$\mathfrak C_p$. The product \(\lambda_f(p^m)\lambda_f(p^n)\) contains at most \(m+1\) terms. Deligne's bound and
the fixed reciprocal polynomials then give $O_g((j+1)^{12})$
for the sum of the absolute projected coefficients of degree $j$.
A partial derivative of total order $i$ contributes at most
$(j\log p)^i$. Thus the terms of degree $j\geq4$ contribute
\begin{align*}
\ll_{g,i} & ~\sum_{j\geq4}(j+1)^{12}(j\log p)^i
p^{-j(1/2-\eta)}
\ll_{g,i}p^{-2+4\eta}(\log p)^i.
\end{align*}
For each of the finitely many primes $p\mid t$, the full projected
local series and its fixed derivatives converge absolutely on
the same strips, by the same polynomial bound in the degree
with constants depending on $t$.

Since $\lambda_g(p)^2=1+\lambda_g(p^2)$ and
\begin{align*}
L_p(v,\sym^2g)= & ~1+\lambda_g(p^2)p^{-v}
+O_g(p^{-2\Re(v)}),\qquad \Re(v)\geq1-2\eta,
\end{align*}
the local zeta and symmetric square factors account for the
quadratic terms. We obtain
\begin{align}\label{eq:moll-Euler-singular}
F_{2,t}(s_1,s_2,w,z)
= & ~\frac{\zeta(1+2s_1)\zeta(1+2s_2)
\zeta(1+s_1+s_2)\zeta(1+w+z)}
{\zeta(1+s_1+w)\zeta(1+s_1+z)
\zeta(1+s_2+w)\zeta(1+s_2+z)}
A_{2,g,t}(s_1,s_2,w,z),
\end{align}
where
\begin{align*}
A_{2,g,t}(s_1,s_2,w,z)
= & ~\frac{L(1+s_1+s_2,\sym^2g)L(1+w+z,\sym^2g)}
{L(1+s_1+w,\sym^2g)L(1+s_1+z,\sym^2g)}
\\
& ~\times\frac{L(1+2w,\sym^2g)L(1+2z,\sym^2g)}
{L(1+s_2+w,\sym^2g)L(1+s_2+z,\sym^2g)}
\\
& ~\times\prod_p\Bigg\{F_{2,t,p}
\frac{\zeta_p(1+s_1+w)\zeta_p(1+s_1+z)}
{\zeta_p(1+2s_1)\zeta_p(1+2s_2)}
\\
& ~\qquad\times
\frac{\zeta_p(1+s_2+w)\zeta_p(1+s_2+z)}
{\zeta_p(1+s_1+s_2)\zeta_p(1+w+z)}
\\
& ~\qquad\times
\frac{L_p(1+s_1+w,\sym^2g)L_p(1+s_1+z,\sym^2g)}
{L_p(1+s_1+s_2,\sym^2g)L_p(1+w+z,\sym^2g)}
\\
& ~\qquad\times
\frac{L_p(1+s_2+w,\sym^2g)L_p(1+s_2+z,\sym^2g)}
{L_p(1+2w,\sym^2g)L_p(1+2z,\sym^2g)}
\Bigg\}.
\end{align*}
Here $\zeta_p(v)=(1-p^{-v})^{-1}$.
For $p\nmid t$ the factors in braces above are $1+O_g(p^{-2+4\eta})$.
Their denominators are non-zero by the Satake bounds and the
positive real parts of the local arguments. The same local bounds
hold for $|\Re(s_i)|\leq\eta$ and $\Re(w),\Re(z)\geq-\eta$, since the corresponding local powers are at most \(p^{-1/2+\eta}\) in modulus. The regular Euler product and its fixed
partial derivatives converge normally and are bounded uniformly on
smaller closed strips.

Since $L(1,\sym^2g)>0$, the global symmetric square quotient is
holomorphic near the origin.
When $s_1=w$ and $s_2=z$, both local $L$-quotients
equal $1$, so we take the coefficient of $\lambda_f(1)$ in
$\lambda_f(p^{v_p(t)})$. The normalizing zeta and symmetric square
quotients also equal $1$. Thus
\begin{align}\label{eq:moll-arithmetic-normalization}
A_{2,g,t}(0,0,0,0)
= & ~\prod_{p\mid t}\mathfrak C_p\bigl(\lambda_f(p^{v_p(t)})\bigr)
=\mathbf 1_{t=1}.
\end{align}
The equality also covers $t=1$ by the empty product convention.
Holomorphy gives
\begin{align*}
A_{2,g,t}(s_1,s_2,w,z)
= & ~\mathbf 1_{t=1}+O_t(|s_1|+|s_2|+|w|+|z|).
\end{align*}

For $s_2=z=\xi$ in the region of absolute convergence, the two
local factors at shift $\xi$ cancel. Put
$F_{1,t}(s_1,w)=F_{2,t}(s_1,\xi,w,\xi)$, which is independent of
$\xi$. Expanding the remaining local quotient gives
\begin{align}\label{F1}
\frac{F_{1,t}(s_1,w)}{\zeta(1+2s_1)}
= & ~\sum_{r,m\geq1}\frac{c_g(r,m)}{r^{1/2+w}m^{1+2w}}
\sum_{e\mid(r,t)}
\frac{\lambda_g(rt/e^2)}{(rt/e^2)^{1/2+s_1}}.
\end{align}
Since $r/t\leq rt/e^2\leq tr$, the series and its derivative
with respect to $s_1$ converge absolutely and uniformly for
$|s_1|\leq1/4$ on $\Re(w)=1$.

On $s_2=z=\xi$, the zeta quotient in \eqref{eq:moll-Euler-singular}
reduces to $\zeta(1+2s_1)/\zeta(1+s_1+w)$, and the global symmetric
square quotient reduces to
\begin{align*}
\frac{L(1+2w,\sym^2g)}{L(1+s_1+w,\sym^2g)}.
\end{align*}
The cancellations hold as meromorphic identities, giving
\begin{align}\label{eq:moll-first-Euler-singular}
F_{1,t}(s_1,w)
= & ~\frac{\zeta(1+2s_1)}{\zeta(1+s_1+w)}
A_{1,g,t}(s_1,w),
\end{align}
where $A_{1,g,t}(s_1,w)=A_{2,g,t}(s_1,0,w,0)$ by continuation
in $\xi$ after cancellation. The local bounds and
\eqref{eq:moll-arithmetic-normalization} yield
\begin{align}\label{eq:moll-first-arithmetic-normalization}
A_{1,g,t}(0,0)=\mathbf 1_{t=1},
\qquad
A_{1,g,t}(s_1,w)=\mathbf 1_{t=1}+O_t(|s_1|+|w|).
\end{align}

\subsection{\texorpdfstring{Proof of the first moment formula in \eqref{eq:mollified-first-moment}}{Proof of the first moment formula}}
\begin{proof}
Fix a positive integer $t$. By \eqref{eq:moll-Hecke-expansion}
and the Hecke relation \eqref{Hecke relation}, we have
\begin{align*}
\lambda_f(t)M_P(f)
= & ~\sum_{r\leq X}\frac{x_r}{\sqrt r}
\sum_{e\mid(r,t)}\lambda_f\left(\frac{rt}{e^2}\right).
\end{align*}
Since $a<\delta$, we have $t\leq K^{\delta-a}$ for sufficiently large
$K$. Thus $rt/e^2\leq tX\leq L^2K^{-(\delta-a)}$, and we may apply
\eqref{eq:fmcor-digamma} at every index in the sum. Each $r$ admits at most
$\tau(t)$ choices of $e$, so the total error is
\begin{align}\label{eq:moll-first-error-sum}
\ll_t K^{-A}\sum_{r\leq X}\frac{|x_r|}{\sqrt r}
\ll_t K^{-A}X^{1/2+\eps}=o_t(L)
\end{align}
for sufficiently large $A$. The main term is
$(L/2)\int_{\mathbb R} h(u)S_t(K+1+Lu)\,\dd u$, where
\begin{align*}
S_t(k)= & ~\sum_{r\leq X}\frac{x_r}{\sqrt r}
\sum_{e\mid(r,t)}\frac{\lambda_g(rt/e^2)}{\sqrt{rt/e^2}}
\bigl(\alpha_k+2\gamma-\log(rt/e^2)\bigr).
\end{align*}
Applying \eqref{eq:moll-Perron} at $rm^2$ to \eqref{F1} gives
\begin{align*}
\frac{1}{2\pi\mathrm{i}}\int_{(1)}
X^w\mathcal H_P(w)\frac{F_{1,t}(s_1,w)}{\zeta(1+2s_1)}\,\dd w
= & ~\sum_{r\leq X}\frac{x_r}{\sqrt r}
\sum_{e\mid(r,t)}\frac{\lambda_g(rt/e^2)}{(rt/e^2)^{1/2+s_1}}.
\end{align*}
Differentiating at $s_1=0$ supplies the factor $-\log(rt/e^2)$.
The absolute convergence of \eqref{F1} and its first $s_1$
derivative on $\Re(w)=1$ justifies differentiation under the
integral. By \eqref{eq:moll-first-Euler-singular}, we obtain
\begin{align}\label{eq:moll-first-integral-kernel}
S_t(k)= & ~\frac1{2\pi\mathrm{i}}
\int_{(1)}X^w\mathcal H_P(w)B_{k,t}(w)\,\dd w, \end{align}
where
\begin{align}\label{Bktw}
B_{k,t}(w)= & ~(\alpha_k+2\gamma)\frac{A_{1,g,t}(0,w)}{\zeta(1+w)}
+\left.\frac{\partial}{\partial s_1}
\frac{A_{1,g,t}(s_1,w)}{\zeta(1+s_1+w)}\right|_{s_1=0}.
\end{align}
By \eqref{eq:common-mellin-residue}, equation \eqref{Bktw} can also
be written as
\begin{align*}
B_{k,t}(w)= & ~2\Res_{s_1=0}
\frac{G(s_1)\Gamma_k(s_1)}{s_1}F_{1,t}(s_1,w).
\end{align*}
Take the residue on $|s_1|=1/(2\log X)$. On this circle
$\zeta(1+2s_1)\ll\log X$ and $G(s_1)\Gamma_k(s_1)\ll_a1$.
For $1/\log X\leq\Re(w)\leq1$, absolute convergence bounds
the global factors with arguments $1+s_1+w$ and $1+2w$ by
powers of $\log X$. The regular Euler product is bounded there,
and $\mathcal H_P(w)\ll_P(\log X)^{-2}|w|^{-3}$ for
$|w|\log X\geq1$. We may therefore move the line in
\eqref{eq:moll-first-integral-kernel} to $\Re(w)=1/\log X$.
On this line, $|s_1|\leq |w|/2$. Thus Remark~\ref{rem:zero-free-region} and the simple zero
of $1/\zeta(1+s_1+w)$ at $s_1+w=0$ give
\[
B_{k,t}(w)\ll_t (\log X)\min\{|w|,1\}
\bigl(\log(3+|w|)\bigr)^C\ll_t(\log X)|w|.
\]
For $w$ near zero, the formula for $B_{k,t}$ and
\eqref{eq:moll-first-arithmetic-normalization} give
\[
B_{k,t}(w)=\mathbf 1_{t=1}(1+\alpha_kw)
+O_t\bigl(|w|+|\alpha_k||w|^2\bigr).
\]
Write $w=(1+\mathrm{i}x)/\log X$. The Perron weight and line element
give $|X^w\mathcal H_P(w)|\,\dd x/\log X\ll_P(1+|x|)^{-3}\dd x$.
Thus both integrands, with $B_{k,t}(w)$ and with
$\mathbf 1_{t=1}(1+\alpha_kw)$, are bounded by a constant times
$(1+|x|)^{-2}\dd x$. Their difference tends to zero on bounded
$x$ intervals uniformly for $k\asymp K$, since $\alpha_k\ll\log K$
and $\log X=2a\log K$. Dominated convergence therefore gives an
error $o_t(1)$ on replacing $B_{k,t}(w)$ by
$\mathbf 1_{t=1}(1+\alpha_kw)$ in the integral. By
\eqref{eq:moll-first-residue-P} and \eqref{eq:moll-first-residue-dP},
\[
S_t(k)=\mathbf 1_{t=1}\left(P(1)+\frac{\alpha_k}{\log X}P'(1)\right)
+o_t(1).
\]
On the support of $h$, we have
$\alpha_{K+1+Lu}/\log X=a^{-1}+O_a(1/\log X)$.
Averaging against $(L/2)h(u)\,\dd u$ and adding the error in
\eqref{eq:moll-first-error-sum} proves \eqref{eq:mollified-first-moment}.
All estimates are uniform on the compact parameter sets in
Proposition~\ref{prop:mollified-moments}, since $a$ and $\delta-a$
stay bounded away from zero.
\end{proof}

\subsection{\texorpdfstring{Proof of the second moment formula in \eqref{eq:mollified-second-moment}}{Proof of the second moment formula}}
\begin{proof}
Fix a positive integer $t$. By \eqref{eq:moll-Hecke-expansion}
and the Hecke relation, the left-hand side of
\eqref{eq:mollified-second-moment} equals
\[
\sum_{u,v\leq X}\frac{x_ux_v}{\sqrt{uv}}
\sum_{d\mid(u,v)}\sum_{e\mid(uv/d^2,t)}
\mathcal M_2\left(\frac{uvt}{d^2e^2}\right).
\]
Since $4a<1$, we have $uvt/(d^2e^2)\leq tX^2\leq K$ for
sufficiently large $K$. By \eqref{eq:general-integer-second}
and \eqref{v6:eq:joint-error}, the error for each index is
\[
\ll LK^\eps\left(
\frac{uvt}{d^2e^2}\left(\frac K{L^3}\right)^{1/2}
+\left(\frac{uvt}{d^2e^2}\right)^{55/64}K^{-25/32}+K^{-1}\right).
\]
There are at most $\tau(t)$ choices of $e$. Since $e\mid t$,
\[
\frac{uv}{d^2t}\leq\frac{uvt}{d^2e^2}\leq\frac{tuv}{d^2}.
\]
Thus $(uvt/(d^2e^2))^\beta\ll_t (uv/d^2)^\beta$ for
$\beta=1,55/64,0$. Writing $u=dm$ and $v=dn$ gives
\[
 \sum_{u,v\leq X}\frac1{\sqrt{uv}}\sum_{d\mid(u,v)}
 \left(\frac{uv}{d^2}\right)^\beta
 =\sum_{d\leq X}\frac1d
   \left(\sum_{m\leq X/d}m^{\beta-1/2}\right)^2
 \ll_\beta X^{2\beta+1}.
\]
Using
$x_u,x_v\ll X^\eps$, we thus bound the total
remainder by
\begin{align}\label{eq:general-mollified-error-sum}
\ll_t LK^\eps\left(K^{6a-(3\delta-1)/2}
                  +K^{87a/16-25/32}+K^{2a-1}\right)=o_t(L).
\end{align}
For the main term, write
\[
S_{2,t}(k):=\sum_{u,v\leq X}\frac{x_ux_v}{\sqrt{uv}}
\sum_{d\mid(u,v)}\sum_{e\mid(uv/d^2,t)}
\mathscr P_k\left(\frac{uvt}{d^2e^2}\right).
\]

For every index $r$, the coefficient $C_r$ of \S\,\ref{sec:moll-euler} gives
\begin{align}\label{eq:moll-ordered-residue}
\mathscr P_k(r)= & ~4\Res_{s_2=0}\Res_{s_1=0}
\frac{G(s_1)\Gamma_k(s_1)G(s_2)\Gamma_k(s_2)}{s_1s_2}
C_r(s_1,s_2),
\end{align}
with the residue at $s_1=0$ taken first.
To verify this identity, fix $s_2\ne0$ near zero and take
the $s_1$-circle about zero with radius at most $|s_2|/2$.
For each $d\mid r$, apply \eqref{eq:common-mellin-residue} with
$x=r/d$ and $H(s_1)=\mathscr D_{d,r/d}(s_1+s_2)$, which is
holomorphic inside this circle. Half the right-hand side of that
formula is exactly the expression in braces in \eqref{eq:general-PV-residue},
evaluated at $s=s_2$.
Multiply by $G(s_2)\Gamma_k(s_2)\zeta(1+2s_2)d^{-s_2}/s_2$
and take the residue at $s_2=0$. Summing over $d\mid r$ and
multiplying by $4/\sqrt r$ gives \eqref{eq:general-PV-residue},
which proves \eqref{eq:moll-ordered-residue}.

On fixed small nested circles avoiding $s_1+s_2=0$,
\eqref{eq:general-D-euler} gives
$C_r(s_1,s_2)\ll r^{-1/2+\eta+\eps}$ for
$|\Re(s_1)|,|\Re(s_2)|\leq\eta<1/16$.
Since $uvt/(d^2e^2)\geq uv/(d^2t)$, choosing $\eps$ sufficiently
small gives
\[
C_{uvt/(d^2e^2)}(s_1,s_2)
\ll_t\left(\frac{uv}{d^2}\right)^{-1/2+\eta+\eps}\ll_t1.
\]
On $\Re(w)=\Re(z)=1$, a summable majorant is
\begin{align*}
\sum_{u,v,m,n\geq1}
\frac{|c_g(u,m)c_g(v,n)|\tau((u,v))}
{u^{3/2}v^{3/2}m^3n^3}<\infty.
\end{align*}
This majorant gives absolute and uniform convergence on the small circles and, with the integrability of $\mathcal H_P$, justifies
termwise integration, including the fixed residue integrals.
Apply \eqref{eq:moll-Perron} at $um^2$ and $vn^2$ in the expansion of
$F_{2,t}$. The sum over $e\mid(uv/d^2,t)$ has at most $\tau(t)$ terms.
Taking the ordered residues in \eqref{eq:moll-ordered-residue} gives
\begin{multline}
S_{2,t}(k)=
\frac4{(2\pi\mathrm i)^2}\int_{(1)}\int_{(1)}
X^{w+z}\mathcal H_P(w)\mathcal H_P(z)\\
\times\Res_{s_2=0}\Res_{s_1=0}
\frac{G(s_1)\Gamma_k(s_1)G(s_2)\Gamma_k(s_2)}{s_1s_2}
F_{2,t}(s_1,s_2,w,z)\dd w\dd z .
\label{eq:moll-second-perron-residues}
\end{multline}

By \eqref{eq:gamma-ratio}, $G(s)=e^{s^2}$ and Stirling's formula,
the second logarithmic derivative of $G(s)\Gamma_k(s)$ is bounded
on a fixed small disc, uniformly for $k\asymp K$. Thus
\begin{align*}
G(s_i)\Gamma_k(s_i)=e^{\alpha_ks_i}(1+O(|s_i|^2)),
\qquad i=1,2,
\end{align*}
where the normalized factor is holomorphic. Since
$v\zeta(1+v)$ is holomorphic with value one at $v=0$,
\eqref{eq:moll-Euler-singular} and
\eqref{eq:moll-arithmetic-normalization} give
\begin{align*}
F_{2,t}(s_1,s_2,w,z)
=\frac{(s_1+w)(s_1+z)(s_2+w)(s_2+z)}
{4s_1s_2(s_1+s_2)(w+z)}
\bigl(\mathbf 1_{t=1}
+O_t(|s_1|+|s_2|+|w|+|z|)\bigr).
\end{align*}
The factor in parentheses is holomorphic near the origin.

Let $N_{k,t}(s_1,s_2,w,z)$ be its product with
$\prod_{i=1}^2G(s_i)\Gamma_k(s_i)e^{-\alpha_ks_i}$. Then
$N_{k,t}(0,0,0,0)=\mathbf 1_{t=1}$. On a smaller polydisc its fixed
derivatives are bounded uniformly for $k\asymp K$, with constants
depending on $t$.
Let $J(A)$ denote the ordered residue on the left of
\eqref{eq:moll-kernel-residue} with $\alpha_k$ replaced by $A$.
The factor $4$ in \eqref{eq:moll-second-perron-residues} cancels
the $4$ in the preceding denominator.
At $A=0$ the inner residue of
$(s_1+w)(s_1+z)/(s_1^2(s_1+s_2))$ is
$(w+z)/s_2-wz/s_2^2$, whence $J(0)=1$.
Differentiation in $A$ cancels $s_1+s_2$ and gives
\[
 J'(A)=\frac1{w+z}
 \left(\Res_{s=0}\frac{e^{As}(s+w)(s+z)}{s^2}\right)^2
 =\frac{(w+z+Awz)^2}{w+z}.
\]
Integration from $0$ to $\alpha_k$ gives
\begin{align}\label{eq:moll-kernel-residue}
\Res_{s_2=0}\Res_{s_1=0}
\frac{e^{\alpha_k(s_1+s_2)}
(s_1+w)(s_1+z)(s_2+w)(s_2+z)}
{s_1^2s_2^2(s_1+s_2)(w+z)}
=1+\alpha_k(w+z)+\alpha_k^2wz
+\frac{\alpha_k^3w^2z^2}{3(w+z)}.
\end{align}

Take the ordered residues on $|s_1|=1/(16\log X)$ and
$|s_2|=1/(8\log X)$. On these circles, $G(s_i)\Gamma_k(s_i)\ll_a1$,
the three zeta factors depending only on $s_1,s_2$ are each $O(\log X)$,
and $L(1+s_1+s_2,\sym^2g)$ is bounded by holomorphy at $1$.
For $1/\log X\leq\Re(w)\leq1$ and $2/\log X\leq\Re(z)\leq1$,
absolute convergence bounds the other global factors by powers
of $\log X$. The regular Euler product is bounded there, and
$\mathcal H_P(v)\ll_P(\log X)^{-2}|v|^{-3}$ for $|v|\log X\geq1$.
Thus we may move the Perron lines in
\eqref{eq:moll-second-perron-residues} to $\Re(w)=1/\log X$
and $\Re(z)=2/\log X$.

On these lines, Remark~\ref{rem:zero-free-region} and the Laurent
expansion at $1$ give
\begin{align*}
\left|\zeta(1+s_i+w)^{-1}\right|
&\ll \frac{|w|}{1+|w|}\bigl(\log(3+|w|)\bigr)^C,\qquad i=1,2,
\end{align*}
and
\begin{align*}
|\zeta(1+w+z)|
&\ll\frac{1+|w+z|}{|w+z|}
\bigl(\log(3+|w|+|z|)\bigr)^C.
\end{align*}
The first bound also holds with $w$ replaced by $z$.
The symmetric square quotient is bounded by a power of
$\log(3+|w|+|z|)$. Since
$(1+|w+z|)(\log(3+|w|+|z|))^C
\ll_C(1+|w|)^2(1+|z|)^2$, the four reciprocal zeta factors
in \eqref{eq:moll-Euler-singular} give
\[
\left|
\Res_{s_2=0}\Res_{s_1=0}
\frac{G(s_1)\Gamma_k(s_1)G(s_2)\Gamma_k(s_2)}{s_1s_2}
F_{2,t}(s_1,s_2,w,z)
\right|
\ll_t (\log X)^3\frac{|w|^2|z|^2}{|w+z|}.
\]
The ordered residue in \eqref{eq:moll-kernel-residue} satisfies
the same bound by using the same $s_i$ circles. Write
$w=(1+\mathrm{i}x)/\log X$ and $z=(2+\mathrm{i}y)/\log X$.
With the Perron weights and line elements, both kernels are
bounded by a constant times
\[
\frac{\dd x\,\dd y}{(1+|x|)(1+|y|)(1+|x+y|)}.
\]
This is integrable, since
\[
\int_{\mathbb R}\frac{\dd y}{(1+|y|)(1+|x+y|)}
\ll\frac{\log(2+|x|)}{1+|x|}.
\]
On each bounded set of $x,y$, the difference between the complete
kernel, including its factor $4$, and the rational kernel multiplied
by $\mathbf 1_{t=1}$ is controlled by
$N_{k,t}(s_1,s_2,w,z)-\mathbf 1_{t=1}=O_t(|s_1|+|s_2|+|w|+|z|)$
on the residue circles. 
Since the remaining factor in the integrand is \(O((\log X)^2)\) 
and each circle has length \(O(1/\log X)\), 
the difference between the two kernels tends to zero uniformly for \(k\asymp K\).
Dominated convergence gives an error $o_t(1)$ on replacing the
complete kernel by $\mathbf 1_{t=1}$ times
\eqref{eq:moll-kernel-residue} in \eqref{eq:moll-second-perron-residues}.

In the Perron integral, the terms $ \mathbf 1_{t=1}(1+\alpha_kw)(1+\alpha_kz)$
give, by \eqref{eq:moll-first-residue-P} and \eqref{eq:moll-first-residue-dP},
\begin{align*}
\mathbf 1_{t=1}\left(P(1)+\frac{\alpha_k}{\log X}P'(1)\right)^2.
\end{align*}
For the last term of \eqref{eq:moll-kernel-residue}, use
\[
\frac1{w+z}=\int_1^\infty v^{-w-z}\frac{\dd v}{v}.
\]
On either positive line, termwise inverse Mellin transformation
of $X^ww^2\mathcal H_P(w)$ gives, for almost every $v>0$,
\[
\frac1{(\log X)^2}
P''\left(\frac{\log(X/v)}{\log X}\right)\boldsymbol 1_{v<X}.
\]
Indeed, integrating each monomial on $0<v<X$ recovers its
Mellin transform, and $X^ww^2\mathcal H_P(w)$ is square
integrable in $\Im(w)$ on each line. Mellin Parseval
gives convergence of the truncated inversions in
$L^2(v^{2\Re(w)}\dd v/v)$, hence in $L^2(\dd v/v)$ for $v\geq1$.
Insert the formula for $(w+z)^{-1}$ on finite segments and
apply Fubini. The products of the truncated inversions converge
in $L^1$ for $v\geq1$ by Cauchy's inequality, and the convolution
bound above gives absolute convergence on the Mellin side.
The remaining contribution is therefore
\begin{align*}
\mathbf 1_{t=1}\frac{\alpha_k^3}{3(\log X)^4}
\int_1^X P''\left(\frac{\log(X/v)}{\log X}\right)^2
\frac{\dd v}{v}=\mathbf 1_{t=1}\frac{\alpha_k^3}{3(\log X)^3}
\int_0^1P''(u)^2\dd u.
\end{align*}
Combining the two contributions with the replacement error gives
\begin{equation}\label{eq:moll-second-twisted-main}
S_{2,t}(k)=\mathbf 1_{t=1}
\left\{\left(P(1)+\frac{\alpha_k}{\log X}P'(1)\right)^2
+\frac{\alpha_k^3}{3(\log X)^3}\int_0^1P''(u)^2\dd u\right\}
+o_t(1),
\end{equation}
uniformly for $k\asymp K$. On the support of $h_{K,L}$,
$\alpha_k/\log X=a^{-1}+O((\log K)^{-1})$, so
$S_{2,t}(k)=\mathbf 1_{t=1}\mu_2(P,a)+o_t(1)$ uniformly in $k$.
Poisson summation gives
\begin{align*}
\sum_{k\equiv0\bmod2}h_{K,L}(k)
&=\frac L2\int_{\mathbb R}h(u)\dd u+O(K^{-A}),
\end{align*}
and
\begin{align*}
\sum_{k\equiv0\bmod2}|h_{K,L}(k)|&\ll_h L.
\end{align*}
Averaging the estimate 
in \eqref{eq:moll-second-twisted-main} and adding the $ o_t(L)$ error in
\eqref{eq:general-mollified-error-sum} proves \eqref{eq:mollified-second-moment}.
For uniformity on the stated compact sets, at each $(\delta,a)$
choose $1/3<\delta_0<\delta$ and $a_0>a$ so that all three exponents
in \eqref{eq:general-mollified-error-sum} are negative at
$(\delta_0,a_0)$. These strict bounds hold in a neighbourhood of
$(\delta,a)$ by using \eqref{v6:eq:joint-error} with $\delta_0$.
A finite cover gives a common choice of $\eps$. The main-term
estimates are uniform since $a$ stays bounded away from zero.
\end{proof}

\subsection{Proof of the non-vanishing theorem}
We maximize $\mu_1(P,a)^2/\mu_2(P,a)$ over non-zero real
polynomials $P$ satisfying $P(0)=P'(0)=0$.
Integration gives
\begin{align*}
P'(1)=\int_0^1P''(u)\dd u,\qquad\text{and}\qquad
P(1)=\int_0^1(1-u)P''(u)\dd u.
\end{align*}
Hence
\begin{align*}
\mu_1(P,a)= & ~\int_0^1(1+a^{-1}-u)P''(u)\dd u.
\end{align*}
The conditions at zero imply that $P''\not\equiv0$.
Cauchy's inequality therefore gives
\begin{align}\label{eq:moll-optimal-functional}
\frac{\mu_1(P,a)^2}{\int_0^1P''(u)^2\dd u}
\leq & ~\int_0^1(1+a^{-1}-u)^2\dd u
=\frac{a^2+3a+3}{3a^2}.
\end{align}
Equality holds precisely when $P''(u)$ is proportional to
$1+a^{-1}-u$. For fixed $a$, corresponding to mollifier length
$X=K^{2a}$, let $P_a$ be the polynomial attaining equality
with $P_a(1)=1$. The ratio $\mu_1(P,a)^2/\mu_2(P,a)$ is
invariant under rescaling $P$ by a non-zero constant.
Integrating twice with $P_a(0)=P_a'(0)=0$ gives
\begin{align*}
P_a(u)= & ~\frac{3(a+1)u^2-au^3}{3+2a}.
\end{align*}
By the definition of $\mu_2(P,a)$ and
\eqref{eq:moll-optimal-functional},
\begin{align*}
\frac{\mu_1(P,a)^2}{\mu_2(P,a)}
\leq & ~\frac{a(a^2+3a+3)}{1+a(a^2+3a+3)}
=1-(1+a)^{-3},
\end{align*}
with equality for $P=P_a$.

\begin{proof}[Proof of Theorem~\ref{thm:intro-nonvanishing}]
Fix $\delta,a$ as in \eqref{eq:intro-parameters} and take
$P=P_a$. This polynomial satisfies $P_a(0)=P_a'(0)=0$
and $P_a(1)=1$, so Proposition~\ref{prop:mollified-moments} applies
with $t=1$.
Also, $\mu_2(P_a,a)>0$ by its definition.
Lemma~\ref{Petersson} and Poisson summation give
\begin{align*}
\sum_{k\equiv0\bmod2}h_{K,L}(k)
\sideset{}{^h}\sum_{f\in H_k}1
= & ~\frac L2\int_{\mathbb R}h(u)\dd u+O(K^{-A}).
\end{align*}
Since $h$ is non-negative and non-zero, its integral is positive.
The second moment is therefore positive for sufficiently large $K$.

The central values and $M_{P_a}(f)$ are real.
Cauchy's inequality for the non-negative harmonic weights,
the two moment formulas, and the total mass formula give
\begin{align*}
& ~\sum_{k\equiv0\bmod2}h_{K,L}(k)
\sideset{}{^h}\sum_{\substack{f\in H_k\\L(1/2,f\otimes g)\ne0}}1
\notag\\
\geq & ~\frac{\left(\displaystyle
\sum_{k\equiv0\bmod2}h_{K,L}(k)
\sideset{}{^h}\sum_{f\in H_k}
L(1/2,f\otimes g)M_{P_a}(f)\right)^2}
{\displaystyle\sum_{k\equiv0\bmod2}h_{K,L}(k)
\sideset{}{^h}\sum_{f\in H_k}
L(1/2,f\otimes g)^2M_{P_a}(f)^2}
\notag\\
= & ~\bigl(1-(1+a)^{-3}+o(1)\bigr)
\sum_{k\equiv0\bmod2}h_{K,L}(k)
\sideset{}{^h}\sum_{f\in H_k}1.
\end{align*}
The errors are uniform for $K^\delta\leq L\leq K$, with
$\delta,a,g,h$ fixed. This proves \eqref{eq:intro-nonvanishing}.
Finally, $1-(1+a)^{-3}$ is continuous and strictly increasing
for $a>0$. Every constant strictly below
$1-(1+\mathfrak a(\delta))^{-3}$ is therefore obtained by
a fixed choice $0<a<\mathfrak a(\delta)$.
\end{proof}

\begin{lemma}
For $Y\geq1$, let $B_Y$ be the least common multiple
of $1,\ldots,\lfloor Y\rfloor$ and set
\[
 Q_Y(f)=\zeta(2)\sum_{n\mid B_Y}\frac{\lambda_f(n^2)}n.
\]
Then $Q_Y(f)\geq0$. For every fixed $Y$,
\begin{equation}\label{eq:sym2-L2-approximation}
\limsup_{\substack{k\to\infty\\2\mid k}}
\sideset{}{^h}\sum_{f\in H_k}|L(1,\sym^2 f)-Q_Y(f)|^2
\ll Y^{-1}.
\end{equation}
\end{lemma}
\begin{proof}
By multiplicativity, $Q_Y(f)/\zeta(2)$ is the product over
$p^J\parallel B_Y$ of the sums on the left below. By
\eqref{Hecke relation}, $\sum_{j=0}^i\lambda_f(p^{2j})=\lambda_f(p^i)^2$.
Partial summation gives
\[
 \sum_{j=0}^J\frac{\lambda_f(p^{2j})}{p^j}
 =p^{-J}\lambda_f(p^J)^2+
 (1-p^{-1})\sum_{j=0}^{J-1}p^{-j}\lambda_f(p^j)^2\geq0,
\]
which proves the assertion about the sign.

For fixed $m,n$, Lemma~\ref{Petersson} gives
\[
\lim_{\substack{k\to\infty\\2\mid k}}
\sideset{}{^h}\sum_{f\in H_k}\lambda_f(m^2)\lambda_f(n^2)
=\mathbf 1_{m=n}.
\]
Expanding $Q_Y(f)^2$, we obtain
\[
\lim_{\substack{k\to\infty\\2\mid k}}
\sideset{}{^h}\sum_{f\in H_k}Q_Y(f)^2
=\zeta(2)^2\sum_{n\mid B_Y}\frac{1}{n^2}.
\]

For the mixed term, the eigenvalue of $T_{n^2}$ on $f$ is
$n^{k-1}\lambda_f(n^2)$. We obtain
$\sum_{f\in H_k}\lambda_f(n^2)=(k-1)/(12n)+O_n(1)$
by dividing the level-one trace formula of
\cite[Proposition~11]{MurtySinha2009} by $n^{k-1}$.
The factor $L(1,\sym^2 f)$ cancels the denominator in the harmonic weight, so
\[
 \sideset{}{^h}\sum_{f\in H_k}L(1,\sym^2 f)\lambda_f(n^2)
 =\frac{2\pi^2}{k-1}\sum_{f\in H_k}\lambda_f(n^2)
 =\frac{\zeta(2)}{n}+O_n(k^{-1}).
\]
Multiplying by $\zeta(2)/n$ and summing over $n\mid B_Y$, we obtain
\[
 \sideset{}{^h}\sum_{f\in H_k}L(1,\sym^2 f)Q_Y(f)
 =\zeta(2)^2\sum_{n\mid B_Y}\frac{1}{n^2}+O_Y(k^{-1}).
\]

It remains to evaluate the harmonic second moment of $L(1,\sym^2 f)$.
Apply \cite[Proposition~4.1 and Eq.\,(4.8)]{Xiao2016} with
$N=1$, $m=2$, $z=2$ and $\eta=1/130$. On the zero-free subfamily defined in \cite[Eq.\,(1.35)]{Xiao2016},
this moment equals, up to $o(1)$, the product over $p$ of
the Sato--Tate averages of $L_p(1,\sym^2 f)^2$.
For local parameters $e^{\pm\mathrm i\theta}$, write
$\chi_j(\theta)=\sin((j+1)\theta)/\sin\theta$.
The local identity in Lemma~\ref{Euler product} gives
\[
 L_p(1,\sym^2 f)=\frac{1}{1-p^{-2}}
 \sum_{j\geq0}\frac{\chi_{2j}(\theta)}{p^j}.
\]
This series converges uniformly in $\theta$ since
$|\chi_{2j}(\theta)|\leq 2j+1$. The characters $\chi_j$ are
orthonormal with respect to
$(2/\pi)\sin^2\theta\,\dd\theta$ on $[0,\pi]$.
Squaring and integrating term by term, we find that the local factor
at $p$ is
\[
 \frac{1}{(1-p^{-2})^2}\sum_{j\geq0}p^{-2j}=(1-p^{-2})^{-3}.
\]
Thus the product is $\zeta(2)^3$.

By \cite[Eq.~(1.36)]{Xiao2016}, there are
$O(k^{65\eta})=O(k^{1/2})$ forms in its complement. Each contributes
$2\pi^2L(1,\sym^2 f)/(k-1)$ to the harmonic second moment. Since
$L(1,\sym^2 f)\ll(\log(k+3))^3$ by \cite[Lemma~2.2]{Xiao2016}, their total
contribution is $O(k^{-1/2}(\log(k+3))^3)=o(1)$.
Hence
\[
\lim_{\substack{k\to\infty\\2\mid k}}
\sideset{}{^h}\sum_{f\in H_k}L(1,\sym^2 f)^2=\zeta(2)^3.
\]
Expanding the square and combining the three estimates then gives
\[
 \lim_{\substack{k\to\infty\\2\mid k}}
 \sideset{}{^h}\sum_{f\in H_k}|L(1,\sym^2 f)-Q_Y(f)|^2
 =\zeta(2)^3-\zeta(2)^2\sum_{n\mid B_Y}\frac{1}{n^2}
 =\zeta(2)^2\sum_{n\nmid B_Y}\frac{1}{n^2}.
\]
Since every positive integer $n\leq Y$ divides $B_Y$, the tail is at most $\sum_{n>Y}n^{-2}\ll Y^{-1}$.
\end{proof}

\begin{proof}[Proof of Corollary~\ref{cor:ordinary-nonvanishing}]
We take $P=P_a$ and write $\mu_j=\mu_j(P_a,a)$ for $j=1,2$. By the
equality case of \eqref{eq:moll-optimal-functional} and the definition
of $\mu_2$, we have $\mu_1^2/\mu_2=1-(1+a)^{-3}$.
We put $V_f=L(1/2,f\otimes g)M_{P_a}(f)$ and
$I_f=\mathbf 1_{\{L(1/2,f\otimes g)\ne0\}}$. For a fixed non-zero
$h\in C_c^\infty((0,1))$ with $h\geq0$, we use the original
harmonic sums below.

Fix $Y\geq1$. For each $n\mid B_Y$, apply
\eqref{eq:mollified-first-moment} and \eqref{eq:mollified-second-moment}
with $t=n^2$. Only $n=1$ contributes to the main term. Likewise,
Lemma~\ref{Petersson} gives
\[
 \begin{aligned}
 &\sum_{k\equiv0\bmod2}h_{K,L}(k)
 \sideset{}{^h}\sum_{f\in H_k}Q_Y(f)=\frac L2\zeta(2)\int_{\mathbb R}h(u)\,\dd u+o_Y(L).
 \end{aligned}
\]
Summing the fixed-twist formulas over $n\mid B_Y$ with coefficients
$\zeta(2)/n$ therefore gives
\begin{equation}\label{eq:ordinary-positive-moments}
 \begin{aligned}
 \sum_{k\equiv0\bmod2}h_{K,L}(k)
 \sideset{}{^h}\sum_{f\in H_k}Q_Y(f)V_f^j=\mu_j\sum_{k\equiv0\bmod2}h_{K,L}(k)
 \sideset{}{^h}\sum_{f\in H_k}Q_Y(f)+o_Y(L),
 \qquad j=1,2.
 \end{aligned}
\end{equation}

Since $\mu_2>0$, the second moment in \eqref{eq:ordinary-positive-moments}
is positive for sufficiently large $K$. Using $V_f=I_fV_f$ and
Cauchy's inequality with non-negative weight $Q_Y$, we obtain
\[
 \begin{aligned}
 \sum_{k\equiv0\bmod2}h_{K,L}(k)
 \sideset{}{^h}\sum_{f\in H_k}Q_Y(f)I_f
 &\geq\frac{\left(\displaystyle\sum_{k\equiv0\bmod2}h_{K,L}(k)
 \sideset{}{^h}\sum_{f\in H_k}Q_Y(f)V_f\right)^2}
 {\displaystyle\sum_{k\equiv0\bmod2}h_{K,L}(k)
 \sideset{}{^h}\sum_{f\in H_k}Q_Y(f)V_f^2}\\
 &=\bigl(1-(1+a)^{-3}\bigr)\sum_{k\equiv0\bmod2}h_{K,L}(k)
 \sideset{}{^h}\sum_{f\in H_k}Q_Y(f)+o_Y(L).
 \end{aligned}
\]
By \eqref{eq:sym2-L2-approximation} and Cauchy's inequality,
\[
 \begin{aligned}
 &\sum_{k\equiv0\bmod2}h_{K,L}(k)
 \sideset{}{^h}\sum_{f\in H_k}|L(1,\sym^2 f)-Q_Y(f)|\\
 &\quad\leq\left(\sum_{k\equiv0\bmod2}h_{K,L}(k)
 \sideset{}{^h}\sum_{f\in H_k}1\right)^{1/2}
 \left(\sum_{k\equiv0\bmod2}h_{K,L}(k)
 \sideset{}{^h}\sum_{f\in H_k}|L(1,\sym^2 f)-Q_Y(f)|^2\right)^{1/2}\\
 &\quad\ll_h LY^{-1/2}+o_Y(L).
 \end{aligned}
\]
For fixed $Y$ the error is uniform in $L$, since $k\asymp K$ on the
support of $h_{K,L}$. As $|I_f-1+(1+a)^{-3}|\leq1$, this bounds the error in
replacing $Q_Y(f)$ by $L(1,\sym^2 f)$ in
\[
 \sum_{k\equiv0\bmod2}h_{K,L}(k)
 \sideset{}{^h}\sum_{f\in H_k}
 Q_Y(f)(I_f-1+(1+a)^{-3}).
\]
Taking the limits as $K$ and then $Y$ tend to infinity, we obtain
\begin{equation}\label{eq:ordinary-sym2-count}
 \begin{aligned}
 &\sum_{k\equiv0\bmod2}h_{K,L}(k)
 \sideset{}{^h}\sum_{f\in H_k}L(1,\sym^2 f)I_f\\
 &\qquad\geq \bigl(1-(1+a)^{-3}\bigr)\sum_{k\equiv0\bmod2}h_{K,L}(k)
 \sideset{}{^h}\sum_{f\in H_k}L(1,\sym^2 f)-o(L).
 \end{aligned}
\end{equation}
uniformly for $K^\delta\leq L\leq K$.

We fix $0<\eta<1/2$ and choose $0\leq h\leq1$ as above with $h=1$
on $[\eta,1-\eta]$. We apply \eqref{eq:ordinary-sym2-count} to $h(u)$
and $u h(u)$, multiply by $K$ and $L$, respectively, and add.
For $u=(k-1-K)/L$, we have $K+Lu=k-1$, and the harmonic weight
times $(k-1)L(1,\sym^2 f)$ is $2\pi^2$. Dividing by $2\pi^2$, we obtain
\[
 \sum_{k\equiv0\bmod2}h_{K,L}(k)\sum_{f\in H_k}I_f
 \geq \bigl(1-(1+a)^{-3}\bigr)\sum_{k\equiv0\bmod2}h_{K,L}(k)\#H_k-o_\eta(KL).
\]
For fixed $\eta$ and sufficiently large $K$, the support lies in
$K<k\leq K+L$. Using the dimension formula $\#H_k=(k-1)/12+O(1)$,
we obtain
\[
 \sum_{\substack{K<k\leq K+L\\k\equiv0\bmod2}}\#H_k
 =\frac{KL}{24}+\frac{L^2}{48}+O(K)\asymp KL.
\]
The two end intervals contain $O(\eta L+1)$ weights, each with
$O(K)$ forms. Thus replacing the smooth total count by the full count
costs $O(\eta KL+K)$, while the full non-vanishing count dominates
the smooth one. After division by the full count, we let
$K$ tend to infinity and then let $\eta$ tend to zero. This proves the
stated bound for fixed $a$. The last assertion follows by taking $a$
arbitrarily close to $\mathfrak a(\delta)$ from below.
\end{proof}

\section{Future directions}
A natural question is whether a positive harmonic proportion of
non-vanishing holds in intervals of length $L=K^{1/3}$. One approach
would be to establish joint square-root cancellation in $n$ and
$\Delta$ in \eqref{eq:ndnz-shift-representation}, uniformly over
the twists in the mollified second moment.

For a result for almost all centres, one could instead seek a
power-saving mean-square bound over $K\leq x\leq2K$ for
$\sqrt r\,\mathcal E(r;x,x^{1/3})/x^{1/3}$, where
$\mathcal E(r;x,x^{1/3})$ is the non-zero-shift error for
centre $x$ and interval length $x^{1/3}$. A bound uniform over a
sufficiently long range of twists could yield a positive harmonic
proportion outside a set of centres of measure $o(K)$.

\bigskip
\noindent{\bf AI assistance.} 
The authors used ChatGPT (OpenAI) to assist with literature searches, checking calculations, 
and polishing the language of the manuscript. 
In particular, ChatGPT assisted with the detailed Euler product computations 
in the proofs of \eqref{eq:mollified-first-moment} and \eqref{eq:mollified-second-moment}. 
The authors take full responsibility for the mathematical arguments, results, and their presentation.

\bigskip

\bibliographystyle{amsplain}

\end{document}